\documentclass[preprint,sort&compress,12pt]{elsarticle}
\usepackage{mathrsfs}
\usepackage{mathrsfs}
\usepackage{amsmath}
\usepackage{mathrsfs}
\usepackage{stmaryrd}
\usepackage{mathrsfs}
\usepackage{bm}
\usepackage{amsmath}
\usepackage{amsfonts}
\usepackage{amssymb}
\usepackage{amsthm}
\usepackage[RGB]{xcolor}
\journal{\quad}
\newcommand{\mf}{\mathbf}
\newcommand{\mm}{\mathrm}

\begin{document}

\begin{frontmatter}
\title{Global Weak Solutions to the Equations of Compressible Flow of
Nematic Liquid Crystals in Two Dimensions\tnoteref{S}}
% \tnotetext[S]{The research of Fei Jiang was supported by the NSFC (Grant No. 11101044),
% the research of Song Jiang by NSFC (Grant No. 40890154) and the National Basic
% Research Program under the Grant 2011CB309705, and the research of
% Dehua Wang by the National Science Foundation under Grant DMS-0906160 and the Office of Naval
% Research under Grant N00014-07-1-0668.}
%%%%%
\author[FJ,SJ]{Fei Jiang}  %\corref{cor1}}
\ead{jiangfei0591@163.com}
\author[SJ]{Song Jiang}
\ead{jiang@iapcm.ac.cn}
% \cortext[cor1]{Corresponding
%author: Tel +86 15001201710.}
\author[DH]{Dehua Wang}
\ead{dwang@math.pitt.edu}
\address[FJ]{College of Mathematics and Computer Science, Fuzhou University, Fuzhou, 350108, China.}
\address[SJ]{Institute of Applied Physics and Computational Mathematics, Beijing, 100088, China.}
\address[DH]{Department of Mathematics, University of Pittsburgh, Pittsburgh, PA, 15260, USA.}

\begin{abstract}
We consider weak solutions to a two-dimensional simplified
Ericksen-Leslie system of compressible flow of nematic liquid crystals.
An initial-boundary value problem is first studied in a bounded domain.
By developing new techniques and estimates to overcome the difficulties induced by the
supercritical nonlinearity $|\nabla{\mathbf d}|^2{\mathbf d}$ in the
equations of angular momentum on the direction field, and adapting the standard
three-level approximation scheme and the weak convergence arguments for the compressible
Navier-Stokes equations, we establish the global existence of weak solutions
under a restriction imposed on the initial energy including the case of small initial energy.
%% The uniform a priori estimates  are obtained and  the weak convergence arguments are adopted.
Then the Cauchy problem with large initial data is investigated, and we prove the global existence
of large weak solutions by using the domain expansion
technique and the rigidity theorem, provided that the second component
of initial data of the direction field satisfies some geometric angle condition.

\iffalse
We study an two-dimensional initial-boundary value problem for the
simplified Ericksen-Leslie system modeling compressible nematic liquid crystal flows.
% in a bounded smooth domain $\Omega\subset \mathbb{R}^2$%  or $\mathbb{R}^2$
We first show the unique solvability of the third approximate problem under certain
restrictions imposed on the modified initial energy
% or on the second component of initial data of $\mf{d}$
 by overcoming the difficulties induced by the the supercritical nonlinearity
$|\nabla\mf{d}|^2\mf{d}$ in the equations of angular momentum on the
macroscopic average of the nematic liquid crystal orientation field
$\mf{d}$. Then, we adapt a standard three-level approximation scheme
to obtain the global existence of weak solutions to the original
problem under certain restrictions imposed on the initial energy,
including the case of small initial energy, by deriving the uniform
a priori estimates and employing weak convergence arguments that
have demonstrated their success for the compressible Navier-Stokes
equations. In addition, using   domain expansion technique and
rigidity theorem on $\mathbb{R}^2$, we also establish the global
existence of large weak solutions to the Cauchy problem provided the
second component of initial data of $\mf{d}$ satisfies some
geometric angle condition.
\fi
\end{abstract}

\begin{keyword} Liquid crystals, compressible flows, weak solutions, Galerkin method,
weak convergence arguments.
\MSC[2000] 35Q35\sep 76D03.
%(2000 is the default)
\end{keyword}
\end{frontmatter}

%%
%% Start line numbering here if you want
%%
% \linenumbers

%% main text
\newtheorem{thm}{Theorem}[section]
\newtheorem{lem}{Lemma}[section]
\newtheorem{pro}{Proposition}[section]
\newtheorem{cor}{Corollary}[section]
\newproof{pf}{Proof}
\newdefinition{rem}{Remark}[section]
\newtheorem{definition}{Definition}[section]
\newcommand{\red}{\color{red}}

\section{Introduction}\label{sec:01}
\label{Intro} \numberwithin{equation}{section}

We consider the existence of weak solutions to an initial-boundary value problem
for the following two-dimensional simplified version of the Ericksen-Leslie model
in a bounded domain $\Omega\subset{\mathbb R}^2$ which describes
the motion of a compressible flow of nematic liquid crystals:
\begin{eqnarray}
&& \label{0101}\partial_t\rho+\mathrm{div}(\rho\,\mathbf{v})=0,\\
 &&  \label{0102} \partial_t(\rho
\mathbf{v})+\mathrm{div}(\rho\mathbf{v}\otimes\mathbf{v})+ \nabla
P(\rho)=\mu \Delta\mathbf{v}+(\mu+\lambda)\nabla
\mm{div}\mf{v}\nonumber\\
&&\qquad \qquad \qquad \qquad  \qquad \qquad  \quad \qquad
-\nu\mm{div}\Big( \nabla \mathbf{d}\odot\nabla
\mathbf{d}-\frac{1}{2}|\nabla \mf{d}|^2\mathbb{I}\Big),  \\
&&  \label{0103}\partial_t\mathbf{d}+\mathbf{v}\cdot \nabla
\mathbf{d}= \theta(\Delta \mathbf{d}+|\nabla \mathbf{d}|^2\mf{d}),
\end{eqnarray}
where
% the spatial variable is $\mathbf{x}=(x_1,x_2)\in{\mathbb R}^2$, $t\ge 0$ is the time variable;
%  on a time interval $[0,T)$,
$\rho$ is the density of the nematic liquid crystals, $\mathbf{v}$ the velocity and
$P(\rho)$ the pressure, $\mathbf{d}\in\mathbb{S}^1:=\{\mf{d}\in
\mathbb{R}^2~|~|\mf{d}|=1\}$ represents the macroscopic average of
the nematic liquid crystal orientation field. The constants $\mu$, $\lambda$, $\nu$, and
$\theta$  denote the shear viscosity, the bulk viscosity, the
competition between kinetic energy and potential energy, and the
microscopic elastic relation time for the molecular orientation
field, respectively. $\mathbb{I}$ denotes the $2\times 2$ identical
matrix. The term $\nabla \mathbf{d}\odot\nabla \mathbf{d}$ denotes the $2\times 2$ matrix
whose $(i,j)$-th entry is given by
$\partial_{x_i}\mathbf{d}\cdot\partial_{x_j}\mathbf{d}$, for $1\leq i,j\leq 2$, i.e.,
 $\nabla \mathbf{d} \odot\nabla \mathbf{d}= (\nabla\mathbf{d})^{\top}\nabla\mathbf{d},$
where $(\nabla\mathbf{d})^{\top}$ denotes the transpose of the $2\times2$ matrix $\nabla\mathbf{d}$.

% \emph{For simplicity, when it is not confusing, the
%notation $d_0|_{\partial\Omega}=\gamma_{\Omega}(d_0)$ is used not
%only for $d_0\in C^{\infty}(\bar{\Omega})$ but also for $d_0\in
%H^k(\Omega)$}, where the function $\gamma_{\Omega}(d_0)(x)$ is
%called the trace of the function $d_0(x)$ on the boundary $\partial
%\Omega$ and $k$ is a positive integer (see [26]).

\iffalse
Nematic liquid crystals are aggregates of molecules which possess
the same orientational order and are made of elongated, rod-like
molecules. The continuum theory of liquid crystals was developed by
Ericksen \cite{EJLHA9} and Leslie \cite{LFMSA2} during the period of
1958 through 1968, see also the book by de Gennes \cite{GPGT1}.
Since then there have been remarkable research developments in
liquid crystals from both theoretical and applied aspects. When a
fluid containing nematic liquid crystal materials is at rest, we
have the well-known Ossen-Frank theory for static nematic liquid
crystals. The reader can refer to the
%%% pioneering
work by Hardt-Kinderlehrer-Lin \cite{HRKDLFEC} on the analysis
of energy minimal configurations
of nematic liquid crystals. In general, fluid motion always takes
place. The so-called Ericksen-Leslie system is a macroscopic
continuum description of the time evolution of materials under the
influence of both flow velocity field $\mf{v}$ and
macroscopic description of the microscopic orientation
configurations $\mf{d}$ of rod-like liquid crystals.
\fi

In 1989, Lin \cite{LFHNC} first derived a simplified
Ericksen-Leslie system modeling liquid crystal flows when the
fluid is  incompressible and viscous. Subsequently, Lin and Liu
\cite{LFHLCNC5,LFHLCPD2} established some analysis results on
the simplified Ericksen-Leslie system, such as the existence of
weak and strong solutions and the partial regularity of suitable
solutions, under the assumption that the liquid crystal  director
field is of varying length by Leslie's terminology, or variable
degree of orientation by Ericksen's terminology. We refer readers to \cite{EJLHA9,LFMSA2,GPGT1,HRKDLFEC}
for details concerning the so-called Ericksen-Leslie system.

When the fluid is allowed to be compressible, the Ericksen-Leslie
system becomes more complicated. From the viewpoint of
partial differential equations, the system
(\ref{0101})--(\ref{0103}) is a highly nonlinear hyperbolic-parabolic coupled system,
and it is challenging to analyze such a system, in particular, as the density function
$\rho$ may vanish. It should be noted that the system
(\ref{0101})--(\ref{0103}) includes two important subsystems. When
$\rho$ is constant and $\mf{v}=\mathbf{0}$, the system
(\ref{0101})--(\ref{0103}) reduces to the equations for heat flow of
harmonic maps into $\mathbb{S}^1$, on which there have been extensive studies
in the past few decades (see, for example, the monograph by Lin and Wang \cite{LFHWCYTA}
and the references therein). When $\mf{d}$ is constant, the system
\eqref{0101}, \eqref{0102} becomes the compressible Naiver-Stokes
equations, which have attracted great interest in the past decades
and there have been many important developments (see, e.g.,
\cite{LPLMTFM98,f2,NASII04,HLJSNH} for a survey of recent developments).
%
% Of cause, since the system (\ref{0101})--(\ref{0103}) contains the
% Navier-Stokes equations as a subsystem, one cannot expect, in
% general, any better results than those for the Navier-Stokes equations.
%

Since the supercritical nonlinearity $|\nabla\mf{d}|^2\mf{d}$
causes significant mathematical difficulties, Lin in \cite{LFHNC} introduced a
Ginzburg-Landau approximation of the simplified Ericksen-Leslie
system, i.e., $|\nabla \mf{d}|^2\mf{d}$ in (\ref{0103}) is replaced
by the Ginzburg-Landau penalty function $(1-|\mf{d}|^2)/\epsilon$ or by a
more general penalty function. Consequently, by establishing some
estimates to deal with the direction field and its coupling/interaction
with the fluid variables, a number of results on the Navier-Stokes equations
can be successfully generalized to such Ginzburg-Landau approximation model.
For examples, when $\rho$ is a constant, i.e., the homogeneous
incompressible case, Lin and Liu \cite{LFHLCNC5} proved the global
existence of weak solutions in two and three dimensions. In
particular, they also obtained the existence and uniqueness of
global classical solutions either in two dimensions or in three dimensions
for large fluid viscosity $\mu$. In addition, the existence of weak
solutions to the density-dependent incompressible flow of liquid
crystals was proved in \cite{LXZZLC3121,JFTZOGWS}. Recently, Wang
and Yu \cite{WDHYCGA}, and Liu and Qin \cite{LXGQJE} independently
established the global existence of weak solutions to
the three-dimensional compressible flow of liquid crystals
with the Ginzburg-Landau penalty function.

In the past a few years, progress has also been made on the analysis of
%% there are also remarkable research developments for
the model  (\ref{0101})--(\ref{0103}) by
overcoming the difficulty induced by the supercritical nonlinearity
$|\nabla\mf{d}|^2\mf{d}$. For the incompressible case,
% (i.e. $\rho$ is constant in (\ref{0101}), (\ref{0102}), thus (\ref{0101})
% reduces to the incompressibility condition $\mm{div}\mf{v}=0$, and $P$ is an unknown pressure
% function in (\ref{0102})), Lin-Lin-Wang \cite{LFHLJYWCY} established
the existence of weak solutions in two dimensions was established in \cite{LFHLJYWCY},
%while  Wang \cite{WCYWA} and Lin \cite{LJYDSJO} proved
and the local and global existence of small strong solutions in three dimensions was proved in
\cite{HW1, LW-JDE2012,LJYDSJO,WCYWA}.
For the compressible case, the existence of strong
solutions have been investigated extensively. For examples,
%Huan et al. \cite{HTWCYWHY,HTWCYWWHYSS} established
the local existence of strong solutions
and a blow-up criterion were obtained in \cite{HTWCYWHY,HTWCYWWHYSS}, while
% Hu and Wu \cite{HXWHGS} showed
the existence and uniqueness of global strong solutions to the Cauchy problem in
critical Besov spaces were proved in \cite{HXWHGS} provided that the initial data are close to an
equilibrium state,
%Li, Xu and Zhang \cite{LJXZHZJWG} established
and the global existence of classical solutions to the Cauchy problem was shown in  \cite{LJXZHZJWG}
with smooth initial data that has small energy but possibly
large oscillations with possible vacuum and constant state as far-field condition.
% which could be either vacuum or non-vacuum.
To our best knowledge, however, there are no results available on
 weak solutions of the multi-dimensional problem
(\ref{0101})--(\ref{0103}) with large initial data due to the
difficulties induced by the compressibility and the supercritical
nonlinearity. It seems that the only global existence of large weak
solutions to (\ref{0101})--(\ref{0103}) was shown in the
one-dimensional case in \cite{DSLJWCWH}.

 In this paper, we will establish the global existence of weak solutions to the
 two-dimensional problem (\ref{0101})--(\ref{0103}) in a domain $\Omega\subset \mathbb{R}^2$ with initial conditions:
\begin{equation} \label{0105}
\rho(\mathbf{x},0)=\rho_0(\mathbf{x}),\ \
\mathbf{d}(\mathbf{x},0)=\mathbf{d}_0(\mathbf{x}),\ \ (\rho
\mathbf{v})(\mathbf{x},0)=\mathbf{m}_0(\mathbf{x})\quad\mbox{ in }\Omega ,
\end{equation}
and boundary conditions:
\begin{equation} \label{0104}
                              \left\{ \begin{array}{ll}
                                \mathbf{d}(\mathbf{x},t)=\mathbf{d}_0(\mathbf{x}), \quad
\mathbf{v}(\mathbf{x},t)=\mathbf{0},\quad \mf{x}\in\partial\Omega
,\; t>0,
           \quad \mbox{if }\Omega\mbox{ is a bounded domain};\\
                                 (\rho,\mathbf{v},\mathbf{d})(\mathbf{x},t)\rightarrow
(\rho_\infty,\mf{0},\mathbf{e}_2)\;\mbox{ as }|\mf{x}|\to\infty,\quad\mbox{if }\Omega=\mathbb{R}^2 ,
                               \end{array}\right.
             \end{equation}
where $\rho_\infty\geq 0$ is a constant and $\mf{e}_2=(0,1)\in\mathbb{R}^2$ is the unit vector.
For simplicity, we call (\ref{0101})--(\ref{0104}) the initial-boundary value problem
when $\Omega\in\mathbb{R}^2$ is a bounded domain, and
the Cauchy problem in the case of  $\Omega=\mathbb{R}^2$.

 Before stating our main result, we remark that the constants $\mu$, $\lambda$, $\nu$, and $\theta$ satisfy
the physical conditions:
\begin{equation*} \label{0107}\mu>0,\quad \lambda +\mu\geq 0,\quad \nu>0,\quad\theta>0.\end{equation*}
And the pressure $P(\rho)$ is usually determined through the
equations of state, here we focus our study on the case of
isentropic flows as in \cite{WDHYCGA} and assume that
\begin{equation*}
\label{0106}P(\rho)=A\rho^\gamma,\quad\mbox{ with }A>0,\ \gamma>1.
\end{equation*}
In addition, for the sake of simplicity, we define
$$I:=I_T:=(0,T),\quad Q_T=\Omega\times I,$$
\begin{equation}\label{dispassive}\mathcal{F}(t):=\mathcal{F}(\rho,\mf{v},\mf{d}):=\int_\Omega
\left(\mu |\nabla \mf{v}|^2+(\lambda+\mu)|\mm{div}\mf{v}|^2 +\theta(|\Delta
\mf{d}+|\nabla \mf{d}|^2\mf{d}|^2)\right)\mathrm{d}\mathbf{x},
\end{equation} and
\begin{equation}\label{jsw0108}
\mathcal{E}(t):= \mathcal{E}(\rho,\mf{m},\mf{d})=\int_\Omega\left(\frac{1}{2}
\frac{|\mf{m}|^2}{\rho}1_{\{\rho>0\}}+ Q(\rho)+\frac{\nu\theta|\nabla
\mf{d}|^2}{2}\right)\mm{d}\mf{x}\;\;\mbox{ with }\;\;\mf{m}=\rho\mf{v} ,
\end{equation}
where $1_{\{\rho>0\}}$ denotes the character function, and
\begin{equation}\label{fzu0114}
Q(s)= \frac{A}{\gamma-1}\left(s^\gamma-\gamma s\rho_\infty^{\gamma-1}
+(\gamma-1)\rho_\infty^\gamma\right) \;\mbox{ for some constant }\rho_\infty\geq 0.
\end{equation}
Obviously, $Q(s)\geq 0$ for any $s\geq 0$, and $\mathcal{E}(t)\geq 0$.
Throughout this paper, we use the bold fonts to denote the product spaces, for examples,
\begin{equation*}  \begin{aligned}&
\mathbf{L}^p(\Omega):=(L^p(\Omega))^2,\quad
\mathbf{H}^k_0(\Omega):=({H}^k_0(\Omega))^2=(W^{k,2}_0(\Omega))^2,\quad
%\mathbf{L}^p(\Omega):=(L^p(\Omega))^2, \\ %\mbox{ for }p=2,\ {2\gamma}/({\gamma+1}),\\
%\mf{L}^\frac{2\gamma}{\gamma+1}_{\mm{weak}}(\Omega):
%=({L}^\frac{2\gamma}{\gamma+1}_{\mm{weak}}(\Omega))^2,\quad
\mathbf{H}^k(\Omega):=(W^{k,2}(\Omega))^2;\end{aligned}\end{equation*}
and the Sobolev space with weak topology is defined as
 $$ C^0(\bar{I},\mf{L}^q_{\mm{weak}}(\Omega)):= \left\{\mf{f}: I\rightarrow
\mf{L}^q(\Omega)~\bigg|~\int_\Omega \mf{f}\cdot\mf{g}\mm{d}\mf{x}\in
C(\bar{I})\mbox{ for any }\mf{g}\in \mf{L}^\frac{q}{q-1}(\Omega)\right\}.$$

Our first result on the existence of weak solutions in a bounded
domain $\Omega$ reads as follows.
%%%%%%%%%%%%%%%%%%%%%%%%%%%%%%%%%%%%%%%%
\begin{thm}\label{thm:0101}  Let the constant $\rho_\infty\geq 0$,
$\Omega$ be a bounded domain of class $C^{2,\alpha}$ with $\alpha\in(0,1)$,
and the initial data $\rho_0,\,\mathbf{m}_0,\,\mathbf{d}_0$
 satisfy the following conditions:
\begin{align}   &\label{jfw0110}
  Q(\rho_0)\in L^1(\Omega),\quad \rho_0\geq0\;\mbox{ a.e. in }  \Omega ,   \\
  &\label{fzu0109} \mathbf{m}_0\in \mf{L}^{\frac{2\gamma}{\gamma+1}}(\Omega),\quad
  \mf{m}_01_{\{\rho_0=0\}}=\mf{0}\;\mbox{ a.e. in }\Omega, \quad \frac{|\mathbf{m}_0|^2}{\rho_0}1_{\{\rho_0>0\}}
  \in L^1(\Omega),  \\
  &\label{jfw0112}  |\mf{d}_0|=1\;\mbox{ in }\Omega,\quad
   \mathbf{d}_0(\mathbf{x})\in \mathbf{H}^2({\Omega}).
                                \end{align}
 %and $$\mathbf{d}_0(\mathbf{x})\in \mathbf{H}^1({\Omega})\mbox{ with
%}\mathbf{d}_0|_{\partial\Omega}(\mathbf{x})\in
%\mathbf{C}^{2,\delta}(\partial \Omega).$$
 Then, there exists a constant $C_0$, such that if $\mathcal{E}_0:=\mathcal{E}(0)$ satisfies
\begin{equation}\label{n0109}\mathcal{E}_0<C_0{\nu},
\end{equation}
 the initial-boundary value problem (\ref{0101})--(\ref{0104})
has a global weak solution $(\rho,\mathbf{v},\mathbf{d})$ on $I=I_T$ for
any given $T>0$, with the following properties:
\begin{enumerate}
            \item[\quad (1)]
Regularity:
\begin{align}%eqnarray}
&0\leq\rho\;\mbox{ a.e. in }Q_T,\quad \rho\in C^0(\bar{I},
{L}^\gamma_{\mm{weak}}(\Omega)
  \cap C^0(\bar{I},L^p(\Omega))\cap
 L^{\gamma+\eta}(Q_T), \label{fz0114} \\
& \mf{v}\in L^2(I,\mf{H}^1_0(\Omega)),\quad  \rho\mf{v}\in
L^\infty(I,\mf{L}^\frac{2\gamma}{\gamma+1}(\Omega))\cap
  C^0(\bar{I},\mf{L}^\frac{2\gamma}{\gamma+1}_{\mm{weak}}(\Omega)), \label{fz0115} \\
&  |\mf{d}|=1\;\mbox{ a.e. in }Q_T,\quad\mathbf{d}\in L^2(I,\mathbf{H}^2(\Omega))\cap
C^0(\bar{I},\mathbf{H}^1(\Omega)),\quad\partial_t\mathbf{d}\in {L}^{\frac{4}{3}}(I, \mathbf{L}^{2}(\Omega)),
\\
& \mathbf{d}|_{\partial\Omega}=\mathbf{d}_0|_{\partial\Omega}\in
 \mf{C}^{0,1}(\partial{\Omega})\;\mbox{ in the sense of trace for a.e. $t\in I$,}
\end{align}% eqnarray}
 where $p\in [1,\gamma)$,  and $\eta\in(0,\gamma-1)$.
    \item[\quad (2)] Equations (\ref{0101}) and (\ref{0102}) hold in $(\mathcal{D}'(Q_T))^3$,
    and equation (\ref{0103}) holds a.e. in $Q_T$.
 \item[\quad (3)] Equation (\ref{0101}) is
satisfied in the sense of renormalized solutions, that is, $\rho$,$\mathbf{v}$ satisfy
%%%
\begin{equation*}\label{0109}
\partial_tb(\rho)+\mathrm{div}[b(\rho)\mathbf{v}]+\left[\rho b'(\rho)-
b(\rho)\right]\mathrm{div}\mathbf{v}=0\mbox{ in }
\mathcal{D}'(\mathbb{R}^2\times I), \end{equation*} provided
($\rho$,$\mathbf{v}$) is prolonged to be zero on
$\mathbb{R}^2\setminus\Omega$, for any $b$ satisfying
\begin{equation}\label{0110}  b\in C^0[0,\infty)\cap C^1(0,\infty),\quad
|b'(s)|\leq cs^{-\lambda_0},\ s\in (0,1],\ \lambda_0<1, \end{equation}
and the growth conditions at infinity:
\begin{equation}\label{0111}|b'(s)|\leq c\,s^{\lambda_1},\
s\geq 1,\ where\ c>0,\ 0<1+\lambda_1<\frac{2\gamma-1}{2}.
\end{equation}
 \item[\quad (4)] The following finite and bounded energy inequalities hold:
\begin{align}
& \frac{d\mathcal{E}(t)}{dt}+\mathcal{F}(t)\leq 0\;\;
\mbox{ in }\mathcal{D}'(I),    \\
& \label{eniq1}\mathcal{E}(t)+\int_0^t\mathcal{F}(s)\mm{d}s\leq
\mathcal{E}_0\mbox{ for a.e. }t\in I.
\end{align}
\end{enumerate}
 \end{thm}

 \iffalse
\begin{rem}  In Theorem \ref{thm:0101}, we have used notations on Sobolev spaces that
\begin{equation*}  \begin{aligned}&
\mathbf{H}^k_0(\Omega):=({H}^k_0(\Omega))^2=(W^{k,2}_0(\Omega))^2,\quad
\mathbf{L}^p(\Omega):=(L^p(\Omega))^2\mbox{ for }p=2,\ {2\gamma}/({\gamma+1}),\\
&\mf{L}^\frac{2\gamma}{\gamma+1}_{\mm{weak}}(\Omega):
=({L}^\frac{2\gamma}{\gamma+1}_{\mm{weak}}(\Omega))^2,\quad
\mathbf{H}^k(\Omega):=(W^{k,2}(\Omega))^2.\end{aligned}\end{equation*}
%%%
Other notations with bold form will appear in the rest of this paper
(to denote the product spaces). For the convenience of readers, we
give the definition of the Sobolev spaces with weak topology: $
C^0(\bar{I},\mf{L}^q_{\mm{weak}}(\Omega)):= \{\mf{f}: I\rightarrow
\mf{L}^q(\Omega)~|~\int_\Omega \mf{f}\cdot\mf{g}\mm{d}\mf{x}\in
C(\bar{I})\mbox{ for any }\mf{g}\in \mf{L}^\frac{q}{q-1}(\Omega)\}$.
\end{rem}
\fi

\begin{rem} It should be noted that for the incompressible case \cite{DSJWHYSN}
the constant $C_0$ in a small energy condition similar to (\ref{n0109})
depends on the domain $\Omega$, since the authors in \cite{DSJWHYSN}
used the following interpolation inequality: Given a
bounded domain $\Omega\subset \mathbb{R}^2$, there exists an
interpolating constant ${c}_1(\Omega)$ depending on $\Omega$, such that
\begin{equation}\label{inp0110}
\|\nabla \mf{d}\|_{\mf{L}^4(\Omega)}^4\leq {c}_1(\Omega)(\|\nabla \mf{d}\|_{\mf{L}^{1}
(\Omega)}^2+\|\nabla^2 \mf{d}\|_{\mf{L}^2(\Omega)}^{2}\|\nabla\mf{d}\|_{\mf{L}^2(\Omega)}^{2})
\end{equation}
for any $\nabla \mf{d}\in \mf{H}^1(\Omega)$ (see \cite[Lemma 2.4]{DSJWHYSN}).
 In this paper, we use another version of interpolation
inequality (see (\ref{jjjw0432ww})), where the interpolating
constant takes the value of $2$. Thus our constant $C_0$ in
(\ref{n0109}) is independent of $\Omega$. In particular, from our
computations (see Proposition \ref{pro:energy}) we can take
$C_0=1/4096$ in (\ref{n0109}). Of course, this value is not optimal.
\end{rem}

\begin{rem}
In the above theorem, we assume that $\mathbf{d}_0\in\mathbf{H}^2({\Omega})$ for simplicity.
If this condition is replaced by $ \mathbf{d}_0\in \mathbf{H}^1({\Omega})$ and
$ \mathbf{d}_0\in\mathbf{C}^{2,\delta}({\partial\Omega})$ as in \cite{LFHLJYWCY}, then the
above theorem still holds. In addition, the proof of Theorem \ref{0101} remains basically unchanged
if the motion of the fluid is driven by a small bounded external force,
i.e., when (\ref{0102}) contains an additional term
$\rho\mf{f}(\mf{x},t)$ with $\mf{f}$ a small, bounded and measurable function.
%%%
% %             \item Theorem \ref{0101} still holds, if
%the condition $\mathbf{d}_0(\mathbf{x})\in \mathbf{H}^2({\Omega})$
%in (\ref{jfw0112}) is replaced by
%\begin{equation}\mathbf{d}_0(\mathbf{x})\in \mathbf{H}^1({\Omega})\mbox{
%and } \mathbf{d}_0(\mathbf{x})\in
%\mathbf{C}^{2,\delta}({\partial\Omega}),\end{equation} which was
%posed by Lin and Wang \cite{LFHLJYWCY} in the existence result of
%weak solutions.
 %           \end{itemize}
 %In addition,
% by tract theorem, for any positive integer
%$k$, there exists a constant $c$ depending on $\Omega$ such that
%$\|\mathbf{d}_0|_{\partial\Omega}\|
%_{\mf{H}^{k-\frac{1}{2}}({\partial\Omega})}\leq \|\mathbf{d}_0\|
%_{\mf{H}^{k}({\Omega})}$
\end{rem}

We now describe the main idea of the proof of Theorem
\ref{thm:0101}. For the Ginzburg-Landau approximation model to
(\ref{0101})--(\ref{0103}), based on some new estimates to deal with
the direction field and its coupling/interaction with the fluid variables,
Wang and Yu in \cite{WDHYCGA} adopted a classical
three-level approximation scheme which consists of the
Faedo-Galerkin approximation, artificial viscosity, artificial
pressure, and the celebrated weak continuity of the effective
viscous flux to overcome the difficulty of possible large
oscillations of the density, and established the existence of weak
solutions. These techniques were developed by Lions, Feireisl, et al
for the compressible Navier-Stokes equations in
\cite{LPLMTFM98,FENAPHOFJ35801,JSZPO}, and we refer to the monograph
\cite{NASII04} for more details.
% later, they are also adopted to investigate the existence of weak solutions
% to other physical models based on compressible Navier-Stokes
% equations, such as heat-conducting flows \cite{FETD2,FENA},
% magnetohydrodynamic flows \cite{HXWDHGELT}, radiative-reacting gas
% \cite{DDTKAM,DDTKOTC} and so on.
Compared with the Ginzburg-Landau approximation model in
\cite{WDHYCGA}, however, the system (\ref{0101})--(\ref{0103}) is
much more difficult to deal with due to the supercritical nonlinearity $|\nabla
\mf{d}|^2\mf{d}$ in (\ref{0103}). Consequently,
%and the strong coupling
%nonlinear term $\mm{div}(\nabla \mathbf{d}\odot\nabla
%\mathbf{d}-|\nabla \mf{d}|^2\mathbb{I}/2)$ in momentum equation
%(\ref{0102}).
one can not deduce sufficiently strong estimates on $\mf{d}$ from
the basic energy inequality (\ref{eniq1}) in the general case, such
as $\nabla^2\mf{d}\in L^2(I,\mf{L}^2(\Omega))$ as in \cite{WDHYCGA};
and on the other hand, one can not establish the global existence of
large solutions to the equation (\ref{0103}) for any given
$\mf{v}\in C^0(\bar{I},\mathbf{C}_0^2(\bar{\Omega}))$ as in the
Ginzburg-Landau approximation model.
%%%%
Recently,  Ding and Wen \cite{DSJWHYSN} obtained the global existence
and uniqueness of strong solutions to the two-dimensional density-dependent
incompressible model with small initial energy and positive initial
density away from zero. In \cite{DSJWHYSN}, they can deduce
$\nabla^2\mf{d}\in L^2(I,\mf{L}^2(\Omega))$ from the basic energy
inequality  in two dimensions under the condition of small initial
energy. In fact, they  first got ${\nu\|\nabla
\mf{d}\|_{\mf{L}^2(\Omega)}^2}/{2}\leq \mathcal{E}_0$ and
$\nu\theta\|\Delta \mf{d}\|_{\mf{L}^2(Q_T)}^2\leq
\mathcal{E}_0+\nu\theta\||\nabla \mf{d}||_{\mf{L}^4(Q_T)}^4$ from
the basic energy inequality.  Then, by interpolation inequality
(\ref{inp0110})  and elliptic estimates, they immediately obtain $\|\nabla^2
\mf{d}\|_{\mf{L}^2(Q_T)}^2\leq
C(\mathcal{E}_0,\|\nabla^2\mf{d}_0\|_{\mathbf{L}^2(\Omega)})$ for
some constant
$C(\mathcal{E}_0,\|\nabla^2\mf{d}_0\|_{\mathbf{L}^2(\Omega)})$
depending on the arguments, if initial energy is sufficiently small.
Motivated by this study,  %We find that this idea  can be applied to
for our compressible case we impose the restriction (\ref{n0109}) on the initial energy,
and consequently deduce the desired energy estimates on $\mf{d}$
from the basic energy inequality as in \cite{WDHYCGA}.
 With these estimates in hand, we also find that the three-level approximation scheme can
be applied to the initial-boundary value problem
(\ref{0101})--(\ref{0104}), if we can construct solutions to the
following third approximate problem:
\begin{eqnarray}&&  \label{n0117}
\rho_t+\mathrm{div}(\rho\mathbf{v})=\varepsilon \Delta \rho, \\
&&\begin{aligned}\label{n0118}&\int_\Omega
(\rho\mathbf{v})(t)\cdot\mathbf{\mathbf{\Psi}}\mathrm{d}\mathbf{x}
-\int_{\Omega}\mathbf{m}_0\cdot\mathbf{\mathbf{\Psi}}\mathrm{d}\mathbf{x} \\
& = \int_0^t\int_\Omega\bigg[\mu\Delta\mathbf{v}+(\mu+\lambda)\nabla
\mm{div}\mf{v}-A\nabla \rho^\gamma -\delta \nabla
\rho^\beta-\varepsilon(\nabla\rho\cdot\nabla
\mathbf{v})\\
&\qquad\qquad -\mathrm{div}(\rho\mathbf{v}\otimes\mathbf{v})
-\nu\mathrm{div}\left(\nabla \mathbf{d}\odot\nabla\mathbf{d}
-\frac{|\nabla\mf{d}|^2\mathbb{I}}{2}\right)\bigg]\cdot\mathbf{\Psi}\mathrm{d}\mathbf{x}\mathrm{d}s\\
&\mbox{\qquad  for all }t\in  I\mbox{ and any functions
}\mathbf{\Psi}\in \mathbf{X}_n,
\end{aligned}\\[1mm]
&&\label{n0119}\partial_t\mathbf{d}+\mathbf{v}\cdot \nabla
\mathbf{d}= \theta(\Delta \mathbf{d}+|\nabla
\mathbf{d}|^2\mf{d}),\end{eqnarray} where the $n$-dimensional
Euclidean space $\mf{X}_n$ will be introduced in Section
\ref{sec:04}, and $\varepsilon,\delta, \beta>0$ are constants. Fortunately, we can show the local existence of strong
solutions $\mf{d}\in C^0(\bar{I},\mf{H}^2(\Omega))$ to the equation
(\ref{n0118}) with smooth initial data for any given
$\mf{v}\in\{\mf{v}\in C^0(\bar{I},\mathbf{H}^{1}(\Omega)\cap
\mf{L}^\infty(\Omega))~|~\partial_t\mf{v}\in
L^{2}(I,\mathbf{H}^{1}(\Omega))\}$, with the aid of the  well-known
 Ehrling-Nirenberg-Gagliardo interpolation
inequality (see Lemma \ref{lem:0201}). With the result of local
existence in hand, adapting some techniques in
\cite{WDHYCGA,NASII04}, and making use of equivalence of norms in
$\mathbf{X}_n$ and the operator forms of $(\mf{v},\partial_t\mf{v})$
(see (\ref{n0512}) and (\ref{n0413nn}) for the operator forms
deduced from approximate momentum equations (\ref{n0118})), and
establishing some additional auxiliary estimates on
$\partial_t\mf{v}$, we can apply a fixed point theorem in some
bounded subspace of the Banach space $\{\mf{v}\in
C([0,T_n^*],\mathbf{X}_n) ~|~\partial_t\mf{v}\in
L^{2}((0,T_n^*),\mathbf{H}^{1}(\Omega))\}$ for some $T_n^*\in (0,T]$
to obtain the local existence of the third approximate problem.
Noting that the regularity of $\partial_t\mf{v}$ depends on
$\partial_t\rho$, this is why we need the properties on
$\partial_t\rho$ in Proportion \ref{pro:0401} on the solvability of
(\ref{n0117}). Then, we can repeat the fixed point argument to
extend the local solution to the whole time interval $I$ by
establishing the uniform-in-time energy estimates for the third
approximate problem.

Next, we state our second result on the existence of global weak
solutions to the Cauchy problem in $\mathbb{R}^2$ with large initial
data:
%%%%%%%%%%%%%%%%%%%%%%%%%%%%%%%%%%%%%%%%%%%%%%
\begin{thm}\label{thm:0102}
 Assume that, for some positive constants $\rho_\infty$ and
$\underline{d}_{02}$,  $(\rho_0,\mf{m}_0)$ satisfies (\ref{jfw0110}), (\ref{fzu0109})
with $\mathbb{R}^2$ in place of $\Omega$,
% Let the initial data ($\rho_0,\,\mathbf{m}_0,\,\mathbf{d}_0$)
% satisfy the compatibility conditions
and \begin{equation}\label{fxww0118}
 |\mf{d}_0|=1\mbox{ a.e. in }\mathbb{R}^2,\quad\mathbf{d}_0(\mathbf{x})-\mf{e}_2\in
\mathbf{H}^1({\mathbb{R}^2}),\quad  {d}_{02}\geq \underline{d}_{02},
                              \end{equation}where $d_{02}$
 denotes the second component of $\mf{d}_0$.
 Then, the Cauchy problem (\ref{0101})--(\ref{0104}) has a
global bounded energy weak solution $(\rho,\mathbf{v},\mathbf{d})$
for any given $T>0$ with the following properties:
\begin{enumerate}
            \item[\quad (1)]
Regularity:
\begin{eqnarray}
  &&0\leq\rho\mbox{ a.e. in }\mathbb{R}^2\times I, \quad
 Q(\rho)\in L^\infty(I,L^1(\mathbb{R}^2)),\nonumber\\
  && \rho\in C^0(\bar{I},
  {L}^\gamma_{\mm{weak}}(\Omega'))\cap C^0(\bar{I},L^p(\Omega'))
\cap L^{\gamma+\eta}_{\mm{loc}}(\mathbb{R}^2\times I)\nonumber,\\
&& \mf{v}\in L^2(I,({D}^{1,2}(\mathbb{R}^2))^2)\cap
L^2(I,\mf{L}^2_{\mm{loc}}(\Omega)),\quad \rho \mf{v}\in
  C^0(\bar{I},\mf{L}^\frac{2\gamma}{\gamma+1}_{\mm{weak}}(\Omega')),\nonumber\\
&& \rho |\mf{v}|^2\in
  L^\infty({I},\mf{L}^1(\mathbb{R}^2))\cap L^2(I,\mf{L}^{q}_\mm{loc}(\mathbb{R}^2)),
  \nonumber \\
&& \label{fz0125}\begin{aligned}&\partial_t\mathbf{d}\in
{L}^{\frac{4}{3}}(I, \mathbf{L}^{2}_{\mm{loc}}(\mathbb{R}^2)),\quad
|\mf{d}|=1\mbox{ a.e. in }\mathbb{R}^2\times I, \\
&\mathbf{d}\in C^0(\bar{I},\mathbf{H}^1_{\mm{loc}}(\mathbb{R}^2)),\quad \nabla\mf{d}\in
 L^2(I, \mf{H}^1({\mathbb{R}^2})\cap \mf{L}^4(\mathbb{R}^2\times I) ,
 \end{aligned}
\end{eqnarray}
where $p\in [1,\gamma)$, $\eta\in (0,\gamma-1)$, $q\geq 1$, and  $\Omega'\subset\mathbb{R}^2$
is an arbitrary bounded subdomain.
    \item[\quad (2)] Equations (\ref{0101}), (\ref{0102}) holds in
$\mathcal{D}'(\mathbb{R}^2\times I)$, equation (\ref{0103}) holds a.e. in $\mathbb{R}^2\times I$.
 Moreover, equation (\ref{0101}) is satisfied in the sense of renormalized solutions.
 \item[\quad (3)] The solution satisfies the energy equality
(\ref{eniq1}) with  $\mathbb{R}^2$ in place of $\Omega$.
%Energy estimate:
%\begin{equation}\label{fusafd0120}\begin{aligned}&\mathcal{E}(t)+\int_0^t\int_{\mathbb{R}^2} \left[\mu |\nabla
%\mf{v}|^2+(\lambda+\mu)|\mm{div}\mf{v}|^2
%+\frac{\theta\varpi_0}{2}(|\Delta \mf{d}|^2+|\nabla
%\mf{d}|^4)\right]\mathrm{d}\mathbf{x}\leq \mathcal{E}_0
%\end{aligned}\end{equation} for a.e. $t\in I$, where the energy $\mathcal{E}(t)$ is defined by (\ref{jsw0108}) with $\mathbb{R}^2$
%in place of $\Omega$, $Q(\rho)$ in $\mathcal{E}(t)$ is defined by
%(\ref{fzu0114}), and the constant $\varpi_0$ depends on
%$\underline{d}_{02}$, $\nu$ and $\mathcal{E}_0$.
\end{enumerate}
 \end{thm}

 \begin{rem} The notation ${D}^{1,2}(\mathbb{R}^2)$ in the above theorem denotes the homogeneous Sobolev
space on $\mathbb{R}^2$, i.e.,
$D^{1,2}(\mathbb{R}^2)=\overline{C_0^\infty(\mathbb{R}^2)}^{\|\nabla
\cdot\|_{\mf{L}^2(\mathbb{R}^2)}}$, where the symbol  ``overline with
norm" means completion with respect to that norm. Please refer to
\cite[Section 1.3.6]{NASII04}  for more basic conclusions concerning
the homogeneous Sobolev spaces.
 \end{rem}

\begin{rem}
In Theorem \ref{thm:0102}, the choice of $\mf{e}_2$ is for convenience only. In general,
one can choose any reference vector $\bf{e}\in \mathbb{S}^1$ and require that the image
of $\mf{d}$ is contained in a hemisphere around $\mf{e}$.
Let $\rho_\infty=0$ and suppose that in addition to the
assumptions above we have also $\rho_0\in L^1(\mathbb{R}^2)$. Then
the result above remains valid, moreover there holds
$\int_{\mathbb{R}^2}\rho\mm{d}\mf{x}=\int_{\mathbb{R}^2}\rho_0\mm{d}\mf{x}$
in $\bar{I}$. In addition, in view of  the proof of Theorem
\ref{0102} and the following interpolation inequality: there exists
a constant $c$ such that
\begin{equation*}  \|v\|_{L^4(B_R)}^4\leq
c\|  v\|_{\mf{H}^{1}(B_R)}^{2}\|v\|_{L^2(B_R)}^{2}\mbox{ for any
}v\in H^1(B_R)\mbox{ with }R\geq 1,
\end{equation*}
 the above theorem still holds  if we impose a  small initial energy
condition similar to (\ref{n0109}) to replace the  geometric
angle condition $ {d}_{02}\geq  \underline{d}_{02}$ in
(\ref{fxww0118}).
Here $B_R\subset \mathbb{R}^2$ is the ball of radius $R$.
%,  energy estimate (\ref{fusafd0120})
%should be modified by
%\begin{equation*}\mathcal{E}(t)+\int_0^t\int_{\mathbb{R}^2} \left[\mu |\nabla
%\mf{v}|^2+(\lambda+\mu)|\mm{div}\mf{v}|^2 +|\Delta \mf{d}|^2+|\nabla
%\mf{d}|^4\right]\mathrm{d}\mathbf{x}\leq C(\nu, \theta,c_0)
%\end{equation*} for a.e. $t\in I$, where the constant $C(\nu, \theta, c_0)$ depends on
%$\nu$, $\theta$, and  $c_0$. Here we have used the notation
%$B_R:=\{\mf{x}\in \mathbb{R}^2~|~|\mf{x}|<R\}$.
%In
%addition, we will also establish the large-time behavior of the
%global weak solutions. %Motivated by [8] and [13],
 \end{rem}

 %\begin{rem}In
%addition, we will also establish the large-time behavior of the
%global weak solutions. %Motivated by [8] and [13],
% \end{rem}

\begin{rem}\label{rem:0106}
Recently, Lin-Lin-Wang \cite{LFHLJYWCY}
established the existence of large weak solutions for the corresponding
incompressible case in a two-dimensional bounded domain. However, it is not clear
whether the method of localized small energy can be applied to the compressible model,
due to lack of regularity of $\|\mf{v}\|_{L^\infty(I,\mf{L}^2(\Omega))}$.
Hence it is still an interesting open problem  whether a general existence
result of large weak solutions holds in bounded domains for the compressible case.
In particular, we can not generalize the above result of the Cauchy problem
to the bounded domain case, since it is not clear whether the rigidity theorem
(see Proposition \ref{pro:0603}) holds in the bounded domain case.
\end{rem}

\begin{rem} When we completed this paper, we were informed that Wu and Tan \cite{WGTZGL}
 just finished the proof of the existence of
global weak solutions to (\ref{0101})--(\ref{0103}) in
$\mathbb{R}^3$ by using Suen and Hoff's method \cite{SADHG}, if the
initial energy is sufficiently small, the coefficients $\mu$ and
$\lambda$ satisfy $0\leq \lambda+\mu< \frac{3+\sqrt{21}}{6}\mu$, and
the initial data $(\mf{v}_0,\mf{v}_0)$ satisfies $\
\|\mathbf{v}_0\|_{\mathbf{L}^p(\mathbb{R}^3)}+\|\mathbf{d}_0\|_{\mathbf{L}^p(\mathbb{R}^3)}<\infty$
with $p>6$. However, it is not clear  whether their proof can be
applied to (\ref{0101})--(\ref{0103}) in a three-dimensional bounded
domain.
\end{rem}

Here we also describe the main idea of the proof of Theorem \ref{thm:0102}.
 Recently Lei, Li and Zhang \cite{LZLDZXY}  proved the rigidity theorem in
$\mathbb{R}^2$ (see Proposition \ref{pro:0603}) and obtained the estimate
on $\|\nabla^2\mf{d}\|_{L^2(I,\mf{L}^2(\mathbb{R}^2))}$ from the energy inequality,
if the second component of initial data of $\mf{d}$
satisfies some geometric angle condition. Motivated by this study, we impose the
geometric angle condition ${d}_{02}\geq  \underline{d}_{02}$ in
(\ref{fxww0118}) in place of (\ref{n0109}).
We first establish the local solvability of the Cauchy problem (\ref{n0119}) on $\mf{d}$
(see Proposition \ref{pro:0602}), which can be shown by following the
proof of the local solvability of $\mf{d}$ for the problem
(\ref{n0119}) defined in a bounded domain and using domain expansion technique.
Then, using the rigidity theorem, elliptic
estimates and interpolation inequalities in $\mathbb{R}^2$, we can
follow the proof of global existence of solutions to the third
approximate problem (\ref{n0117})--(\ref{n0119}) defined in a bounded
domain to establish the global existence of solutions to the approximate problem
(\ref{n0117}), (\ref{n0118}) defined in a ball $B_R$ and (\ref{n0119}) defined in
 $\mathbb{R}^2$ (see Proposition \ref{pro:0605}). Finally, we adapt the
 proof in \cite[Section 7.11]{NASII04} on the Navier-Stokes equations to prove
Theorem \ref{thm:0102}, by using Proposition \ref{pro:0605} with
the mass and momentum equations defined in the bounded invading
domain $B_R$ and letting $R\rightarrow \infty$.

 The rest of paper is organized as follows. In Section \ref{sec:02}, we deduce the basic energy equalities
 from (\ref{0101})--(\ref{0103}) and derive more energy estimates on
$\mf{d}$ under the assumption (\ref{n0109}). In Section
\ref{sec:04}, we list some preliminary results on solvability of two
sub-systems (\ref{n0117}) and (\ref{n0118}) in the third approximate
problem, while in Sections \ref{sec:05} we establish the unique
solvability of the third approximate problem. In Sections
\ref{sec:03}, exploiting the standard three-level approximation
scheme, we complete the proof of Theorem \ref{thm:0101}. Finally we
provide the proof of Theorem \ref{thm:0102} in Section \ref{Cauchy}.

\section{Energy estimates}\label{sec:02}

This section is devoted to formally deriving the energy estimates from
(\ref{0101})--(\ref{0103}) in a bounded domain $\Omega$, which will
play a crucial role in the proof of existence. Some of our results
have been  established  in \cite{LXWDHG,DSJWHYSN} for the
incompressible case.

 \subsection{Energy equalities}
We consider a classical solution $(\rho, \mathbf{v}, \mathbf{d})$ of
the problem (\ref{0101})--(\ref{0103}) with initial and boundary
conditions (\ref{0105}) and (\ref{0104}). First we verify that
\begin{equation} \label{0201}|\mf{d}|\equiv1\mbox{ in }Q_T, \mbox{ if }|\mf{d}_0|=1\mbox{ in }\Omega.
\end{equation}
 Multiplying the $\mf{d}$-system (\ref{0103}) by $\mf{d}$, we obtain
\begin{equation*} \label{0202}\begin{aligned}
\frac{1}{2}\partial_t|\mf{d}|^2+\frac{1}{2}\mf{v}\cdot\nabla
|\mf{d}|^2=\theta(\Delta\mf{d}\cdot\mf{d}+|\nabla
\mf{d}|^2|\mf{d}|^2).\end{aligned}\end{equation*}
Since
\begin{equation*} \label{0203}\begin{aligned}
\Delta|\mf{d}|^2=2|\nabla
\mf{d}|^2+2\Delta\mf{d}\cdot\mf{d},\end{aligned}\end{equation*}
it follows that
\begin{equation} \label{0205}\begin{aligned}
\partial_t(|\mf{d}|^2-1)-\theta\Delta(|\mf{d}|^2-1)+\mf{v}\cdot\nabla
(|\mf{d}|^2-1)-2\theta|\nabla\mf{d}|^2(|\mf{d}|^2-1)=0.\end{aligned}
\end{equation}
Multiplying (\ref{0205}) by $(|\mf{d}|^2-1)$ and then integrating over $\Omega$,
we use the boundary conditions to get
\begin{equation*} \label{0206}\begin{aligned}
\frac{d}{dt}\int_\Omega(|\mf{d}|^2-1)^2\mm{d}\mf{x}\leq & \int_\Omega(4\theta|\nabla
\mf{d}|^2+\mm{div}\mf{v})(|\mf{d}|^2-1)^2\mm{d}\mf{x}\\
\leq & \|4\theta|\nabla\mf{d}|^2+|\mm{div}\mf{v}|\|_{L^\infty(\Omega)}
\int_\Omega(|\mf{d}|^2-1)^2\mm{d}\mf{x}.
\end{aligned}\end{equation*}
We assume that $(\mf{v},\mf{d})$ satisfies the following
regularity:
\begin{equation*}
\|4\theta|\nabla
\mf{d}|^2+|\mm{div}\mf{v}|\|_{L^1(I,L^\infty(\Omega))}<\infty,
\end{equation*}
 then, using Gronwall's inequality, we immediately
verify (\ref{0201}). We shall see that all the couples
$(\mf{v}_n,\mf{d}_n)$ in the third approximate solutions constructed
in Section \ref{sec:05} satisfy the regularity above.

Similarly, we can verify
\begin{equation}\label{ji0615}\begin{aligned}
d_2\geq \underline{d}_{02}\;\mbox{ if }\; {d}_{02}\geq \underline{d}_{02}
\;\mbox{ for some given constant }\;\underline{d}_{02},
\end{aligned}\end{equation} which will play a crucial role in the
existence of large weak solutions to the Cauchy problem. Here we have
denoted the second component of $\mf{d}$ and $\mf{d}_0$ by ${d}_2$
and $d_{02}$, respectively.

Next, we give the deduction of \eqref{ji0615} for the reader's convenience.
Let
$$\omega=d_2-\underline{d}_{02}, \quad \omega^-=\min\{\omega,0\},$$
 we can deduce from (\ref{0103}) that
\begin{equation}\label{0615}\partial_t\omega-\theta \Delta\omega =\theta |\nabla
\mf{d}|^2(\omega+\underline{d}_{02})-\mf{v}\cdot\nabla \omega.
\end{equation}
Multiplying (\ref{0615}) by $\omega^-$ and integrating over
$\Omega$, we get
\begin{equation*}\begin{aligned}
& \frac{1}{2}\frac{d}{dt}\|\omega^-\|_{L^2(\mathbb{R}^2)}^2 +\theta \|\nabla
\omega^-\|_{L^2(\Omega)}^2  \\
& = \int_{\Omega} [\theta |\nabla\mf{d}|^2(\omega+\underline{d}_{02})
-\mf{v}\cdot\nabla\omega]\omega^-\mm{d} \mf{x}  \\
& \leq \left\|\theta |\nabla
\mf{d}|^2+\frac{1}{2}\mm{div}\mf{v}\right\|_{L^\infty(\Omega)}
\|\omega^-\|_{L^2(\Omega)}^2+\underline{d}_{02}\theta\||\nabla
\mf{d}|\|_{L^\infty(\Omega)} \int_{\Omega}|\nabla \omega^-|\omega^-\mm{d} \mf{x} \\
& \leq \left(\left\|\theta |\nabla
\mf{d}|^2+\frac{1}{2}\mm{div}\mf{v}\right\|_{L^\infty(\Omega)}
+\frac{1}{2}\underline{d}_{02}^2\theta\||\nabla
\mf{d}|\|_{L^\infty(\Omega)}^2 \right)\|\omega^-\|_{L^2(\Omega)}^2
+\frac{\theta }{2}\|\nabla \omega^-\|_{L^2(\Omega)}^2,
\end{aligned} \end{equation*}
which yields
\begin{equation*}\begin{aligned}\frac{d}{dt}\|\omega^-
\|_{L^2(\Omega)}^2 \leq \left(2\left\|\theta |\nabla
\mf{d}|^2+\frac{1}{2}\mm{div}\mf{v}\right\|_{L^\infty(\Omega)}
+\underline{d}_{02}^2\theta\||\nabla \mf{d}|\|_{L^\infty(\Omega)}^2
\right)\|\omega^-\|_{L^2(\Omega)}^2 .
\end{aligned}\end{equation*}
Hence, by Gronwall's  inequality, one obtains
\begin{equation*}\begin{aligned}\|\omega^-(t)
\|_{L^2(\Omega)}^2 \leq \|\omega^-(0) \|_{L^2(\Omega)}^2
e^{\int_0^t(2\left\|\theta |\nabla
\mf{d}|^2+\frac{1}{2}\mm{div}\mf{v}\right\|_{L^\infty(\Omega)}
+\underline{d}_{02}^2\theta\||\nabla
\mf{d}|\|_{L^\infty(\Omega)}^2)\mm{d}s}=0,
\end{aligned}\end{equation*}which implies (\ref{ji0615}).

Now we derive the energy equality. With the help of (\ref{0201}), we can deduce
the basic energy equality. To this end, we multiply equation (\ref{0102}) by $\mathbf{v}$
and integrate over $\Omega$ to deduce that
\begin{equation}\label{0210}\begin{aligned}
&\frac{d}{dt}\int_\Omega\left(\frac{1}{2}\rho|\mf{v}|^2+
Q(\rho)\right)\mm{d}\mf{x}+\int (\mu
|\nabla \mf{v}|^2+(\lambda+\mu)|\mm{div}\mf{v}|^2)\mathrm{d}\mathbf{x}\\
&=-\nu\int_\Omega (\nabla \mf{d})^T\Delta
\mf{d}\cdot\mf{v}\mathrm{d}x,\end{aligned}\end{equation}
where we have used integration by parts, the boundary condition
of $\mf{v}$, the mass equation (\ref{0101}) and the equality
\begin{equation*}
\mm{div}\left(\nabla \mathbf{d}\odot\nabla
\mathbf{d}-\frac{1}{2}|\nabla \mf{d}|^2\mathbb{I}\right)=(\nabla
\mf{d})^T\Delta \mf{d}:=(\partial_id_j)_{2\times 2}\Delta
\mf{d}.\end{equation*}

On the other hand, multiplying (\ref{0103}) by $-(\Delta
\mf{d}+|\nabla \mf{d}|^2\mf{d})$ and integrating over $\Omega$, we obtain
\begin{equation}\label{0211}-\int_\Omega\partial_t\mf{d}\cdot\Delta\mf{d}
\mm{d}\mf{x}-\int_\Omega(\mf{v}\cdot\nabla
\mf{d})\cdot\Delta\mf{d}\mm{d}\mf{x}=-\theta\int_\Omega |\Delta
\mf{d}+|\nabla \mf{d}|^2\mf{d}|^2\mm{d}\mf{x}, \end{equation}
where we have used the fact that $|\mf{d}|=1$ to get
\begin{equation*}\label{0212}(\partial_t\mf{d}+\mathbf{v}\cdot\nabla \mf{d})
\cdot|\nabla \mf{d}|^2\mf{d}=\frac{1}{2}(|\nabla
\mf{d}|^2\partial_t|\mf{d}|^2+|\nabla \mf{d}|^2\mf{v}\cdot\nabla
|\mf{d}|^2)=0.\end{equation*}
Using the boundary condition of $\mf{d}$, and integrating by parts,
we deduce from (\ref{0211}) that
\begin{equation}\label{0213}\frac{1}{2}\frac{d}{dt}
\int_\Omega|\nabla \mf{d}|^2\mm{d}\mf{x}+\theta\int_\Omega |\Delta
\mf{d}+|\nabla
\mf{d}|^2\mf{d}|^2\mm{d}\mf{x}=\int_\Omega(\mf{v}\cdot\nabla
\mf{d})\cdot\Delta\mf{d}\mm{d}\mf{x}. \end{equation}

Consequently, (\ref{0210})+(\ref{0213})$\times \nu$ gives the energy
equality in the differential form:
\begin{equation}\label{0214}\begin{aligned}&\frac{d}{dt}\mathcal{E}(t)+
\mathcal{F}(t)=0,
\end{aligned}\end{equation}
where the energy $\mathcal{E}(t)$ and dissipation  term
$\mathcal{F}(t)$ are defined by (\ref{jsw0108}) and
(\ref{dispassive}). Finally, integrating (\ref{0214}) with respect
to time over $(0,t)$, we get the energy equality
\begin{equation}\label{0215}\begin{aligned}&\mathcal{E}(t)
+\int_0^t\mathcal{F}(s)\mathrm{d}s=\mathcal{E}_0,\quad t\in I,
\end{aligned}\end{equation}
where $\mathcal{E}_0=\mathcal{E}(0)$ denotes the initial energy. Due
to the weak lower semicontinuity of norms on the right-hand side of
(\ref{0215}), for weak solutions one expects an inequality rather
than an equality.

\subsection{More \emph{a priori} estimates on $\mf{d}$ under the condition (\ref{n0109})}

In this subsection, we will deduce more estimates of $\mf{d}$ under
the assumption (\ref{n0109}), including the case of small initial
energy. For simplicity in the deduction, we always use the positive
constant $c_1(\Omega)$ to denote various constants depending on $\Omega$.

 % To this end, we first introduce the following well-known Gagliardo-Nirenberg inequality, the proof of
 % which can be found in \cite{VVAKAVOS3,DSJWHYSN} for example.
%\begin{lem}\label{lem:0201}
%Let $\Omega\subset \mathbb{R}^N$ be an arbitrary bounded domain
%satisfying the cone condition. Then the following inequality is
%valid for every function $v\in W^{1,m}(\Omega)$:
%\begin{equation}\label{n0216}\|v\|_{L^q(\Omega)}\leq \tilde{c}(\|v\|_{L^1(\Omega)}
%+\|\nabla
%v\|_{\mathbf{L}^m(\Omega)}^\alpha\|v\|_{L^r(\Omega)}^{1-\alpha}),
%\end{equation}
%where $m\geq 1$, $r\leq 1$ and $\alpha=(1/r-1/q)(1/r-1/m+1/N)^{-1}$,
%$\tilde{c}_1$ is a positive constant depending on $N$, $m$, $r$, $\alpha$ and
%$\Omega$ but not on $v$. Moreover, if $m<N$ and $N>1$,
%then $q\in [r,mN/(N-m)]$ for $r\leq mN/(N-m)$, and $q\in
%[mN/(N-m),r]$ for $r\geq mN(N-m)$. If $m\geq N>1$, then $q\in
%[r,+\infty)$ is arbitrary. If $m>N$, then $q\in [r,+\infty]$.
%If $m\leq N=1$, then $q\in [r,+\infty]$. In particular,
%if $v\in W_0^{1,m}(\Omega)$, then the term $\|v\|_{L^1(\Omega)}$ in
%(\ref{n0216}) can be deleted.
%\end{lem}
%To estimate the last two terms in (\ref{0430}),

Multiplying (\ref{0103}) by $\Delta \mf{d}$, integrating the
resulting equation over $\Omega$, and using integration by parts, we find that
\begin{equation}\label{0218}\frac{1}{2}\frac{d}{dt}
\int_\Omega|\nabla \mf{d}|^2\mm{d}\mf{x}+ \theta\int_\Omega |\Delta
\mf{d}|^2\mm{d}\mf{x}=\theta\int_\Omega |\nabla
\mf{d}|^4\mm{d}\mf{x}+ \int_\Omega(\mf{v}\cdot\nabla
\mf{d})\cdot\Delta\mf{d}\mm{d}\mf{x}. \end{equation}
Thus, (\ref{0210})$+$(\ref{0218})$\times \nu$ gives
\begin{equation}\label{0219}\begin{aligned}&\frac{d}{dt}\mathcal{E}(t)+\int_\Omega
\left[\mu |\nabla \mf{v}|^2+(\lambda+\mu)|\mm{div}\mf{v}|^2 +\nu\theta|\Delta
\mf{d}|^2\right]\mathrm{d}\mathbf{x}=\nu\theta \int_\Omega |\nabla
\mf{d}|^4\mathrm{d}\mathbf{x}.
\end{aligned}\end{equation}

Recalling $\Omega\subset \mathbb{R}^2$, exploiting \cite[Lemma 3.3]{NSRT1}, i.e.,
\begin{equation}\label{jjjw0432ww}\|v\|_{L^4(\Omega)}^4\leq
2\| \nabla v\|_{\mf{L}^{2}(\Omega)}^{2}\|v\|_{L^2(\Omega)}^{2}\;
\mbox{ for any }v\in H_0^1(\Omega),
\end{equation}
we can infer that
\begin{equation}\begin{aligned}\label{jjjwww1}\|\partial_id_j\|_{L^4(\Omega)}^4\leq
& \left(\|\partial_id_{0j}\|_{L^4(\Omega)}+2^{1/4}\|
\partial_i\nabla
(d_j-d_{0j})\|_{\mf{L}^{2}(\Omega)}^{1/2}\|\partial_i(d_j-
d_{0j})\|_{L^2(\Omega)}^{1/2}\right)^4\\
& \leq 8\left(\|\partial_id_{0j}\|_{L^4(\Omega)}^4+2\|\partial_i\nabla(d_j
-d_{0j})\|_{\mf{L}^{2}(\Omega)}^{2}\|\partial_i(d_j-d_{0j})\|_{L^2(\Omega)}^{2}\right) \\
& \leq 32\Big(c_1 (\Omega)\|\partial_id_{0j}\|_{H^1(\Omega)}^2
(\|\partial_id_{0j}\|_{H^1(\Omega)}^2 +\|\partial_id_j\|_{L^2(\Omega)}^{2}) \\
& \quad +\| \nabla\partial_i d_j\|_{\mf{L}^{2}(\Omega)}^{2}\|
\partial_i(d_j-d_{0j})\|_{L^2(\Omega)}^{2}\Big)\quad \mbox{for }1\leq i,j\leq 2.
\end{aligned}\end{equation}
Consequently,  making  use of (\ref{0215}), (\ref{jjjwww1}),
H\"{o}lder's inequality, and  the elliptic estimate
 \begin{equation}\begin{aligned}\label{0221}
 \|\nabla^2 {d}_j \|_{\mathbf{L}^2(\Omega)}\leq 2\|\Delta
 {d}_j\|_{{L}^2(\Omega)}+6\|\nabla^2{d}_0^j\|_{\mathbf{L}^2(\Omega)}
  \end{aligned}\end{equation}
 deduced from \cite[Collary 9.10]{GDTR},
%and the trace theorem, i.e., there exists a constant ${c}(\Omega)$
%depending on $\Omega$, such that
%\begin{equation*}\begin{aligned}\label{n0221}
%\|{d}_{0j}\|_{{H}^2(\Omega)}\leq
%{c}(\Omega)\|{d}_{0j}\|_{{H}^{\frac{3}{2}}(\partial\Omega)}\quad\mbox{
%for }j=1,\ 2,
%\end{aligned}\end{equation*}
we see that the last term in (\ref{0219}) can be estimated as follows:
%%%
\begin{equation}\begin{aligned}\label{n0228j}
&\nu\theta \int_\Omega |\nabla
\mf{d}|^4\mathrm{d}\mathbf{x}
\leq 4\nu\theta\sum_{1\leq i,j\leq 2}
\int_\Omega |\partial_i {d}_j|^4\mathrm{d}\mathbf{x}\\
 \leq & 128 \nu\theta\sum_{1\leq i,j\leq 2}\Big(c_1
(\Omega)\|\partial_id_{0j}\|_{H^1(\Omega)}^2(\|\partial_id_{0j}\|_{H^1(\Omega)}^2+\|\partial_id_j\|_{L^2(\Omega)}^{2})\\
&\qquad\qquad \qquad  +\| \nabla\partial_i
d_j\|_{\mf{L}^{2}(\Omega)}^{2}\|\partial_i(d_j-d_{0j})\|_{L^2(\Omega)}^{2}\Big)\\
\leq &128\nu\theta \sum_{1\leq i,j\leq 2}\Big(c_1
(\Omega)\|\partial_id_{0j}\|_{H^1(\Omega)}^2(\|\partial_id_{0j}\|_{H^1(\Omega)}^2+\|\partial_id_j\|_{L^2(\Omega)}^{2})\\
&+6\|\nabla^2{d}_{0j}\|_{\mathbf{L}^2(\Omega)}^2\|\partial_i(d_j-d_{0j})\|_{L^2(\Omega)}^{2}
+ 2\|\Delta
{d}_j\|^2_{L^2(\Omega)}\|\partial_i(d_j-d_{0j})\|_{L^2(\Omega)}^{2}\Big)\\
\leq & 128\nu\theta
(c_1(\Omega)+12)\|\mf{d}_0\|_{\mathbf{H}^2(\Omega)}^2\left(\|\mf{d}_0\|_{\mathbf{H}^2(\Omega)}^2+\frac{2\mathcal{E}_0}{\nu}
\right)+1024\theta\mathcal{E}_0\|\Delta
 \mf{d}\|_{L^2(\Omega)}^2.%\\
 %\leq &128\nu\theta
%(c(\Omega)+12)^3\|\mf{d}_0\|_{\mathbf{H}^{\frac{3}{2}}(\partial\Omega)}^2\left(\|\mf{d}_0\|_{\mathbf{H}^{\frac{3}{2}}(\partial\Omega)}^2+\frac{2\mathcal{E}_0}{\nu}
%\right)+1024\theta\mathcal{E}_0\|\Delta
% \mf{d}\|_{L^2(\Omega)}^2.
\end{aligned}\end{equation}

Choosing $\mathcal{E}_0>0$, such that
%\begin{equation*}
$\mathcal{E}_0 < \nu/4096$, we have
\begin{equation}\begin{aligned}\label{02241n}&\frac{d}{dt}\mathcal{E}(t)+\int_\Omega \left(\mu |\nabla
\mf{v}|^2+(\lambda+\mu)|\mm{div}\mf{v}|^2
+\frac{3\nu\theta}{4}|\Delta
\mf{d}|^2\right)\mathrm{d}\mathbf{x}\\
&\quad \leq 128\nu\theta
(c_1(\Omega)+12)\|\mf{d}_0\|_{\mathbf{H}^2(\Omega)}^2\left(\|\mf{d}_0\|_{\mathbf{H}^2(\Omega)}^2+\frac{2\mathcal{E}_0}{\nu}
\right):=g_1.
\end{aligned}\end{equation}
By (\ref{0221}), we obtain
\begin{equation*}\label{0224}\begin{aligned}\|\nabla^2
\mf{d}\|^2_{\mf{L}^2(Q_T)}\leq &
8\left(\frac{4(g_1T+\mathcal{E}_0)}{3\nu\theta}+9
\|\nabla^2\mf{d}_0\|_{\mathbf{L}^2(\Omega)}^2\right),
\end{aligned}\end{equation*}
where $$\|\nabla^2\mf{d}\|_{\mf{L}^2(\Omega)}^2:=\sum_{1\leq
i,j,l\leq 2}\|\partial_i\partial_j{d}_l\|_{{L}^2(\Omega)}^2.$$
Moreover, recalling $|\mf{d}|=1$, (\ref{02241n}) and (\ref{0215}),
we see that
\begin{equation}\begin{aligned}\label{n0225}\int_0^T\int_\Omega|\nabla
\mf{d}|^4\mathrm{d}\mathbf{x}\mm{d}t\leq
2\left(\sqrt{\frac{\mathcal{E}_0}{\nu\theta}}+\sqrt{\frac{4(g_1T+\mathcal{E}_0)}{3\nu\theta}}\right)^2=:g_2.
\end{aligned}\end{equation}

Finally, utilizing (\ref{0215}), (\ref{n0225}), H\"{o}lder's and Poincar\'e's inequalities,
 we get from the equation (\ref{0103}) that
\begin{equation*}\label{0226}
\begin{aligned}
 \|\partial_t{d}_j\|_{L^{4/3}(I,{L}^2(\Omega))}\leq &\|\mathbf{v}\cdot
 \nabla{d}_j
\|_{L^{4/3}(I,{L}^2(\Omega))}+\theta \|\Delta{d}_j
+|\nabla \mathbf{d}|^2{d}_j\|_{L^{4/3}(I,{L}^2(\Omega))}\\
\leq &{c}_1(\Omega)\|\nabla \mathbf{v}\|_{\mathbf{L}^{2}(Q_T)}\|
\nabla {d}_j\|_{\mathbf{L}^4(Q_T)}+\theta T^{1/4}\|\Delta {d}_j
+| \nabla \mathbf{d}|^2{d}_j\|_{{L}^2(Q_T)}  \\
\leq & {c}_1(\Omega)\sqrt{\frac{\mathcal{E}_0g_2}{\mu}}+
T^{1/4}\sqrt{\frac{\mathcal{E}_0\theta}{\nu}}.
\end{aligned}\end{equation*}
In addition, we can also deduce that
\begin{equation*}\label{0226d}\begin{aligned}
 \|\partial_t{d}_j\|_{L^{2}(I,({H}^{1}(\Omega))^*)}\leq & {c}_1(\Omega)\sqrt{\frac{\mathcal{E}_0}{\mu}}+
\sqrt{\frac{\mathcal{E}_0\theta}{\nu}},
\end{aligned}\end{equation*}
where $(H^{1}(\Omega))^*$ denotes the dual space of $H^1(\Omega)$.
%for some constant $c(\Omega)$ depending on $\Omega$.

Summing up the above estimates, we conclude that
%%%%%%%%%%%%%%%%%%%%%%%%%%%%%%%%%%%%%%%%
\begin{pro}\label{pro:energy}
Let $\Omega$ be a bounded domain. We have the following {a
priori} estimate for the initial-boundary value problem
(\ref{0101})--(\ref{0104}) with initial data $|\mf{d}_0|\equiv 1$:
\begin{equation}
\begin{aligned}\label{0226b}
%& |\mf{d}|\equiv1\\
 &\sup\limits_{t\in I}(\|\sqrt{\rho}{\bf
v}\|_{\mathbf{L}^2(\Omega)}+\|\rho\|_{{L}^\gamma(\Omega)}+\|\nabla
\mathbf{d}\|_{\mathbf{L}^2(\Omega)}) +\|\nabla
{\mathbf{v}}\|_{L^2(I,\mf{L}^2(\Omega))}\leq c_2(\mathcal{E}_0).
\end{aligned}\end{equation}
Moreover, if the initial energy satisfies $\mathcal{E}_0< \nu/4096$, then
\begin{equation}\label{0227}
\begin{aligned}
&\|\nabla^2 \mf{d}\|_{\mf{L}^2(Q_T)}+\|\nabla
\mf{d}\|_{\mf{L}^4(Q_T)}+\|\partial_t\mathbf{d}\|_{L^{4/3}(I,\mathbf{L}^2(\Omega))}
+\|\partial_t\mathbf{d}\|_{L^{2}(I,(\mf{H}^{1}(\Omega))^*)}\\
&\leq c_3(\|\mf{d}_0\|_{\mathbf{H}^{2}(\Omega)},\mathcal{E}_0).
\end{aligned}\end{equation}
Here $c_2$ and $c_3$ are positive constants which are, in
particular, nodecreasing in their variables. Moreover, $c_2$ also
depends on the given physical parameters $A$, $\mu$, $\nu$ and
$\theta$; and $c_3$ on $\nu$, $\theta$ and the domain $\Omega$.
\end{pro}

\section{Strong solvability of sub-systems in the third approximate problem}\label{sec:04}

In order to get a weak solution to the problem
(\ref{0101})--(\ref{0104}) in a bounded domain $\Omega$, we shall
first investigate the existence of solutions to the third
approximate problem of the original problem (\ref{0101})--(\ref{0104}):
\begin{eqnarray}&&\label{n301}\partial_t\rho+\mathrm{div}(\rho\mathbf{v})
=\varepsilon\Delta\rho ,   \\
&& \label{n302}\partial_t\mathbf{d}+\mathbf{v}\cdot \nabla \mathbf{d}
= \theta(\Delta \mathbf{d}+|\nabla \mathbf{d}|^2\mf{d}),  \\[1mm]
&& \begin{aligned}\label{n303}&\int_\Omega
(\rho\mathbf{v})(t)\cdot\mathbf{\mathbf{\Psi}}\mathrm{d}\mathbf{x}
-\int_{\Omega}\mathbf{m}_0\cdot\mathbf{\mathbf{\Psi}}\mathrm{d}\mathbf{x} \\
& = \int_0^t\int_\Omega\bigg[\mu\Delta\mathbf{v}+(\mu+\lambda)\nabla
\mm{div}\mf{v}-A\nabla \rho^\gamma -\delta \nabla
\rho^\beta-\varepsilon(\nabla\rho\cdot\nabla \mathbf{v})  \\
&\qquad\qquad -\mathrm{div}(\rho\mathbf{v}\otimes\mathbf{v})
-\nu\mathrm{div}\left(\nabla \mathbf{d}\otimes\nabla\mathbf{d} -\frac{|\nabla\mf{d}|^2
\mathbb{I}}{2}\right)\bigg]\cdot\mathbf{\Psi}\mathrm{d}\mathbf{x}\mathrm{d}s
 \end{aligned} \end{eqnarray}
 for all $t\in  I$  and any $\mathbf{\Psi}\in \mathbf{X}_n$, with boundary conditions
\begin{eqnarray}\nabla\rho\cdot\mf{n}|_{\partial\Omega}=0,\quad
\mathbf{v}|_{\partial \Omega}=\mathbf{0},\quad
\mathbf{d}|_{\partial\Omega}=\mathbf{d}_0,\end{eqnarray} and
modified initial data
\begin{eqnarray}\label{n0305}&&\rho(\mathbf{x},0)=\rho_0\in W^{1,\infty}(\Omega),
\quad 0<\underline{\rho}\leq\rho_0\leq \bar{\rho}<\infty,  \\
&&\label{n0306}\mathbf{d}(\mathbf{x},0)=\mathbf{d}_0\in\mf{H}^{3}({\Omega}),\quad
 \mathbf{v}(\mathbf{x},0)=\mf{v}_0\in \mathbf{X}_n,\end{eqnarray}
where  $\mf{n}$ denotes the outward normal to $\partial \Omega$,
and $\varepsilon,\delta, \beta, \underline\rho, \bar{\rho}>0$ are constants.
%
%$\widetilde{H}^{3/2}(\Omega):=\widetilde{W}^{2-\frac{1}{2},2}(\Omega)$
%is a completion of the space $\{\mathbf{z}\in C^\infty(\Omega)~|~\nabla
%\mathbf{z}\cdot\mf{n}|_{\partial\Omega}=0\}$ in
%$W^{2-\frac{1}{2},q}(\Omega)$ (see \cite[Section 1.3.5.10]{NASII04}
%for the Sobolev-Slobodetskii spaces of functions with ``fractional derivatives").
%
 Here we briefly introduce the finite dimensional space $\mf{X}_n$.
We know from \cite[Section 7.4.3]{NASII04} that there exist
countable sets
\begin{equation*}\begin{aligned}&\{\lambda_i\}_{i=1}^\infty,\
0<\lambda_1\leq \lambda_2\leq \cdots ,\mbox{ and }  \\
&\{\mathbf{\Psi}_i\}_{i=1}^\infty\subset\mathbf{W}_0^{1,p}(\Omega)\cap\mathbf{W}^{2,p}(\Omega),\
1\leq p<\infty,
\end{aligned}\end{equation*}
such that
\begin{equation*} \begin{aligned}
&-\mu\Delta\mathbf{\Psi}_i-(\mu+\lambda)\nabla
\mm{div}\mathbf{\Psi}_i=\lambda_i\mathbf{\Psi}_i,\qquad i=1,2,\cdots ,
\end{aligned}\end{equation*}
and $\{\mathbf{\Psi}_i\}_{i=1}^\infty$ is an orthonormal basis in
$\mf{L}^2(\Omega)$ and an orthogonal basis in $\mf{H}_0^1(\Omega)$
with respect to the scalar product
$\int_\Omega[\mu\partial_j\mf{u}\cdot\partial_j\mf{v}
+(\mu+\lambda)\mm{div}\mf{u}\,\mm{div}\mf{v}]\mm{d}\mf{x}$.
We define a $n$-dimensional Euclidean space $\mf{X}_n$ with scalar
product $<\cdot,\cdot >$ by
\begin{equation*}\begin{aligned}
\mf{X}_i=\mm{span}\{\mathbf{\Psi}_i\}_{i=1}^n,\quad
<\mf{u},\mf{v}> =\int_\Omega\mf{u}\cdot\mf{v}\mm{d}\mf{x},\quad \mf{u},\mf{v}\in\mf{X}_n,
\end{aligned}\end{equation*}
and denote by $\mathscr{P}_n$ the orthogonal projection of $\mf{L}^2(\Omega)$ onto $\mf{X}_n$.

In the next section, we will show the unique solvability of the
third approximate problem based on the idea from \cite[Section 7.7]{NASII04}
with certain restrictions imposed on the initial
approximate energy. Moreover, the unique solution enjoys the energy
estimates such as (\ref{0226b}), (\ref{0227}). As a preliminary, this
section is devoted to  the global solvability of the
Neumann problem (\ref{n301}) for the density and the local
solvability of the non-homogeneous Derichlet problem (\ref{n302})
when $\mf{v}$ is given. Due to the difficulty of the supercritical
nonlinearity $|\nabla \mf{d}|^2\mf{d}$, we need here higher
regularity imposed on the initial data $\mf{d}_0$ (see
(\ref{n0306})) in order to get the local strong solution $\mf{d}\in
{L^2(I_d^*,\mf{H}^3(\Omega))}$ for some time interval $I_d^*$. Thus,
using the embedding theorem, $\nabla \mf{d}\in
{L^2(I_d^*,\mf{L}^\infty(\Omega))}$. Such regularity is very
important to deduce the pointwise property of $\mf{d}$ (see (\ref{0201})).

\subsection{The Neumann problem for the density}

We consider the following problem:
\begin{eqnarray}\label{0401}\partial_t\rho+\mathrm{div}(\rho\mathbf{v})
=\varepsilon \Delta \rho, \end{eqnarray}
with initial and boundary conditions:
\begin{eqnarray} \label{0402}
&& \rho(\mathbf{x},0)=\rho_0(\mathbf{x})\in W^{1,\infty}(\Omega),
 \quad 0<\underline{\rho}\leq\rho_0(x)\leq \bar{\rho}<\infty,  \\
\label{0403}
&& \nabla\rho\cdot\mf{n}|_{\partial\Omega}=0,
\end{eqnarray}
which can be uniquely solved in terms of
$\mf{v}\in L^\infty (I,\mathbf{W}_0^{1,\infty}(\Omega))$.
The global well-posedness of the above problem can be read as follows
(see, e.g., \cite[Proposition 7.39]{NASII04}).
%%%%%%%%%%%%%%%%%%%%%%%%%%%%%%%%%%%%%
\begin{pro}\label{pro:0401}
Let $0<\alpha<1$, $\Omega$ be a bounded domain of class $C^{2,\alpha}$,
and $\rho_0$ satisfy (\ref{0402}). Then there exists a unique mapping
\begin{equation*}\begin{aligned}
\mathscr{S}_{\rho_0}:
L^{\infty}(I,\mathbf{W}^{1,\infty}_0(\Omega))\rightarrow\ C^0(\bar{I},H^{1}(\Omega)),
\end{aligned}\end{equation*}
such that
\begin{enumerate}
  \item[ \qquad (1)] $\mathscr{S}_{\rho_0}(\mathbf{v})$ belongs to the function class
  \begin{equation}\label{n0404} %\begin{aligned}
\mathcal{R}_T:=\left\{\rho|\,\rho\in L^2(I,W^{2,q}(\Omega))\cap C^0(\bar{I},W^{1,q}(\Omega)),
\partial_t\rho\in L^2(I,L^q(\Omega))\right\},1<q<\infty.%\end{aligned}
\end{equation}
  \item[ \qquad (2)]  The function $\rho=\mathscr{S}_{\rho_0}(\mathbf{v})$ satisfies (\ref{0401})
a.e. in $Q_T$, (\ref{0402}) a.e. in $\Omega$ and (\ref{0403}) in the
sense of traces a.e. in $I$.
\item[ \qquad (3)] $\mathscr{S}_{\rho_0}(\mathbf{v})$ is pointwise bounded, i.e.,
 \begin{eqnarray}\label{0404}\underline{\rho}e^{-\int_0^t\|\mf{v}\|_{\mf{W}^{1,\infty}(\Omega)}\mm{d}\tau}
 \leq \mathscr{S}_{\rho_0}(\mf{v})(\mf{x},t)\leq \bar{\rho}
 e^{\int_0^t\|\mf{v}(\tau)\|_{\mf{W}^{1,\infty}(\Omega)}\mm{d}\tau},\ t\in \bar{I}
 \;\mbox{ for a.e. }\mf{x}\in \Omega.\end{eqnarray}
\item[ \qquad (4)] If
$\|\mathbf{v}\|_{L^\infty(I,\mathbf{W}^{1,\infty}(\Omega))}\leq\kappa_v$, then
\begin{eqnarray}\label{0405}&&\|\mathscr{S}_{\rho_0}(\mathbf{v})\|_{L^\infty(I_t,H^{1}(\Omega))}\leq C_1\|\rho_0
\|_{H^1(\Omega)}e^{\frac{C_1}{2\varepsilon}(\kappa_{v}+\kappa^2_{v})t},\\
&&\|\nabla^2\mathscr{S}_{\rho_0}(\mathbf{v})\|_{\mf{L}^2(Q_t)}\leq
\frac{C_1}{\varepsilon}\sqrt{t} \|\rho_0\|_{H^1(\Omega)}\kappa_{v}
e^{\frac{C_1}{2\varepsilon}(\kappa_{v}+\kappa^2_{v})t}\\
&&\label{0407}\|\partial_t\mathscr{S}_{\rho_0}(\mathbf{v})\|_{L^2(Q_t)}\leq
C_1\sqrt{t}\|\rho_0\|_{H^1(\Omega)} \kappa_{v}
e^{\frac{C_1}{2\varepsilon}(\kappa_{v}+\kappa^2_{v})t}\end{eqnarray}
for any $t\in \bar{I}$, where $I_t:=(0,t)$, and $Q_t:=\Omega\times
I_t$. %Furthermore,
%\begin{equation}\begin{aligned}&\label{n0409}
%\varepsilon^{1-\frac{1}{r}}\|\mathscr{S}_{\rho_0}(\mf{v})\|_{L^\infty(I,
%H^{2-\frac{2}{r}}(\Omega))}+\|\partial_t\mathscr{S}_{\rho_0}(\mf{v})\|_{L^{r}(I,L^2(\Omega))}\\
%& \quad \leq C_1(\kappa_v,\varepsilon,r,T)\left( \|\rho_0\|_{{H}^{2-\frac{2}{r}}(\Omega)}+
%\|\rho_0\|_{H^1(\Omega)}
%e^{\frac{c}{2\varepsilon}(\kappa_v+\kappa^2_v)t}\right).\end{aligned}\end{equation}
 The constant $C_1$ in (\ref{0405})--(\ref{0407}) depends at most on $\Omega$, and
is independent of $\kappa_v$, $\varepsilon$, $T$, $\rho_0$ and $\mathbf{v}$.
%The constant $C_1$ is nondecreasing in the first variable and may depend on $\Omega$.
  \item[ \qquad (5)]  $\mathscr{S}_{\rho_0}(\mathbf{v})$ depends continuously  on $\mf{v}$,
  i.e.,
\begin{equation}\begin{aligned}\label{0409}\|[\mathscr{S}_{\rho_0}
(\mathbf{v}_1)-\mathscr{S}_{\rho_0}(\mathbf{v}_2)](t)\|_{L^2(\Omega)}\leq
C_2(\kappa_v,\varepsilon,T)t\|\rho_0\|_{H^1(\Omega)}
\|\mathbf{v}_1-\mathbf{v}_2\|_{L^\infty(I_t,\mathbf{W}^{1,\infty}(\Omega))},\end{aligned}
\end{equation}
\begin{equation}\begin{aligned}\label{n04091}\|\partial_t[\mathscr{S}_{\rho_0}
(\mathbf{v}_1)-\mathscr{S}_{\rho_0}(\mathbf{v}_2)](t)\|_{L^2(Q_t)}\leq
C_2(\kappa_v,\varepsilon,T)\sqrt{t}\|\rho_0\|_{H^1(\Omega)}
\|\mathbf{v}_1-\mathbf{v}_2\|_{L^\infty(I_t,\mathbf{W}^{1,\infty}(\Omega))}\end{aligned}
\end{equation}
for any $t\in \bar{I}$, and for any
$\|\mathbf{v}_1\|_{L^\infty(I,\mathbf{W}^{1,\infty}(\Omega))}\leq
\kappa_v$ and $\|\mathbf{v}_2\|_{L^\infty(I,\mathbf{W}^{1,\infty}(\Omega))}\leq
\kappa_v$. The constant $C_2$ is nondecreasing in the first
variable and may depend on $\Omega$.
\end{enumerate}
\end{pro}

\begin{pf} All the results above have been shown in \cite[Proposition 7.39]{NASII04}, except (\ref{n04091}).
Hence it suffices to verify (\ref{n04091}). For simplicity, we always
use the positive constant $C_2(\kappa_v,\varepsilon,T)$ to denote
various constants depending on $\kappa_v$, $\varepsilon$ and $T$.

Let $\rho_1,\rho_2\in {\cal R}_T$ be two solutions of the problem
(\ref{0401})--(\ref{0403}) with $\mf{v}$ replaced by $\mf{v}_1$ and $\mf{v}_2$ respectively.
Then after a straightforward calculation, we find that
\begin{equation}   \begin{aligned}\label{n0412}
&\frac{1}{2}\frac{d}{dt}\int_\Omega|\nabla
(\rho_1-\rho_2)|^2\mm{d}\mf{x}+\varepsilon\int_\Omega|\Delta(\rho_1-\rho_2)|^2\mm{d}\mf{x} \\
& =\int_\Omega\Big(\mf{v}_2\cdot\nabla
(\rho_1-\rho_2)+(\rho_1-\rho_2)\mathrm{div}\mf{v}_2+\rho_1\mm{div}(\mf{v}_1-\mf{v}_2) \\
&\qquad\qquad +(\mf{v}_1-\mf{v}_2)\cdot\nabla \rho_1\Big)
(\Delta\rho_1-\Delta\rho_2)\mm{d}\mf{x}.\end{aligned}
\end{equation}
Using Cauchy-Schwarz's and Poincar\'e's inequalities, the right hand
side of (\ref{n0412}) can be bounded from above by
%%%%%
\begin{equation*} \begin{aligned}
&\frac{C_1(\Omega)}{\varepsilon}\left(\|\mf{v}_2\|_{\mf{W}^{1,\infty}(\Omega)}^2\|\nabla
(\rho_2-\rho_1)\|_{\mf{L}^2(\Omega)}^2+\|\mf{v}_1-\mf{v}_2\|_{\mf{W}^{1,\infty}(\Omega)}^2
\|\rho_1\|_{H^1(\Omega)}^2\right)  \\
&\quad +\frac{\varepsilon}{2}\|\Delta
(\rho_1-\rho_2)\|_{L^2(\Omega)}^2\mbox{\quad for some constant
}C_1(\Omega)\mbox{ depending on }\Omega,\end{aligned}
\end{equation*}
which, together with (\ref{n0412}), yields
\begin{equation*}\begin{aligned}
&\frac{d}{dt}\int_\Omega|\nabla
(\rho_1-\rho_2)|^2\mm{d}\mf{x}+\varepsilon\int_\Omega|\Delta(\rho_1-\rho_2)|^2\mm{d}\mf{x}
\\
& \quad \leq
\frac{2C_1(\Omega)}{\varepsilon}\left(\|\mf{v}_2\|_{\mf{W}^{1,\infty}(\Omega)}^2\|\nabla
(\rho_2-\rho_1)\|_{\mf{L}^2(\Omega)}^2+\|\mf{v}_1-\mf{v}_2\|_{\mf{W}^{1,\infty}(\Omega)}^2
\|\rho_1\|_{H^1(\Omega)}^2\right).\end{aligned}
\end{equation*}
Therefore, by Gronwall's inequality and (\ref{0405}),
\begin{equation}\begin{aligned}\label{n0415}
\|\nabla (\rho_1-\rho_2)(t)\|_{\mf{L}^2(\Omega)}\leq
 C_2(\kappa_v,\varepsilon,T)\sqrt{t}\|\rho_0
\|_{H^1(\Omega)}\|\mf{v}_1-\mf{v}_2\|_{L^\infty(I_t,\mf{W}^{1,\infty}(\Omega))}.
\end{aligned}
\end{equation}
Moreover, we have
\begin{equation}\begin{aligned}\label{n0416}
\|\Delta (\rho_1-\rho_2)\|_{L^2(Q_t)}\leq
C_2(\kappa_v,\varepsilon,T)
\sqrt{t}\|\rho_0\|_{H^1(\Omega)}\|\mf{v}_1-\mf{v}_2\|_{L^\infty(I_t,\mf{W}^{1,\infty}(\Omega))}.
\end{aligned}
\end{equation}

Consequently, using (\ref{n0415}), (\ref{n0416}), we get from (\ref{0401}) that
\begin{equation*}\begin{aligned}
&\|\partial_t(\rho_1-\rho_2)\|_{L^2(Q_t)}\\
&=\|(\mf{v}_2\cdot\nabla
(\rho_1-\rho_2)+(\rho_1-\rho_2)\mathrm{div}\mf{v}_2+\rho_1\mm{div}(\mf{v}_1-\mf{v}_2)\\
&\quad +(\mf{v}_1-\mf{v}_2)\cdot\nabla \rho_1)-\varepsilon\Delta (
\rho_1- \rho_2)\|_{L^2(Q_t)}\\
&\leq C_1(\Omega)\left(\kappa_v\|\nabla
(\rho_1-\rho_2)\|_{L^2(Q_t)}+\sqrt{t}\|\rho_1\|_{L^\infty(I_t,H^1(\Omega))}
\|\mathbf{v}_1-\mathbf{v}_2\|_{L^\infty(I_t,\mathbf{W}^{1,\infty}(\Omega))}\right)\\
&\quad +\varepsilon \|\Delta( \rho_1 - \rho_2)\|_{L^2(Q_t)}  \\
&\leq C_2(\kappa_v,\varepsilon,T)
\sqrt{t}\|\rho_0\|_{H^1(\Omega)}\|\mf{v}_1-\mf{v}_2\|_{L^\infty(I_t,\mf{W}^{1,\infty}(\Omega))},
\end{aligned}
\end{equation*} which implies (\ref{n04091}).
\end{pf}

\subsection{The local solvability of the direction vector}\label{nsectionjf:0302}

Now we turn to  show the local solvability of the direction vector. More precisely, we will show
that for
\begin{equation*}\begin{aligned}\mathbf{v}\in  \mathbb{V}:=
\{\mathbf{v}~|~&\mathbf{v}\in
C^0(I,\mathbf{H}^{1}(\Omega)\cap\mf{L}^\infty(\Omega)),\
\partial_t\mf{v} \in
L^{2}(I,\mathbf{H}^{1}(\Omega))\},\end{aligned}\end{equation*} there
exists a $T_d^*\in (0,T]$, such that the equation
\begin{eqnarray}\label{0411}\partial_t\mathbf{d}+\mathbf{v}\cdot \nabla \mathbf{d}= \theta(\Delta
\mathbf{d}+|\nabla \mathbf{d}|^2\mf{d})\end{eqnarray} has a unique
strong solution $\mf{d}(\mf{x},t)$ on $[0,T^*_d)$ satisfying the
initial and boundary conditions
\begin{eqnarray}\label{0412}
&&\mathbf{d}(\mf{x},0)=\mathbf{d}_0(\mf{x})\in \mf{H}^{3}({\Omega}),\\ \label{0413}
&&\mathbf{d}|_{\partial\Omega}=\mathbf{d}_0(\mf{x})\quad\mbox{ for
}t\in (0,T^*_d).\end{eqnarray}
This local existence can be shown by modifying
the arguments in \cite[Section 3]{DSJWHYSN}. Here we give its proof,
since some arguments in our proof are very different
from those in \cite[Section 3]{DSJWHYSN}, and can be applied to the Cauchy problem
in Section \ref{Cauchy} by using the Ehrling-Nirenberg-Gagliardo interpolation
inequality and domain expansion technique. Moreover, we shall show that the solution
$\mf{d}$ continuously depends on $\mf{v}$, which will be also used in the proof
of the local existence of solutions to the third approximate problem.

\subsubsection{Linearized problems}
Denote
\begin{equation*}\left\{\begin{aligned}
 &\mathbb{D}^*:=\{\mf{b}\in C^0(\bar{I}_d^*,\mf{H}^2(\Omega))\; |\; \partial_t\mf{b}
 \in {C^0(\bar{I}_d^*,\mf{L}^2(\Omega))}\cap
{L^2(I_d^*,\mf{H}^1_0(\Omega))}\}, \\
&\mathbb{V}_K:=\{\mf{v}\in \mathbb{V}|\|\mf{v}\|_{\mathbb{V}}\leq K\},
\\
& {\mathbb{D}}^*_{\kappa_d}:=\{\mf{b}\in
\mathbb{D}^*~|~\mf{b}(\mf{x},0)=\mf{d}_0(\mf{x}),\
\|\mf{b}\|_{\mathbb{D}^*}\leq\kappa_d\}, \end{aligned}\right. \end{equation*}
where $K$ and $\kappa_d$ are positive constants,
\begin{equation}\label{0415}\begin{aligned}
&\|\mf{v}\|_{\mathbb{V}}:= \left(\|\mathbf{v}\|_{C^0
(\bar{I},\mathbf{H}^{1}(\Omega))}^2+\|\mathbf{v}\|_{C^0 (\bar{I},\mathbf{L}^{\infty}(\Omega))}^2
+\|\partial_t\mf{v}\|_{L^{2}(I,\mathbf{H}^{1}(\Omega))}^2\right)^\frac12, \\
&\|\mf{b}\|_{\mathbb{D}^*}:=\left(\|\partial_t\mf{b}\|_{L^2(I_d^*,\mf{H}^1(\Omega))}^2+\|
\partial_t\mf{b}\|_{L^\infty(I_d^*,\mf{L}^2(\Omega))}^2+\|\mf{b}\|_{L^\infty(I_d^*,\mf{H}^2(\Omega))}^2
+\|\mf{b}\|_{L^2(I_d^*,\mf{H}^3(\Omega))}^2\right)^\frac12,
\end{aligned}\end{equation}
and $I_d^*:=(0,T_d^*)\subset I$.  Without loss of generality, we assume $\kappa_d\geq 1$.

 To show the local existence of
solutions, we will construct a sequence of approximate solutions to
(\ref{0411})--(\ref{0413}) and use the technique of iteration based
on the set ${\mathbb{D}}^*_{\kappa_d}$. First we linearize the
original system (\ref{0411}) for $\mf{b}\in
{\mathbb{D}}^*_{\kappa_d}$ as follows:
\begin{eqnarray}\label{0416}
\partial_t\mf{d}-\theta\Delta\mf{d}=\theta|\nabla \mf{b}|^2\mf{b}
-\mf{v}\cdot\nabla\mf{b}\end{eqnarray}
with a given function $\mf{v}\in \mathbb{V}_K$, and initial and boundary conditions
\begin{eqnarray}\label{0417}
&&\mf{d}|_{t=0}=\mf{d}_0,\ \mf{x}\in \Omega, \\
&&\label{0418} \mf{d}|_{\partial\Omega}=\mf{d}_0,\ t\in
I_d^*.\end{eqnarray} Since $F(\mf{b}):=|\nabla
\mf{b}|^2\mf{b}-\mf{v}\cdot\nabla\mf{b}\in L^2
(I_d^*,\mf{H}^1(\Omega))$ and $\partial_tF(\mf{b})\in L^2
(I_d^*,(\mf{H}^{1}(\Omega))^*)$, we can apply a standard method,
such as a semidiscrete Galerkin method in \cite{ELGP}, to prove that
the initial-boundary value problem (\ref{0416})--(\ref{0418}) has a
unique solution
\begin{eqnarray*}\begin{aligned}
&\mf{d}\in {C^0(\bar{I}_d^*,\mf{H}^2(\Omega)) \cap
L^2(I_d^*,\mf{H}^3(\Omega))},\quad \partial_{tt}^2\mf{d}\in
{L^2(I_d^*,(\mf{H}^{1}(\Omega))^*)},\\
& \partial_t\mf{d}\in {C^0(\bar{I}_d^*,\mf{L}^2(\Omega))}\cap
{L^2(I_d^*,\mf{H}^1_0(\Omega))}.\end{aligned}\end{eqnarray*} Here we
have used $(\mf{H}^{1}(\Omega))^*$ to denote the dual space of
$\mf{H}^{1}(\Omega)$. Therefore, we can get a solution $\mf{d}^1$ to
(\ref{0416}) with $\mf{b}$ replaced by  some given $\mf{d}^0\in
{\mathbb{D}}^*_{\kappa_d}$. Assuming $\mf{d}^{k-1}\in
{\mathbb{D}}^*_{\kappa_d}$ for $k\geq 1$, we can construct an
approximate solution $\mf{d}^k$ satisfying
\begin{eqnarray}\label{0420}
\partial_t\mf{d}^k-\theta\Delta\mf{d}^k=\theta|\nabla \mf{d}^{k-1}|^2\mf{d}^{k-1}-
\mf{v}\cdot\nabla\mf{d}^{k-1}\end{eqnarray}
with initial and boundary conditions:
\begin{eqnarray*}
&&\mf{d}^k|_{t=0}=\mf{d}_0,\ \mf{x}\in \Omega,\\
&& \mf{d}^k|_{\partial\Omega}=\mf{d}_0,\ t\in I_d^*.\end{eqnarray*}
%%%%%%%%%%%%%%%%%

\subsubsection{Uniform estimates}\label{030202}

Next we derive the uniform estimates:
\begin{eqnarray}\label{0423}&&\|\mf{d}^k\|_{\mathbb{D}^*}\leq \kappa_d
\end{eqnarray}
for some constant $\kappa_d$ depending on $\mf{d}_0$ and the given
constant $K$. We mention that in the following estimates the letter
$c \geq 1$ will denote various positive constants independent of
$\kappa_d$, $K$, $\mf{v}$, $\mf{d}_0$ and $\Omega$; the letter
$C_1(\Omega)$ will denote various positive constants depending on
$\Omega$, and the letter $C(\ldots)$ various positive constants
depending on its variables (it may depend on
 $\Omega$ sometimes, however we omit them for simplicity), and is
nondecreasing in its variables. Of course,
they may also depend on the fixed value $\theta$.

First, we can deduce from (\ref{0420})  and integration by parts that
\begin{equation*}\begin{aligned}  &\frac{1}{2}\frac{d}{dt}\int_\Omega| \mf{d}^k|^2\mm{d}\mf{x}
+\theta\int_\Omega|\nabla \mf{d}^k|^2\mm{d}\mf{x}\\
&=\theta\int_\Omega|\nabla \mf{d}^{k-1}|^2\mf{d}^{k-1}\cdot
\mf{d}^k\mm{d}\mf{x}-\int_\Omega\mf{v}\cdot \nabla \mf{d}^{k-1}\cdot
  \mf{d}^k\mm{d}\mf{x}+\theta\int_{\partial\Omega} \mf{n}\cdot\nabla \mf{d}_0\cdot\mf{d}_0 \mm{d}\mf{s}
   \end{aligned}\end{equation*} and
\begin{equation*}\begin{aligned}
&\frac{1}{2}\frac{d}{dt}\int_\Omega|\nabla \mf{d}^k|^2\mm{d}\mf{x}
+\theta\int_\Omega|\Delta\mf{d}^k|^2\mm{d}\mf{x}\\
&=-\theta\int_\Omega|\nabla \mf{d}^{k-1}|^2\mf{d}^{k-1}\cdot\Delta
\mf{d}^k\mm{d}\mf{x}+\int_\Omega\mf{v}\cdot \nabla \mf{d}^{k-1}\cdot
\Delta \mf{d}^k\mm{d}\mf{x},\end{aligned}\end{equation*}
which, together with Cauchy-Schwarz's and H\"older's  inequalities, imply that
\begin{equation*}\begin{aligned}
\frac{1}{2}\frac{d}{dt}\| \mf{d}^k\|^2_{\mf{H}^1(\Omega)}
&\leq \left(\theta\||\nabla
\mf{d}^{k-1}|\|^4_{{L}^4(\Omega)}\||\mf{d}^{k-1}|\|^2_{{L}^\infty(\Omega)}+\||\mf{v}|\|_{{L}^4(\Omega)}^2
\||\nabla \mf{d}^{k-1}|\|_{{L}^4(\Omega)}^2\right)\\
& \qquad +\|\mf{d}^k\|^2_{\mf{L}^2(\Omega)}+\theta\|\nabla
\mf{d}_0\|_{\mf{L}^2(\partial\Omega)}\|\mf{d}_0\|_{\mf{L}^2(\partial\Omega)}.
\end{aligned}\end{equation*}
Hence, using the embeddings $H^2(\Omega)\hookrightarrow
L^\infty(\Omega)$ for $\mf{d}$, $H^1(\Omega)\hookrightarrow
L^4(\Omega)$ for $(\nabla \mf{d},\mf{v})$, and trace theorem, we have
\begin{equation}\begin{aligned}\label{jjwfsaf123}
\frac{d}{dt}\| \mf{d}^k\|^2_{\mf{H}^1(\Omega)}
&\leq C(K,\|\mf{d}_0\|_{\mf{H}^2(\Omega)})(\kappa_d^6+\kappa_d^2K^2 )
+2\|\mf{d}^k\|^2_{\mf{H}^1(\Omega)}, \end{aligned}\end{equation}
which, together with Gronwall's inequality, implies
\begin{equation*}\begin{aligned}
\| \mf{d}^k(t)\|^2_{\mf{H}^1(\Omega)} & \leq \left(\|\mf{d}_0\|^2_{\mf{H}^1(\Omega)}
+C(K,\|\mf{d}_0\|_{\mf{H}^2(\Omega)})(\kappa_d^6+\kappa_d^2K^2)t\right)e^{2t}
\end{aligned}\end{equation*}
for any $t\in I_d^*$. Taking $T_1^*=\kappa_d^{-6}\leq 1$ and letting $T_d^*\leq T^*_1$,
we conclude
\begin{equation}\begin{aligned}\label{js0330n}\| \mf{d}^k\|^2_{\mf{H}^1(\Omega)}
&\leq C(K,\|\mf{d}_0\|_{\mf{H}^2(\Omega)}).
\end{aligned}\end{equation}

Next we derive bounds on $\partial_t\mf{d}$. It follows from
(\ref{0420}) and integration by parts that
\begin{equation}\begin{aligned}\label{enes0331}
&\frac{1}{2}\frac{d}{dt}\int_\Omega| \partial_t\mf{d}^k|^2\mm{d}\mf{x}
+\theta\int_\Omega|\nabla \partial_t\mf{d}^k|^2\mm{d}\mf{x} \\
&=-2\theta\bigg(\int_\Omega(\Delta \mf{d}^{k-1}\cdot
\partial_t\mf{d}^{k-1})\mf{d}^{k-1}\cdot
\partial_t\mf{d}^k\mm{d}\mf{x}+\int_\Omega\nabla \mf{d}^{k-1}_i
\cdot\nabla \mf{d}^{k-1}_j\partial_t\mf{d}^{k-1}_i\partial_t\mf{d}^k_j \mm{d}\mf{x}\\
&\quad +\int_\Omega \partial_t\mf{d}^{k-1}_i\mf{d}^{k-1}_j \nabla \mf{d}^{k-1}_i\cdot\nabla
\partial_t\mf{d}^k_j\mm{d}\mf{x}\bigg)+\theta\int_\Omega|\nabla
\mathbf{d}^{k-1}|^2\partial_t\mathbf{d}^{k-1}\cdot \partial_t\mf{d}^{k}\mm{d}\mf{x}\\
&\quad -\int_\Omega\mf{v}\cdot \nabla \partial_t\mf{d}^{k-1}\cdot
\partial_t\mf{d}^k\mm{d}\mf{x}-\int_\Omega\partial_t\mf{v}\cdot
\nabla \mf{d}^{k-1}\cdot  \partial_t\mf{d}^k\mm{d}\mf{x}.
&\end{aligned}\end{equation}
Noting that $(\partial_t\mf{d},\partial_t\mf{v},\mf{v})=\mf{0}$ on $\partial\Omega$,
applying Cauchy-Schwarz's and H\"older's inequalities, and
(\ref{jjjw0432ww}), the right hand side of (\ref{enes0331}) can be
bounded from above by
\begin{equation*}\begin{aligned}\label{0333jwd}&
 c\Big((\|\Delta
\mf{d}^{k-1}\|^2_{\mf{L}^2(\Omega)}
\|\mf{d}^{k-1}\|^2_{\mf{L}^\infty(\Omega)}+\|\nabla
\mf{d}^{k-1}\|^4_{\mf{L}^4(\Omega)} )\|\partial_t
\mf{d}^{k-1}\|_{\mathbf{L}^2(\Omega)}\|\nabla\partial_t
\mf{d}^{k-1}\|_{\mathbf{L}^2(\Omega)}^\frac{1}{2}\\
&\quad +\|\mf{d}^{k-1}\|^2_{\mf{L}^\infty(\Omega)}\|\nabla
\mf{d}^{k-1}\|^2_{\mf{L}^4(\Omega)}\|\partial_t
\mf{d}^{k-1}\|_{\mathbf{L}^2(\Omega)}\|\nabla\partial_t
\mf{d}^{k-1}\|_{\mathbf{L}^2(\Omega)} \\
&\quad +\|\mf{v}\|_{\mf{L}^\infty(\Omega)}^2 \|\nabla\partial_t
\mf{d}^{k-1}\|_{\mathbf{L}^2(\Omega)}+\|\nabla
\mf{d}^{k-1}\|_{\mathbf{L}^4(\Omega)}^2 \\
&\quad +(\|\nabla\partial_t
\mf{d}^{k-1}\|_{\mathbf{L}^2(\Omega)}+\|\partial_t\mf{v}\|_{\mf{L}^2(\Omega)}\|\nabla
\partial_t\mf{v}\|_{\mf{L}^2(\Omega)})\|\partial_t
\mf{d}^{k}\|_{\mathbf{L}^2(\Omega)}^2\Big) +\frac{\theta}{2}\|\nabla\partial_t
\mf{d}^{k}\|_{\mathbf{L}^2(\Omega)}^2. &\end{aligned}\end{equation*}
Using $H^2(\Omega)\hookrightarrow L^\infty(\Omega)$ and
$H^1(\Omega)\hookrightarrow L^4(\Omega)$ for $\mf{d}$ and $\nabla\mf{d}$
respectively, the above term can be further bounded from above by
\begin{equation}\begin{aligned}\label{nsfs0333jwd}&
 C_1(\Omega)\Big((\kappa^5+K^2)(\|\nabla\partial_t
\mf{d}^{k-1}\|_{\mathbf{L}^2(\Omega)}+\|\nabla\partial_t
\mf{d}^{k-1}\|_{\mathbf{L}^2(\Omega)}^{\frac{1}{2}}) +\kappa^2_d \\
&\quad +(\|\nabla\partial_t
\mf{d}^{k-1}\|_{\mathbf{L}^2(\Omega)}+\|\partial_t\mf{v}\|_{\mf{L}^2(\Omega)}\|\nabla
\partial_t\mf{v}\|_{\mf{L}^2(\Omega)})\|\partial_t
\mf{d}^{k}\|_{\mathbf{L}^2(\Omega)}^2\Big)+\frac{\theta}{2}\|\nabla\partial_t
\mf{d}^{k}\|_{\mathbf{L}^2(\Omega)}^2.\end{aligned}\end{equation}
Hence, we have
\begin{equation}\begin{aligned}\label{js0334}&\frac{d}{dt}\| \partial_t\mf{d}^k\|^2_{\mf{L}^2(\Omega)}
+\theta\|\nabla \partial_t\mf{d}^k\|^2_{\mf{L}(\Omega)}\\ &\leq
C_1(\Omega)\Big((\kappa^5+K^2)(\|\nabla\partial_t
\mf{d}^{k-1}\|_{\mathbf{L}^2(\Omega)}+\|\nabla\partial_t
\mf{d}^{k-1}\|_{\mathbf{L}^2(\Omega)}^{\frac{1}{2}})+\kappa^2_d \\&
\qquad\qquad +(\|\nabla\partial_t
\mf{d}^{k-1}\|_{\mathbf{L}^2(\Omega)}+\|\partial_t\mf{v}\|_{\mf{L}^2(\Omega)}\|\nabla
\partial_t\mf{v}\|_{\mf{L}^2(\Omega)})\|\partial_t
\mf{d}^{k}\|_{\mathbf{L}^2(\Omega)}^2\Big).\end{aligned}\end{equation}
Consequently, using Gronwall's and H\"older's inequalities, we get
\begin{equation*}\begin{aligned}\label{nn0337}
\| \partial_t\mf{d}^k(t)\|^2_{\mf{L}^2(\Omega)} \leq
\Big(\|\partial_t\mf{d}_0\|^2_{\mf{L}^2(\Omega)}+C_1(\Omega)(\kappa_d^{5}+K^2)( t
+\kappa_d\sqrt{t}+\sqrt{\kappa_d}t^{\frac{3}{4}})\Big)e^{C_1(\Omega)(\kappa_d\sqrt{t}+K^2)}
\end{aligned}\end{equation*}
for any $t\in I_d^*$. Now, taking $T_2^*=\kappa_d^{-12}\leq 1$,
letting $T_d^*\leq T^*_2$, and noting that
$$\|\partial_t\mf{d}(0)\|_{\mathbf{L}^2(\Omega)}=\|\theta(\Delta
\mathbf{d}_0+|\nabla \mathbf{d}_0|^2\mf{d}_0)-\mathbf{v}_0\cdot
\nabla \mathbf{d}_0\|_{\mf{L}^2(\Omega)},$$
 we conclude
\begin{equation}\begin{aligned}\label{0428}
\sup_{t\in I_d^*} \|\partial_t\mf{d}^k\|_{\mathbf{L}^2(\Omega)} \leq
C(K,\|\mf{d}_0\|_{\mf{H}^2(\Omega)}). \end{aligned}
\end{equation}
Moreover, using  (\ref{0420}), (\ref{js0334}), (\ref{0428}),
Cauchy's inequality,  and the elliptic theory, we  obtain the following estimates,
 \begin{equation}\begin{aligned}\label{0429}
\|\partial_t\mf{d}^k\|_{L^2(I_d^*,\mf{H}^1(\Omega)}\leq
C(K,\|\mf{d}_0\|_{\mf{H}^2(\Omega)}),\end{aligned}\end{equation}
\begin{equation}\begin{aligned}\label{0430}
\|\nabla^2\mf{d}^k\|_{\mf{L}^2(\Omega)}^2\leq&
c\Big(\|\nabla^2\mf{d}_0\|_{\mf{L}^2(\Omega)}^2+\|\partial_t\mf{d}^k\|_{\mf{L}^2(\Omega)}^2\\
&\quad+\|\mf{v}\cdot\nabla \mf{d}^{k-1}\|_{\mf{L}^2(\Omega)}^2+\||\nabla
\mf{d}^{k-1}|^2\mf{d}^{k-1}\|_{\mf{L}^2(\Omega)}^2\Big)  \\
\leq &C(K,\|\mf{d}_0\|_{\mf{H}^2(\Omega)}) +c\|\mf{v}\cdot\nabla
\mf{d}^{k-1}\|_{\mf{L}^2(\Omega)}^2 +c\| |\nabla
\mf{d}^{k-1}|^2\mf{d}^{k-1}\|_{\mf{L}^2(\Omega)}^2,\end{aligned}\end{equation}
and
\begin{equation}\label{nestf0429}\begin{aligned}
\|\nabla^3 \mf{d}^k\|_{L^2(I_d^*,\mf{L}^2(\Omega)}^2\leq &
C(K,\|\mf{d}_0\|_{\mf{H}^2(\Omega)},\|\nabla^3
\mf{d}_0\|_{\mf{L}^2(\Omega)}) +c\|\nabla (\mf{v}\cdot\nabla
\mf{d}^{k-1})\|_{L^2(I_d^*,\mf{L}^2(\Omega))}^2 \\
&\quad +c\| \nabla (|\nabla
\mf{d}^{k-1}|^2\mf{d}^{k-1})\|_{L^2(I_d^*,\mf{L}^2(\Omega))}^2,\end{aligned}\end{equation}
where
$$\|\nabla^3\mf{d}^k\|_{L^2(I_d^*,\mf{L}^2(\Omega))}^2:=\sum_{1\leq
i,j,l,m\leq 2}\|\partial_i\partial_j\partial_m{d}_l^k\|_{L^2(I_d^*,{L}^2(\Omega))}^2.$$
To complete the derivation of (\ref{0423}), it suffices to deduce the
uniform estimates of the last two terms in (\ref{0430}) and
(\ref{nestf0429}) as follows.
%Next we start with the  deduction of uniform estimates of the four terms.

%To this purpose, we shall introduce the following
%lemmaApplying this lemma to $\mf{d}$, we have
%\begin{equation*}\begin{aligned}&\|\mf{d}(t)\|_{\mf{H}^2(\Omega)}^2-\|\mf{d}_0\|_{\mf{H}^2(\Omega)}^2\\
%=&2\int_0^t\int_\Omega \left(\partial_t\mf{d}\cdot \mf{d}+\nabla
%\partial_t\mf{d}:\nabla
%\mf{d}+\sum_{1\leq i,j,k\leq 2}\partial_{ij}\partial_t
%d_k\partial_{ij}d_k\right)(\mf{x},\tau)\mm{d}\mf{x}\mm{d}\tau\\
%\leq & c\sup_{\tau\in I_d^*}\|\mf{d}\|_{\mf{H}^2(\Omega)}
%\int_0^t\|\partial_t\mf{d}\|_{\mf{H}^2(\Omega)}\leq
%c\kappa_d^2\sqrt{t}\mbox{ for any }t\in
%I_d^*.\end{aligned}\end{equation*}
%
Since $\mf{d}\in \mathbb{D}^*$, in view of \cite[Lemmas 1.65, 1.66]{NASII04}, we have
\begin{equation*}\begin{aligned}  &\mf{d}(t)-\mf{d}_0
=\int_0^t\partial_t\mf{d}(\tau)\mm{d}\mf{\tau}\mbox{ for any }t\in
I_d^*,\end{aligned}\end{equation*}
which implies that
\begin{equation}\label{jjqwww0339}\begin{aligned}&\|\mf{d}(t)-\mf{d}_0\|_{\mathbf{H}^1(\Omega)}
\leq
\sqrt{t}\|\partial_t\mf{d}\|_{L^2(I_d^*,\mathbf{H}^1(\Omega))}.\end{aligned}\end{equation}
%Similarly, we
%have\begin{equation}\begin{aligned}\label{jjqwww0340}&\| \mf{v}(t)-
%\mf{v}_0\|_{\mathbf{H}^1(\Omega)}\leq
%\sqrt{t}\|\partial_t\mf{v}\|_{H^1(I_d^*,\mathbf{H}^1(\Omega))}.\end{aligned}\end{equation}
%
%\begin{equation*}\begin{aligned}&\|\mf{d}(t)\|_{\mf{H}^2(\Omega)}^2-\|\mf{d}_0\|_{\mf{H}^2(\Omega)}^2\\
%=&2\int_0^t\int_\Omega \left(\partial_t\mf{d}\cdot \mf{d}+\nabla
%\partial_t\mf{d}:\nabla
%\mf{d}+\sum_{1\leq i,j,k\leq 2}\partial_{ij}\partial_t
%d_k\partial_{ij}d_k\right)(\mf{x},\tau)\mm{d}\mf{x}\mm{d}\tau\\
%\leq & c\sup_{\tau\in I_d^*}\|\mf{d}\|_{\mf{H}^2(\Omega)}
%\int_0^t\|\partial_t\mf{d}\|_{\mf{H}^2(\Omega)}\leq
%c\kappa_d^2\sqrt{t}\mbox{ for any }t\in
%I_d^*.\end{aligned}\end{equation*}
 Using (\ref{jjqwww0339}), we obtain
\begin{equation*}\begin{aligned}   \|\mf{v}\cdot\nabla
\mf{d}^{k-1}\|_{\mf{L}^2(\Omega)}^2 \leq&
c\|\mf{v}\|^2_{\mf{L}^\infty(\Omega)}\Big(\|\nabla \mf{d}^{k-1}-\nabla
\mf{d}^{k-1}_0\|_{\mf{L}^2(\Omega)}^2+\|\nabla
\mf{d}^{k-1}_0\|_{\mf{L}^2(\Omega)}^2\Big) \\
\leq &cK^2 ({t}\kappa_d^2+\|\mf{d}_0\|_{\mathbf{H}^1(\Omega)}^2).
\end{aligned}
\end{equation*}
Recalling that $0\leq t\leq T^*_d\leq \kappa_d^{-12}$, we have
\begin{equation}\label{0434}
 \|\mf{v}\cdot\nabla\mf{d}^{k-1}\|_{\mf{L}^2(\Omega)}^2\leq
C(K,\|\mf{d}_0\|_{\mathbf{H}^1(\Omega)}).
\end{equation}

To estimates the last term in (\ref{0430}) and the the last two
terms in (\ref{nestf0429}), we shall introduce the following well-known
Ehrling-Nirenberg-Gagliardo interpolation inequality:
\begin{lem}\label{lem:0201}
Let $\Omega$ be a domain in $\mathbb{R}^N$ satisfying the cone
condition.
%\begin{enumerate}
%  \item [\quad \quad (1)] For each $\varepsilon_0>0$ there is a constant $K'$
%depending on $N$, $m$, $p$ and $\varepsilon_0$, such that if
%$0<\varepsilon\leq \varepsilon_0$, $0\leq j\leq m$, and $u\in
%W^{m,p}(\Omega)$, then
%\begin{equation}\label{0431}
%\sum_{|\alpha|=j}\int_\Omega|D^\alpha u(x)|^p\mathrm{d}x\leq
%K'\left(\varepsilon\sum_{|\alpha|=m}\int_\Omega|D^\alpha
%u(x)|^p\mathrm{d}x +\varepsilon^{-j/(m-j)}\int_\Omega
%|u|^p\mathrm{d}x\right).
%\end{equation}
  %\item [\quad \quad (2)]
  If $mp>N$, let $p\leq q\leq \infty$; if
  $mp=N$, let $p\leq q<\infty$; if $mp<N$, let $p\leq q\leq
  p^*=Np/(N-mp)$. Then there exists a constant ${c}_1$ depending on
  $m$, $N$, $p$, $q$ and the dimension of the cone $\mathbf{C}$,
%%  providing the cone condition for $\Omega$,
  such that for all $v\in W^{m,p}(\Omega)$,
\begin{equation}\label{jiangfeiweiwei0432}
\|v\|_{L^q(\Omega)}\leq {c}_1\|v\|_{W^{m,p}}^\alpha\|v\|_{L^p(\Omega)}^{1-\alpha},
\end{equation}
where $\alpha=(N/mp)-(N/mq)$.
%\end{enumerate}
\end{lem}
\begin{pf} The proof can be found in \cite[Chapter 5]{ARAJJFF}.
\hfill $\Box$
\end{pf}

In particular, we can deduce from
(\ref{jiangfeiweiwei0432}) that
\begin{eqnarray}
\label{jjjw0432}&&\|v\|_{L^\infty(\Omega)}^2\leq {c}_1(\Omega)\|
v\|_{H^{2}(\Omega)}\|v\|_{L^2(\Omega)}\\
&&\label{njjjw0432}\|v\|_{L^4(\Omega)}^2\leq {c}_1(\Omega)\|
v\|_{H^{1}(\Omega)}\|v\|_{L^2(\Omega)}
\end{eqnarray}for some constants ${c}_1(\Omega)$ depending on $\Omega$.
 %This interpolation
 %inequality will be used in the Section \ref{Cauchy}.

Making use of (\ref{jjqwww0339}), (\ref{jjjw0432}) and (\ref{njjjw0432}), it is easy to see that
\begin{equation*}\begin{aligned}
 \||\nabla \mf{d}^{k-1}|^2\mf{d}^{k-1}(t)\|_{\mf{L}^2(\Omega)}^2\leq&
c  \| \nabla\mf{d}^{k-1}\|^4_{\mf{L}^{4}(\Omega)}\|
\mf{d}^{k-1}\|_{\mf{L}^\infty(\Omega)}^2
\\
\leq &c\Big(\| \nabla\mf{d}^{k-1}\|^4_{\mf{L}^{4}(\Omega)}\|
\mf{d}^{k-1}\|_{\mf{L}^\infty(\Omega)}^2+\|\nabla\mf{d}^{k-1}_0\|^4_{\mf{L}^{4}(\Omega)}
\|\mf{d}^{k-1}_0\|_{\mf{L}^\infty(\Omega)}^2  \\
& +\|\nabla\mf{d}^{k-1}-\nabla\mf{d}^{k-1}_0\|^4_{\mf{L}^{4}(\Omega)}\|
\mf{d}^{k-1}_0\|_{\mf{L}^\infty(\Omega)}^2\Big) \\
%%%%%%%
 \leq & C_1(\Omega)\Big(
\kappa^5 \| \mf{d}^{k-1}-\mf{d}_0\|_{\mf{H}^{1}(\Omega)} +\|\mf{d}_0\|_{\mathbf{H}^2(\Omega)}^6\Big)\\
%%%%%
\leq & C_1(\Omega)\Big(\kappa^5
\|\partial_t\mf{d}\|_{L^2(I_d^*,\mathbf{L}^2(\Omega))}\sqrt{t}+\| \mf{d}_0\|^6_{\mf{H}^{2}(\Omega)}\Big)
\end{aligned}\end{equation*}
for any $t\in I_d^*$. Since $T_d^*\leq \kappa_d^{-12}$, we have
\begin{equation}\begin{aligned}\label{0436}
 \||\nabla \mf{d}^{k-1}|^2\mf{d}^{k-1}\|_{\mf{L}^2(\Omega)}\leq
C_1(\Omega)\left(1+\| \mf{d}_0\|_{\mf{H}^{2}(\Omega)}^6\right).
\end{aligned} \end{equation}

Using H\"older's inequality, (\ref{jjjw0432}) and (\ref{njjjw0432}), we infer that
\begin{equation*}\begin{aligned}\label{nest0429}&
\|\nabla (\mf{v}\cdot\nabla
\mf{d}^{k-1})\|_{L^2(I_d^*,\mf{L}^2(\Omega))}^2 +c\| \nabla (|\nabla
\mf{d}^{k-1}|^2\mf{d}^{k-1})\|_{L^2(I_d^*,\mf{L}^2(\Omega))}^2\\
&\leq c\Big(\|\nabla \mf{v}\|_{L^\infty(I_d^*,\mf{L}^2(\Omega))}^2\|\nabla
\mf{d}^{k-1}\|_{L^2(I_d^*,\mf{L}^\infty(\Omega))}^2+\|
\mf{v}\|_{\mf{L}^\infty(Q_T)}^2\|\nabla^2
\mf{d}^{k-1}\|_{L^2(I_d^*,\mf{L}^2(\Omega)}^2\\
&\quad\ + \|\mf{d}^{k-1}\|_{\mf{L}^\infty(Q_T)}^2\||\nabla\mf{d}^{k-1}| \nabla^2
\mf{d}^{k-1}\|_{L^2(I_d^*,\mf{L}^2(\Omega))}^2+ \int_{I_d^*}\|\nabla
\mf{d}^{k-1}\|_{\mf{L}^\infty(\Omega)}^2\|\nabla
\mf{d}^{k-1}\|_{\mf{L}^4(\Omega)}^4\mm{d}s\Big) \\
&\leq C_1(\Omega)\Big(K^2\|\nabla
\mf{d}^{k-1}\|_{L^\infty(I_d^*,\mf{L}^2(\Omega))}\|\mf{d}^{k-1}\|_{L^1(I_d^*,\mf{H}^3(\Omega))}
+K^2\|\nabla \mf{d}^{k-1}\|_{L^\infty(I_d^*,\mf{L}^2(\Omega))}^2T_d^*\\
&\quad \ +\| \mf{d}^{k-1}\|_{L^\infty(I_d^*,\mf{H}^2(\Omega)}^5\|\mf{d}^{k-1}|
\|_{L^1(I_d^*,\mf{H}^3(\Omega))}\Big)\\
&\leq C_1(\Omega)\Big(K^2\kappa^2\sqrt{T_d^*} +K^2\kappa^2T_d^*+
\kappa^6\sqrt{T_d^*}\Big),\end{aligned}\end{equation*}
which, together with the fact $T_d^*\leq \kappa_d^{-12}\leq 1$, yields
\begin{equation}\begin{aligned}\label{nf0436}
\|\nabla (\mf{v}\cdot\nabla
\mf{d}^{k-1})\|_{L^2(I_d^*,\mf{L}^2(\Omega))}^2 +c\| \nabla (|\nabla
\mf{d}^{k-1}|^2\mf{d}^{k-1})\|_{L^2(I_d^*,\mf{L}^2(\Omega))}^2\leq C_1(\Omega)(1+K^2).
\end{aligned} \end{equation}

 Finally, putting (\ref{js0330n}), (\ref{0428})--(\ref{nestf0429}), (\ref{0436})
 and (\ref{nf0436}) together, we then obtain that
\begin{eqnarray*}&&\|\mf{d}^k\|_{\mathbb{D}^*}\leq C(K,\|\mf{d}_0\|_{\mathbf{H}^2(\Omega)},\|\nabla^3
\mf{d}_0\|_{\mf{L}^2(\Omega)}).\end{eqnarray*} Choosing
$\kappa_d\geq
\max\{C(K,\|\mf{d}_0\|_{\mathbf{H}^3(\Omega)},\|\nabla^3
\mf{d}_0\|_{\mf{L}^2(\Omega)}),\|\mf{d}_0\|_{\mathbf{H}^3(\Omega)},1\}$,
we get (\ref{0423}).  %Based on (\ref{0423}), using
%(\ref{jjjw0432ww}), we can further deduce higher uniform estimate
%from (\ref{0420}), i.e., we also have
%\begin{eqnarray}\label{fzn0344}\|\mf{d}^k\|_{L^2(I_d^*,\mf{H}^{3}_2(\Omega)}
% \leq C_3(K,\|\mf{d}_0\|_{\mathbf{H}^3(\Omega)})(1+\|\mf{v}\|_{L^2(I, \mf{H}^2(\Omega))}).\end{eqnarray}

\subsubsection{Taking limits}
In order to take limits in (\ref{0420}), we shall further show that
$\{\mf{d}^k\}_{k=1}^\infty$ is a Cauchy sequence. To this end,  we
define $$\bar{\mf{d}}^{k+1}={\mf{d}}^{k+1}-{\mf{d}}^{k}$$  which satisfies
\begin{eqnarray}\label{0439}
\partial_t\bar{\mf{d}}^{k+1}-\theta\Delta\bar{\mf{d}}^{k+1}
=\theta[\nabla \bar{\mf{d}}^k:(\nabla \mf{d}^k+\nabla
\mf{d}^{k-1})]\mf{d}^{k}+\theta|\nabla
\mf{d}^{k-1}|^2\bar{\mf{d}}^k-\mf{v}\cdot\nabla\bar{\mf{d}}^{k}  \end{eqnarray}
with initial and boundary conditions:
\begin{eqnarray*}
&&\mf{d}^{k+1}|_{t=0}=\mf{0},\quad \mf{x}\in \Omega ,  \\
&& \mf{d}^{k+1}|_{\partial\Omega}=\mf{0},\quad t\in I_d^*.
\end{eqnarray*}

Using Cauchy-Schwarz's and H\"older's inequalities, (\ref{jjjw0432}) and (\ref{njjjw0432})
again, it follows from (\ref{0439}) and the boundary condition that
%\begin{equation*}  \begin{aligned}&\frac{1}{2}\frac{d}{dt}\|\bar{\mf{d}}^{k+1}\|^2_{\mf{L}^2(\Omega)}+\|\nabla
%\bar{\mf{d}}^{k+1}\|^2_{\mf{L}^2(\Omega)}\leq C_4(K,{\kappa}_d)\|
%\bar{\mf{d}}^k\|^2_{\mf{H}^1(\Omega)}+\frac{1}{2}\|
%\bar{\mf{d}}^{k+1}\|^2_{\mf{L}^2(\Omega)}
%\end{aligned}\end{equation*} and
%\begin{eqnarray*}
%\frac{1}{2}\frac{d}{dt}\|\nabla
%\bar{\mf{d}}^{k+1}\|^2_{\mf{L}^2(\Omega)}+\|\Delta
%\bar{\mf{d}}^{k+1}\|^2_{\mf{L}^2(\Omega)} \leq C_4(K,{\kappa}_d)\|
%\bar{\mf{d}}^k\|^2_{\mf{H}^1(\Omega)},\end{eqnarray*} which imply
\begin{equation}\begin{aligned}\label{fz0346}
&\frac{d}{dt}\|\bar{\mf{d}}^{k+1}\|^2_{\mf{H}^1(\Omega)}+2\left(\|\nabla
\bar{\mf{d}}^{k+1}\|^2_{\mf{L}^2(\Omega)}+\|\Delta
\bar{\mf{d}}^{k+1}\|^2_{\mf{L}^2(\Omega)}\right)\\
\leq& c\Big(\|\nabla \bar{\mf{d}}^k\|^2_{\mf{L}^2(\Omega)}\|\nabla(
{\mf{d}}^k+ {\mf{d}}^{k-1})\|^2_{\mf{L}^\infty(\Omega)}\|
{\mf{d}}^{k}\|^2_{\mf{L}^\infty(\Omega)}\\
&\quad  +\|\bar{\mf{d}}^k\|^2_{\mf{L}^4(\Omega)}\|\nabla
{\mf{d}}^{k-1}\|^2_{\mf{L}^\infty(\Omega)}\| \nabla
{\mf{d}}^{k}\|^2_{\mf{L}^4(\Omega)}+\|{\mf{v}}^k\|^2_{\mf{L}^\infty(\Omega)}\|\nabla
\bar{\mf{d}}^k\|^2_{\mf{L}^2(\Omega)}\Big)
+\|\bar{\mf{d}}^{k+1}\|^2_{\mf{H}^1(\Omega)}    \\
\leq& C_1(\Omega)\|\bar{\mf{d}}^k\|^2_{\mf{H}^1(\Omega)}\Big(K^2+\kappa_d^3
(\|\mf{d}^{k}\|_{\mf{H}^3(\Omega)}+\|\mf{d}^{k-1}\|_{\mf{H}^3(\Omega)})\Big)
+\|\bar{\mf{d}}^{k+1}\|^2_{\mf{H}^1(\Omega)}.
%\\
%\leq& C_4(K,{\kappa}_d)\|\bar{\mf{d}}^k\|^2_{\mf{H}^1(\Omega)}(1+\|\mf{d}^{k}\|_{\mf{H}^3(\Omega)}
+\|\mf{d}^{k-1}\|_{\mf{H}^3(\Omega)})
%+\|\bar{\mf{d}}^{k+1}\|^2_{\mf{H}^1(\Omega)}.
\end{aligned} \end{equation}
Then, using Gronwall's and H\"older inequalities, we have
\begin{equation*}\begin{aligned}
\frac{d}{dt}\|\bar{\mf{d}}^{k+1}\|^2_{\mf{H}^1(\Omega)}
\leq & C_1(\Omega)\left(K^2t +\kappa_d^3
\int_0^t(\|\mf{d}^{k}\|_{\mf{H}^3(\Omega)}+\|\mf{d}^{k-1}\|_{\mf{H}^3(\Omega)})\mm{d}s\right)
\sup_{0\leq s\leq t}\|\bar{\mf{d}}^{k}\|^2_{\mf{H}^1(\Omega)} \\
\leq & C(K,{\kappa}_d)(t+\sqrt{t})e^{t}\sup_{0\leq s\leq t}\|\bar{\mf{d}}^{k}\|^2_{\mf{H}^1(\Omega)}.
\end{aligned}\end{equation*}
Taking
$$T^*_3=\min\{T^*_2,C^{-2}(K,{\kappa}_d)/(64e^2),1/2\}$$ and letting $T_d^*\leq T^*_3$,
we have
\begin{eqnarray} \label{fz0347}
\sup_{t\in {I}_d^*}\|\bar{\mf{d}}^{k+1}\|^2_{\mf{H}^1(\Omega)}
\leq \frac{1}{4}\sup_{t\in {I}_d^*}\|\bar{\mf{d}}^{k}\|^2_{\mf{H}^1(\Omega)}.
\end{eqnarray}
Furthermore, from (\ref{fz0346}), (\ref{fz0347}) and the elliptic estimates we get
\begin{eqnarray*}\|\bar{\mf{d}}^{k+1}\|^2_{L^2(I_d^*,\mf{H}^2(\Omega))}
\leq \frac{1}{4}\sup_{t\in
{I}_d^*}\|\bar{\mf{d}}^{k}\|^2_{\mf{H}^1(\Omega)}. \end{eqnarray*}
By iteration, we have
\begin{eqnarray*}
\sup_{t\in {I}_d^*}\|\bar{\mf{d}}^{k+1}\|^2_{\mf{H}^1(\Omega)}
+\|\bar{\mf{d}}^{k+1}\|^2_{L^2(I_d^*,\mf{H}^2(\Omega))}\leq\frac{1}{2^{k-1}}
\sup_{t\in {I}_d^*}\|\bar{\mf{d}}^{k}\|^2_{\mf{H}^1(\Omega)}\; \mbox{ for any
}k\geq 2.\end{eqnarray*}
In particular,
\begin{eqnarray*}\label{n0438}
\sum_{k=2}^\infty\Big(\|\bar{\mf{d}}^k\|_{L^\infty(I_d^*, \mathbf{H}^1(\Omega))}+
\|\bar{\mf{d}}^k\|_{L^2(I_d^*, \mathbf{H}^2(\Omega))}\Big)<\infty,\end{eqnarray*}
which implies
$\{\mf{d}_k\}_{k=1}^\infty$ is a Cauchy sequence in
${L^\infty(I_d^*, \mathbf{H}^1(\Omega))}\cap {L^2(I_d^*,
\mathbf{H}^2(\Omega))}$. Therefore,
\begin{eqnarray}\label{0447}\mf{d}^k\rightarrow {\mf{d}}
\mbox{ strongly in }{L^\infty(I_d^*, \mathbf{H}^1(\Omega))}\cap
L^2(I_d^*,\mf{H}^2(\Omega)).\end{eqnarray}

As a consequence, $\mf{d}$ is accurately a solution to
problem (\ref{0411})--(\ref{0413}) satisfying the regularity
$\|\mf{d}\|_{\mathbb{D}}\leq \kappa_d$, by virtue of (\ref{0423})
and (\ref{0447}), the lower semi-continuity, the elliptic estimate and the
compactness theorem with time (see \cite[Lemma 2.5]{DSJWHYSN}).
% and the following lemma (see \cite[Theorem 1.67]{NASII04}):
%\begin{lem}Let $H$ be a Hilbert space and $V\hookrightarrow H$
%be dense in $H$. If $u$, $v\in L^p(I,V)$ with $1<p<\infty$, and
%$u'$, $v'\in L^q(I, V^*)$, $\frac{1}{p}+\frac{1}{q}=1$, then $u$,
%$v\in C(\bar{I}, H)$ and
%\begin{eqnarray*}(u(t),v(t))-(u(s),v(s))=\int_s^t
%(<u'(\tau),v(\tau)>+<v'(\tau),u(\tau)>)\mm{d} \tau.\end{eqnarray*} (
%Here $\hookrightarrow$ and $<\cdot,\cdot>$ denote the continuous
%embedding and the duality between $V$ and $V^*$, respectively.)
%\end{lem}

\subsubsection{Continuous dependence}
Finally we show that the solution $\mf{d}$ continuously depends on
$\mf{v}$. Let $T_d^*\in (0, T^*_3]$, $\mathbf{v}_1$, $\mathbf{v}_2\in\mf{V}_K$, $\mathbf{d}_1$
and $\mathbf{d}_2$ be two solutions of the initial-boundary value problem
(\ref{0411})--(\ref{0413}), corresponding to $\mathbf{v}=\mathbf{v}_1$ and
$\mathbf{v}=\mathbf{v}_2$ respectively. Moreover, two solutions satisfy
\begin{eqnarray*}\label{0448}
\|\mf{d}_1\|_{\mathbb{D}^*} +\|\mf{d}_2\|_{\mathbb{D}^*}<\infty.\end{eqnarray*}

 Multiplying the difference of the equations for $\mathbf{d}_1-\mathbf{d}_2$ by
 $\mathbf{d}_1-\mathbf{d}_2$ and integrating by parts, we obtain
\begin{eqnarray*}\begin{aligned}
&\displaystyle\frac{d}{dt}\int_\Omega|\mathbf{d}_1-\mathbf{d}_2|^2\mathrm{d}\mathbf{x}+2\theta
\int_{\Omega}|\nabla(\mathbf{d}_1-\mathbf{d}_2)|^2\mathrm{d}\mathbf{x}\qquad \qquad \\
&=2\theta\left(\displaystyle\int_\Omega|\nabla
\mf{d}_1|^2|\mathbf{d}_1-\mathbf{d}_2|^2\mathrm{d}\mf{x}+\int_\Omega(\nabla
(\mf{d}_1-\mf{d}_2):\nabla
(\mf{d}_1+\mf{d}_2))\mathbf{d}_2\cdot(\mathbf{d}_1-\mathbf{d}_2)\mm{d}\mf{x}\right)\\
&\quad
-2\displaystyle\int_{\Omega}(\mathbf{v}_1-\mathbf{v}_2)\cdot\nabla
\mathbf{d}_1\cdot(\mathbf{d}_1-\mathbf{d}_2)\mathrm{d}\mathbf{x}
-2\displaystyle\int_{\Omega}\mathbf{v}_2\cdot
\nabla(\mathbf{d}_1-\mathbf{d}_2)\cdot(\mathbf{d}_1-\mathbf{d}_2)\mathrm{d}\mathbf{x},
\end{aligned}\end{eqnarray*}
Similarly to that in the derivation of (\ref{fz0346}), we use (\ref{jjjw0432}) and (\ref{njjjw0432})
to see that the right-hand side of the above identity can be bounded from above by
\begin{equation*}\begin{aligned}
&C(\|\mf{d}_1\|_{\mathbb{D}^*},\|\mf{d}_2\|_{\mathbb{D}^*},K)
\left(1+\|\mf{d}_1\|_{\mf{H}^3(\Omega)}+\|\mf{d}_2\|_{\mf{H}^3(\Omega)}\right)\\
&\qquad \times \left(\| \mf{d}_1-\mf{d}_2\|_{\mathbf{L}^2(\Omega)}^2
+\|\mathbf{v}_1-\mathbf{v}_2\|_{\mf{L}^2(\Omega)}\|
\mf{d}_1-\mf{d}_2\|_{\mf{L}^2(\Omega)}\right)
 +\theta\|\nabla (\mf{d}_1-\mf{d}_2)\|_{\mf{L}^2(\Omega)}^2,\end{aligned}
\end{equation*}
Hence,
\begin{equation*}\begin{aligned}
&\frac{d}{dt}\| \mf{d}_1-\mf{d}_2\|_{\mf{L}^2(\Omega)}^2+\theta\|
\nabla (\mf{d}_1-\mf{d}_2)\|_{\mf{L}^2(\Omega)}^2 \\
&\leq
C(\|\mf{d}_1\|_{\mathbb{D}^*},\|\mf{d}_2\|_{\mathbb{D}^*},K)\left(1+\|\mf{d}_1\|_{\mf{H}^3(\Omega)}
+\|\mf{d}_2\|_{\mf{H}^3(\Omega)}\right) \\
&\qquad \times \left(\|
\mf{d}_1-\mf{d}_2\|_{\mathbf{L}^2(\Omega)}^2+\|\mathbf{v}_1-\mathbf{v}_2\|_{\mf{L}^2(\Omega)}\|
\mf{d}_1-\mf{d}_2\|_{\mf{L}^2(\Omega)}\right),\end{aligned}
\end{equation*}
in particular,
\begin{equation*}\begin{aligned} \frac{d}{dt}\| \mf{d}_1-\mf{d}_2\|_{\mf{L}^2(\Omega)}
\leq&
C(\|\mf{d}_1\|_{\mathbb{D}^*},\|\mf{d}_2\|_{\mathbb{D}^*},K)
\left(1+\|\mf{d}_1\|_{\mf{H}^3(\Omega)}+\|\mf{d}_2\|_{\mf{H}^3(\Omega)}\right)
\\&\quad \times \left(\| \mf{d}_1-\mf{d}_2\|_{\mathbf{L}^2(\Omega)} +
\|\mathbf{v}_1-\mathbf{v}_2\|_{\mf{L}^2(\Omega)}\right).\end{aligned}
\end{equation*}
Thus, applying Gronwall's and Cauchy-Schwarz's inequalities, we conclude
\begin{equation}\begin{aligned}\label{n0461}
\| (\mf{d}_1-\mf{d}_2)(t)\|_{\mf{L}^2(\Omega)}  \leq
(t+\sqrt{t})C(\|\mf{d}_1\|_{\mathbb{D}^*},\|\mf{d}_2\|_{\mathbb{D}^*},K,T_d^*)
\|\mathbf{v}_1-\mathbf{v}_2\|_{L^\infty(I_t,\mf{L}^2(\Omega))}\end{aligned}\end{equation}
 for any $t\in I_d^*$.
Moreover,
\begin{equation*}\begin{aligned}
\|\nabla(\mathbf{d}_1-\mathbf{d}_2)\|_{\mf{L}^2(Q_t)}\leq
(t+\sqrt{t})C(\|\mf{d}_1\|_{\mathbb{D}^*},\|\mf{d}_2\|_{\mathbb{D}^*}
,K,T_d^*)\|\mathbf{v}_1-\mathbf{v}_2\|_{L^\infty(I_t,\mf{L}^2(\Omega))}.
\end{aligned}\end{equation*}

Obviously, the uniqueness of local solutions in the function class
$\mathbb{D}^*$ follows from (\ref{n0461}) immediately. In
particular, $\|\mf{d}_1\|_{\mathbb{D}^*}$, $\|\mf{d}_2\|_{\mathbb{D}^*}\leq \kappa_d$.
 As the end of this subsection, we summarize our previous results on the local existence of $\mathbf{d}$.
 %%%%%%%%%%%%%%%%%%%%%%%%%%%%%%%%%%%
\begin{pro}\label{pro:0402}
Let $K>0$, $0<\alpha<1$, $\Omega$ be a bounded domain of class
$C^{2,\alpha}$, $\mf{v}\in \mathbb{V}_K$, and $\mf{d}_0\in \mf{H}^3(\Omega)$.
Then there exist a finite time
\begin{equation*}{T_d^K}:=h_1(K,\|\mf{d}_0\|_{\mathbf{H}^2(\Omega)},\|\nabla^3
\mf{d}_0\|_{\mf{L}^2(\Omega)})\in (0, \min\{1,T\}),\end{equation*}
and a corresponding unique mapping
\begin{equation*}\mathscr{D}^K_{\mathbf{d}_0}:
\mathbb{V}_K (\mbox{with }T^K_d\mbox{ in place of }T)\rightarrow\
C^0(\bar{I}^K_d,\mf{H}^2(\Omega)),\end{equation*}
where $h_1$ is nonincreasing in its first two variables and
${Q}_{T^K_d}:={\Omega}\times I^K_d:=\Omega\times (0,T^K_d)$, such that
\begin{enumerate}
  \item[\qquad (1)] $\mathscr{D}^K_{\mathbf{d}_0}(\mathbf{v})$ belongs to the following function class
\begin{equation}\label{0465}\begin{aligned}\mathcal{R}_{T_d^K}:=\{\mf{d}~|~%|\mf{d}|\equiv 1
& \mf{d}\in L^2(I_d^K,\mf{H}^3(\Omega)),\ \partial_{tt}^2\mf{d}\in
{L^2(I_d^K,(\mf{H}^{1}(\Omega))^*)},\\
& \partial_t\mf{d}\in {C^0(\bar{I}_d^K,\mf{L}^2(\Omega))}\cap
{L^2(I_d^K,\mf{H}^1_0(\Omega))} \}.\end{aligned}\end{equation}
  \item[\qquad (2)]  $\mf{d}=\mathscr{D}^K_{\mathbf{d}_0}(\mathbf{v})$ satisfies
  (\ref{0411}) a.e. in $Q_{T_d^K}$, (\ref{0412}) in $\Omega$ and (\ref{0413})
 in $I_d^K$. Moreover, $\mf{d}|_{t=0}=\mf{d}_0$.
   \item[\qquad (3)]  $\mathscr{D}^K_{\mathbf{d}_0}(\mathbf{v})$ enjoys the following estimate:
\begin{equation}\label{0455}
\|\mathscr{D}^K_{\mathbf{d}_0}(\mathbf{v})\|_{\mathbb{D}^*}\leq
C(K,\|\mf{d}_0\|_{\mathbf{H}^2(\Omega)},\|\nabla^3
\mf{d}_0\|_{\mf{L}^2(\Omega)}),
\end{equation}where $C$ is
nondecreasing in its variables, and $T_d^*$ in the definition of
$\mathbb{D}^*$ should be replaced by $T^K_{d}$. Moreover (in view of
(\ref{0201}) and (\ref{ji0615})),
\begin{eqnarray}
&&|\mf{d}(t)|\equiv 1\;\mbox{ for any }t,\;\;\mbox{ if } |\mf{d}_0|\equiv
1;\nonumber\\
&&  d_2(t)\geq \underline{d}_{02}\;\mbox{ for any }t,\;\;\mbox{ if }{d}_{02}
\geq \underline{d}_{02}.\nonumber
\end{eqnarray}
  \item[\qquad (4)]   $\mathscr{D}^K_{\mathbf{d}_0}(\mathbf{v})$ continuously depends on
  $\mf{v}$ in the following sense:
  \begin{equation}\label{n0456}
  \begin{aligned}&
\|[\mathscr{D}^K_{\mathbf{d}_0}(\mathbf{v}_1)-\mathscr{D}^K_{\mathbf{d}_0}
(\mathbf{v}_2)](t)\|_{\mathbf{L}^2(\Omega)} +
\|\nabla(\mathscr{D}^K_{\mathbf{d}_0}(\mathbf{v}_1)-\mathscr{D}^K_{\mathbf{d}_0}
(\mathbf{v}_2))\|_{\mf{L}^2(Q_{t})}
\\
&\leq \sqrt{t}C(K,\|\mf{d}_0\|_{\mathbf{H}^2(\Omega)},\|\nabla^3
\mf{d}_0\|_{\mf{L}^2(\Omega)})
\|\mathbf{v}_1-\mathbf{v}_2\|_{L^\infty(I_t,\mf{L}^2(\Omega))}\end{aligned}\end{equation}
for any $t\in \bar{I}_d^K$, where $C$ is nondecreasing in its
variables.
\end{enumerate}
\end{pro}
\begin{rem}\label{rem:0301} We can further show that $\mathscr{D}^K_{\mathbf{d}_0}(\mathbf{v})$
continuously depends on  $\mf{v}$ as the form of (\ref{n04091}). However, the estimate (\ref{n0456})
  is sufficient to prove the local existence of solutions to the third approximate problem.
It should be noted that the constants $h_1$ and $C$ above depend on the domain $\Omega$.
However, in Section \ref{Cauchy}, we will see that the above result in bounded domains
can be generalized to the Cauchy problem.
\end{rem}

\section{Unique solvability of the third approximate problem}\label{sec:05}

In this section, we establish the global existence of a unique solution
to the third approximate problem \eqref{n301}--\eqref{n0306}. We first use the iteration technique
and fixed point theorem to show the local existence. For this
purpose, with the help of Proposition \ref{pro:0401} and \ref{pro:0402},
we shall introduce the operator form of the approximate
momentum equations (\ref{n303}) to construct a contractive mapping.

 \subsection{Operator form}
 Given
\begin{equation}\begin{aligned}\label{0504}\rho\in C^0(\bar{I},L^1(\Omega)),
\quad \partial_t \rho\in L^1(\Omega_{T}),\quad  \underset{(\mathbf{x},t)\in \Omega_T}
{\mathrm{ess\ inf}}\,\rho (\mathbf{x},t)\geq \underline{\rho}>0, \end{aligned}\end{equation}
we define, for all $t\in \bar{I}$
\begin{equation*}\begin{aligned}
\mathscr{M}_{\rho(t)}:\; \mathbf{X}_n \to \mathbf{X}_n
\end{aligned}\end{equation*}
by
$$\ <\mathscr{M}_{\rho(t)}\mathbf{v},\,
\mathbf{w}>\equiv\int_{\Omega}\rho(t)\mathbf{v}\cdot\mathbf{w}\mathrm{d}x,\quad
\mathbf{v},\mathbf{w} \in \mathbf{X}_n.$$
Recall that any norms of $\mathbf{X}_n$ are equivalent, in particular,
\begin{equation}\label{0506}
\mathbf{W}^{k_1,p_1}(\Omega)\mbox{ and }
(\mathbf{W}^{k_2,p_2}_0(\Omega))^*\mbox{-norms are equivalent on }\mathbf{X}_n,\end{equation}
where $(\mathbf{W}^{k_2,p_2}_0(\Omega))^*$ denotes the dual space of
$(\mathbf{W}^{k_2,p_2}_0(\Omega))$, $k_1$ and $k_2$ are integers, and
$0\leq k_2<\infty$, $1\leq p_2\leq \infty$, $0\leq k_1\leq1$
and $1\leq p_1\leq \infty$ (or $k_1=2$, $1\leq p_1<\infty$).
Note that this property will be repeatedly used in the estimates below. First one has
\begin{equation}\begin{aligned}\label{0507}
\|\mathscr{M}_{\rho(t)}\|_{\mathscr{L}(\mathbf{X}_n,\mathbf{X}_n)}\leq c(n)\int_\Omega\rho(t)
\mathrm{d}x,\quad t\in \bar{I}.\end{aligned}\end{equation}
It is easy to observe that $\mathscr{M}^{-1}_{\rho(t)}$ exists for all
$t\in \bar{I}$ and
\begin{equation}\label{0508}
\|\mathscr{M}^{-1}_{\rho(t)}\|_{\mathscr{L}(\mathbf{X}_n,\mathbf{X}_n)}
\leq \frac{1}{\underline{\rho}},\end{equation}
where $\mathscr{L}(\mathbf{X}_n,\mathbf{X}_n)$ denotes the set of
all continuous linear operators mapping $\mathbf{X}_n$ to $\mathbf{X}_n$.

 By virtue of (\ref{0507}) and (\ref{0508}), we have
\begin{equation}\begin{aligned}\label{0509}\|\mathscr{M}^{-1}_{\rho(t)}\mathscr{M}_{\rho_1(t)}
\mathscr{M}^{-1}_{\rho(t)} \|_{\mathscr{L}(\mathbf{X}_n,\mathbf{X}_n)}\leq\frac{c(n)}
{\underline{\rho}^2}\|\rho_1(t)\|_{L^1(\Omega)},\quad t\in
\bar{I}.\end{aligned}\end{equation}
For the difference, the following inequality
\begin{equation}\label{n0510}\begin{aligned}\|\mathscr{M}_{\rho_1(t)}-\mathscr{M}_{\rho_2(t)}
\|_{\mathscr{L}(\mathbf{X}_n,\mathbf{X}_n)}\leq c(n)\|(\rho_2-\rho_1)(t)\|_{L^1(\Omega)},
\quad t\in \bar{I},\end{aligned}\end{equation}
holds. Due to the identity
$\mathscr{M}^{-1}_{\rho_2(t)}-\mathscr{M}^{-1}_{\rho_1(t)}
=\mathscr{M}^{-1}_{\rho_2(t)}(\mathscr{M}_{\rho_1(t)}-\mathscr{M}_{\rho_2(t)})
\mathscr{M}^{-1}_{\rho_1(t)}$,
we find that
\begin{equation}\begin{aligned}\label{0511}\|\mathscr{M}^{-1}_{\rho_2(t)}
-\mathscr{M}^{-1}_{\rho_1(t)}\|_{\mathscr{L}(\mathbf{X}_n,\mathbf{X}_n)}\leq
\frac{c(n)}{\underline{\rho}^2}\|(\rho_2 -\rho_1)(t)\|_{L^1(\Omega)},
\quad t\in \bar{I} \end{aligned}\end{equation}
for $\rho_1(t),\rho_2(t)$ satisfying (\ref{0504}).

 Next, we shall look for $T_n^*\subset (0,T^{\tilde{K}}_d]$ and
\begin{equation*}
\mathbf{v}\in \mathbb{A}:=\{\mf{v}\in C(\bar{I}_n^*,
\mathbf{X}_n)~|~\partial_t\mf{v}\in L^2(I_n^*,\mf{X}_n)\}, \quad
I_n^*:=(0,T_n^*)\subset (0,T^{\tilde{K}}_d)\end{equation*}
with
$\|\mathbf{v}\|_{C(\bar{I}_n^*,\mathbf{H}^2(\Omega))}
+\|\partial_t\mathbf{v}\|_{L^2({I}_n^*, \mathbf{H}^1(\Omega))}\leq \tilde{K}$
for some $\tilde{K}$, satisfying
\begin{equation}\begin{aligned}\label{0513}
&\int_\Omega(\rho\mathbf{v})(t) \cdot\mathbf{\mathbf{\Psi}}\mathrm{d}\mathbf{x}
-\int_{\Omega}\mathbf{m}_0\cdot\mathbf{\mathbf{\Psi}}\mathrm{d}\mathbf{x} \\
&= \int_0^t\int_\Omega\bigg[\mu\Delta\mathbf{v}+(\mu+\lambda)\nabla
\mm{div}\mf{v}-A\nabla \rho^\gamma -\delta \nabla
\rho^\beta-\varepsilon(\nabla\rho\cdot\nabla \mathbf{v}) \\
& \qquad\qquad -\mathrm{div}(\rho\mathbf{v}\otimes\mathbf{v})
-\nu\mathrm{div}\left(\nabla \mathbf{d}\otimes\nabla\mathbf{d}-\frac{|\nabla
\mf{d}|^2\mathbb{I}}{2}\right)\bigg]\cdot\mathbf{\Psi}\mathrm{d}\mathbf{x}\mathrm{d}s
\end{aligned}\end{equation}
for all $t\in [0,T_n]$ and any $\mathbf{\Psi}\in\mathbf{X}_n$,
where $\rho(t)=[\mathscr{S}_{\rho_0}(\mathbf{v})](t)$ is the solution of the
problem (\ref{0401})--(\ref{0403}) constructed in Proposition \ref{pro:0401},
$\mathbf{d}(t)=\mathscr{D}_{\mf{d}_0}^{\tilde{K}}(\mathbf{v})(t)$
is the solution of the problem (\ref{0411})--(\ref{0413}) constructed in
Proposition \ref{pro:0402}. By the regularity of
$(\rho,\mf{d})$ in Propositions \ref{pro:0401}, \ref{pro:0402} and
 the operator $\mathscr{M}_{\rho(t)}$, the equations (\ref{0513}) can be rephrased as
\begin{equation}\label{n0512}
 {\mathbf v}(t) =\mathscr{M}^{-1}_{[\mathscr{S}_{\rho_0}({\mf{
v}})](t)}\bigg((\mathscr{P}\mathbf{m}_0+\int^t_0
\mathscr{P}[{\mathscr{N}(\mathscr{S}_{\rho_0}({\bf v})},{\mf{ v}},
\mathscr{D}_{\mf{d}_0}^{\tilde{K}}(\mathbf{v}))]\mathrm{d}s\bigg)
\end{equation}
with $\mf{m}_0=(\rho\mf{v})(0)$, where $\mathscr{P}:=\mathscr{P}_n$
is the orthogonal projection of $\mathbf{L}^2(\Omega)$ to $\mathbf{X}_n$, and
\begin{equation}\label{0515}\begin{aligned}
&\mathscr{N}(\rho,\mathbf{v},\mathbf{d})=\mu\Delta\mathbf{v}+(\mu+\lambda)\nabla
\mm{div}\mf{v}-A\nabla \rho^\gamma -\delta \nabla
\rho^\beta-\varepsilon(\nabla\rho\cdot\nabla \mathbf{v})\\
&\quad\quad\quad\quad\quad -\mathrm{div}(\rho\mathbf{v}\otimes\mathbf{v})
-\nu\mathrm{div}\left(\nabla \mathbf{d}\otimes\nabla
\mathbf{d}-\frac{|\nabla \mf{d}|^2\mathbb{I}}{2}\right).
\end{aligned}
\end{equation}
Moreover, one has
\begin{equation}\label{n0413nn}\begin{aligned}
\partial_t\mathbf{v}(t) & =\mathscr{M}_{[\mathscr{S}_{\rho_0}({\mf{v}})](t)}^{-1}
\mathscr{M}_{\partial_t[\mathscr{S}_{\rho_0}({\mf{v}})](t)}
\mathscr{M}_{[\mathscr{S}_{\rho_0}({\mf{v}})](t)}^{-1}
\bigg\{\mathscr{P}\mathbf{m}_0 +\int_0^t[\mathscr{P}\mathscr{N}(\mathscr{S}_{\rho_0}(\mathbf{v}),
\mathbf{v},\mathscr{D}_{\mf{d}_0}(\mf{v}))](s)\mathrm{d}s\bigg\} \\
&\quad +\mathscr{M}_{[\mathscr{S}_{\rho_0}({\mf{ v}})](t)}^{-1}[\mathscr{P}
\mathscr{N}(\mathscr{S}_{\rho_0}(\mathbf{v}),\mathbf{v},\mathscr{D}^{\tilde{K}}_{\mf{d}_0}(\mf{v}))](t).
\end{aligned}\end{equation}

\subsection{Auxiliary estimates}
We shall derive some auxiliary estimates on
$(\rho=\mathscr{S}_{\rho_0}(\mathbf{v}),\mathbf{v},\mf{d}=\mathscr{D}_{\mf{d}_0}^{\tilde{K}}(\mathbf{v}_n))$
and
$(\rho_k=\mathscr{S}_{\rho_0}(\mathbf{v}_k),\mathbf{v}_k,\mf{d}=\mathscr{D}_{\mf{d}_0}^{\tilde{K}}(\mathbf{v}_k))$,
$k=1,2$, where $\mathbf{v}$ and $\mf{v}_k$ belong to the class $\mathbb{A}$, and
 $$\|\mathbf{v}\|_{C(\bar{I}_n^*,\mathbf{X}_n)}+\|\partial_t\mathbf{v}\|_{L^2({I}_n^*,\mathbf{X}_n)}\leq K, \quad
 \|\mathbf{v}_k\|_{C(\bar{I}_n^*,\mathbf{X}_n)}+\|\partial_t\mathbf{v}_k\|_{L^2({I}_n^*,\mathbf{X}_n)} \leq K,$$
 with $K$ being a positive constant. By the equivalence of norms in (\ref{0506}), we have
\begin{equation*}
\left(\|\mathbf{v}\|_{C(\bar{I}_n^*,\mathbf{H}^1(\Omega))}^2
+\|\mathbf{v}\|_{C(\bar{I}_n^*,\mathbf{L}^\infty(\Omega))}^2 +
+\|\partial_t\mathbf{v}\|_{L^2({I}_n^*,\mathbf{H}^1(\Omega))}^2\right)^\frac12
\leq c(n)K=:\tilde{K}.
\end{equation*}
We denote $(\mathscr{S}(\mathbf{v})$,
$\mathscr{D}(\mathbf{v})) = (\mathscr{S}_{\rho_0}(\mathbf{v})$,
$\mathscr{D}_{\mf{d}_0}^{\tilde{K}}(\mathbf{v}))$ for simplicity.

From (\ref{0515}) and (\ref{0506}) we get
\begin{equation}\label{0516}\begin{aligned}
\|\mathscr{P}\mathscr{N}(\rho,\mathbf{v},\mathbf{d})\|_{\mathbf{X}_n}
\leq & c(n)\Big(\|\mathbf{v}\|_{\mathbf{X}_n}
+\|\rho\|_{L^\infty(\Omega)}(\|\mathbf{v}\|_{\mathbf{X}_n}
+\|\mathbf{v}\|_{\mathbf{X}_n}^2 ) \\
& \quad\quad\quad\  +\| \rho\|_{L^\infty(\Omega)}^\gamma+\|
\rho\|_{L^\infty(\Omega)}^\beta +\|\nabla \mathbf{d}\|_{\mathbf{L}^2(\Omega)}^2\Big).
 \end{aligned}
\end{equation}
From (\ref{0404}), (\ref{0415}), (\ref{0455}) and (\ref{0516}), it follows that
\begin{equation}\begin{aligned}\label{0517}
&\|\mathscr{P}\mathscr{N}(\mathscr{S}(\mathbf{v}),
\mathbf{v},\mathscr{D}(\mathbf{v}))(t)\|_{\mathbf{X}_n} \leq
h_2(K,\bar{\rho},\|\mf{d}_0\|_{\mf{H}^2(\Omega)},\|\nabla^3
\mf{d}_0\|_{\mf{L}^2(\Omega)},T,n),\quad t\in I_n^{*},
 \end{aligned}
\end{equation}
where the constant $h_2$ is nondecreasing in its first three
variables.

Employing (\ref{0515}) and the formula
$F(z_1)-F(z_2)=\int_{z_1}^{z_2}F'(s)\mm{d}s$ with $F(s)=s^\gamma$
and $F(s)=s^\beta$, we obtain
\begin{equation}\begin{aligned}\label{0518}
&<\mathscr{N}(\rho_1,\mathbf{v}_1,
\mathbf{d}_1)-\mathscr{N}(\rho_2,\mathbf{v}_2,\mathbf{d}_2),\mathbf{\Phi}>\\
&=\int_\Omega
[\mu\Delta(\mathbf{v}_1-\mathbf{v}_2)+(\mu+\lambda)\nabla
\mm{div}(\mf{v}_1-\mf{v}_2)]\cdot\mathbf{\Phi}\mathrm{d}\mathbf{x}
\\
& \quad + \int_\Omega[(\rho_1-\rho_2)u_1^iu_1^j+
\rho_2(u_1^i-u_2^i\big)u_1^j+
\rho_2u_2^i(u_1^j-u_2^j\big)]\partial_j\Phi^i\mathrm{d}\mathbf{x}\\
&\quad
+\varepsilon\int_\Omega\rho_2[\Delta(\mathbf{v}_1-\mathbf{v}_2)\cdot\mathbf{\Phi}
+\partial_j(u_1^i-u_2^i)\partial_j\Phi^i]\mathrm{d}\mathbf{x}+
\varepsilon\int_\Omega\big(\rho_1-\rho_2\big)(\Delta
\mathbf{v}_1\cdot
\mathbf{\Phi}+\partial_ju_1^i\partial_j\Phi^i)\mathrm{d}\mathbf{x}\\
&\quad+\nu\int_\Omega\partial_id^k_2(
\partial_jd^k_1-\partial_jd^k_2)\partial_j\Phi^i\mathrm{d}\mathbf{x}+
\nu\int_\Omega\partial_jd^k_1(\partial_id^k_1-\partial_id^k_2))\partial_j\Phi^i\mathrm{d}\mathbf{x}\\
&\quad+\frac{\nu}{2}\int_\Omega(\partial_jd^k_1+\partial_jd^k_2)
(\partial_jd^k_2-\partial_jd^k_1)\partial_i\Phi^i\mathrm{d}\mathbf{x}+\int_\Omega\left[\int_{\rho_2}^{\rho_1}(\gamma
A s^{\gamma-1}+\delta\beta
s^{\beta-1})\mm{d}s\right]\mm{div}\mf{\Phi}\mm{d}\mf{x}
 \end{aligned}
\end{equation} for any $\Phi\in \mf{W}_0^{1,\infty}(\Omega)$.
 Making use of (\ref{0404}), (\ref{0455}), (\ref{0506}), (\ref{0518}),
 and the elementary properties of the projection
$\mathscr{P}:=\mathscr{P}_n$ (see \cite[Exercise 7.33]{NASII04}), we get
\begin{equation}\begin{aligned}\label{0520}
&\|[\mathscr{P}\mathscr{N}(\rho_1,\mathbf{v}_1,\mathbf{d}_1)
-\mathscr{P}\mathscr{N}(\rho_2,\mathbf{v}_2,\mathbf{d}_2)](t)\|_{\mathbf{X}_n}\\
&\quad \leq h_2(K,\bar{\rho},\|\mathbf{d}_0\|_{\mathbf{H}^2(\Omega)},\|\nabla^3
\mf{d}_0\|_{\mf{L}^2(\Omega)}, T,n)\{\|
(\mathbf{v}_1-\mathbf{v}_2)(t)\|_{\mathbf{X}_n}   \\
&\qquad +\|(\rho_1-\rho_2)(t)\|_{L^1(\Omega)}+
 \|\nabla (\mathbf{d}_1-\mathbf{d}_2)(t)\|_{\mathbf{L}^2(\Omega)}\}, \quad\; t\in I_n^*,
 \end{aligned}  \end{equation}
where $h_2$ is again nondecreasing in its first three variables.

Thanks to (\ref{0511}), (\ref{0409}), (\ref{0404}), and H\"older's inequality, we obtain
\begin{equation}\label{n0521}
\|\mathscr{M}^{-1}_{[\mathscr{S}(\mathbf{v}_1)](t)}-
\mathscr{M}^{-1}_{[\mathscr{S}(\mathbf{v}_2)](t)}\|_{\mathscr{L}(\mf{X}_n,\mf{X}_n)}\leq
h_3(K,\|\rho_0\|_{H^1(\Omega)},\underline{\rho},T,n)
t\|\mathbf{v}_1-\mathbf{v}_2\|_{C^0(\bar{I_t},\mathbf{X}_n)}, \end{equation}
for any $\mathbf{v}_1$, $\mathbf{v}_2\in C^0(\bar{I}_n^*,\mathbf{X}_n)$, and $t\in I_n^*$.
Here the constant $h_3$ is nondecreasing in its first two variables, and is
nonincreasing in its third variable. Since
\begin{equation*}\begin{aligned}&
\mathscr{M}^{-1}_{[\mathscr{S}(\mathbf{v})](t_1)}\int_0^{t_1}[
\mathscr{P}\mathscr{N}(\mathscr{S}(\mathbf{v}),\mathbf{v},\mathscr{D}(\mathbf{v}))](s)\mathrm{d}s-
\mathscr{M}^{-1}_{[\mathscr{S}(\mathbf{v})](t_2)}\int_0^{t_2}[\mathscr{P}\mathscr{N}(\mathscr{S}(\mathbf{v}),
\mathbf{v},\mathscr{D}(\mathbf{v}))](s)\mathrm{d}s\\
&=\mathscr{M}^{-1}_{[\mathscr{S}(\mathbf{v})](t_1)}\int_{t_2}^{t_1}[\mathscr{P}(\mathbf{v})
\mathscr{N}(\mathscr{S}(\mathbf{v}),\mathbf{v},\mathscr{D}(\mathbf{v}))](s)\mathrm{d}s\\
&\quad +\left(\mathscr{M}^{-1}_{[\mathscr{S}(\mathbf{v})](t_1)}
-\mathscr{M}^{-1}_{[\mathscr{S}(\mathbf{v})](t_2)}\right) \int_{0}^{t_2}
[\mathscr{P}\mathscr{N}(\mathscr{S}(\mathbf{v}),\mathbf{v},\mathscr{D}(\mathbf{v}))](s)\mathrm{d}s,
\end{aligned} \end{equation*}
we make use of (\ref{0404}), (\ref{0508}), (\ref{0511}) and (\ref{0517}) to find that
$$\mathscr{M}^{-1}_{\mathscr{S}(\mathbf{v})(t)}\int_0^{t}[\mathscr{P}\mathscr{N}(\mathscr{S}(\mathbf{v}),
\mathbf{v},\mathscr{D}(\mathbf{v}))](s)\mathrm{d}s\in C^0(\bar{I}_n^*,\mathbf{X}_n).$$
 And notice that
\begin{equation}\label{0522}
\bigg\|\mathscr{M}^{-1}_{\mathscr{S}(\mathbf{v})}\int_0^t
[\mathscr{P}\mathscr{N}(\mathscr{S}(\mathbf{v}),
\mathbf{v},\mathbf{d}](s)\mathrm{d}s\bigg\|_{C(I_t,\mathbf{X}_n)}
\leq h_4(K,\bar{\rho},\|\mathbf{d}_0\|_{\mathbf{H}^2(\Omega)},\underline{\rho},
\|\nabla^3 \mf{d}_0\|_{\mf{L}^2(\Omega)} ,T,n)t,\end{equation}
for any $t\in I_n^*$, where $h_4$ is nondecreasing in its first three variables
and nonincreasing in its fourth variable. Similarly,
$\mathscr{M}^{-1}_{[\mathscr{S}(\mathbf{v})](t)}(\mathscr{P}\mathbf{m}_0)\in
C(\bar{I}_n^*,\mathbf{X}_n)$ and
\begin{equation}\label{0523}
\|\mathscr{M}^{-1}_{\mathscr{S}(\mathbf{v})}(\mathscr{P}\mathbf{m}_0)
\|_{C(I_t,\,\mathbf{X}_n)}\leq \underline{\rho}^{-1}e^{Kt}
\|\mathscr{P}\mathbf{m}_0\|_{\mf{X_n}}, \quad t\in I_n^*.  \end{equation}

Finally, we denote
\begin{equation}\label{jiangfei0591}\begin{aligned}
\mathbf{u} &=\mathscr{M}_{\mathscr{S}\mf{(v)}}^{-1} \mathscr{M}_{\partial_t\mathscr{S}\mf{
(v)}}\mathscr{M}_{\mathscr{S}\mf{(v)}}^{-1}\bigg(\mathscr{P}\mathbf{m}_0+\int_0^t
[\mathscr{P}\mathscr{N}(\mathscr{S}(\mathbf{v}),\mathbf{v},\mathscr{D}(\mf{v}))](s)\mathrm{d}s\bigg)\\
&\quad +\mathscr{M}_{\mathscr{S}\mf{ (v)}}^{-1}[
\mathscr{P}\mathscr{N}(\mathscr{S}(\mathbf{v}),\mathbf{v},\mathscr{D}(\mf{v}))].
\end{aligned}\end{equation}
Thus, we can use (\ref{0404}), (\ref{0508}), (\ref{0509}) and
(\ref{0517}) to evaluate (\ref{jiangfei0591}) and obtain that
\begin{equation*}
\|\mathbf{u}\|_{\mathbf{X}_n}\leq \tilde{h}_4(K,\bar{\rho},
\|\mf{d}_0\|_{\mf{H}^2(\Omega)},\underline{\rho},\|\nabla^3
\mf{d}_0\|_{\mf{L}^2(\Omega)},T,n)\big(\|\partial_t\mathscr{S}(\mathbf{v})
\|_{L^1(\Omega)}(1+\|\mathscr{P}\mf{m}_0\|_{\mf{X}_n })+1\big),
\end{equation*}
which, together with (\ref{0407}) and H\"older's inequality, yields
\begin{equation}\label{n0527}
\|\mathbf{u}\|_{L^2(I_t,\mathbf{X}_n)}\leq
h_5(K,\bar{\rho},\|\rho_0\|_{H^1(\Omega)},
\|\mf{d}_0\|_{\mf{H}^2(\Omega)}, \underline{\rho},\|\nabla^3
\mf{d}_0\|_{\mf{L}^2(\Omega)},T,n)\sqrt{t}(1+\|\mathscr{P}\mf{m}_0\|_{\mf{X}_n
}), \end{equation} where $h_5$ is nondecreasing in its first four
variables and nonincreasing in its fifth variable.

 \subsection{Local existence}

 Now we are in a position to prove existence of a local solution by applying a fixed point theorem. To this
 end, we assume that
\begin{equation}\label{n0528}
5\mathrm{max}\bigg\{\frac{\|\mathscr{P}\mathbf{m}_0\|_{\mathbf{X}_n}}{\underline{\rho}},
\|\mathbf{v}(0)\|_{\mathbf{X}_n}\bigg\}<K,\end{equation}
and take
\begin{equation}\label{n0529}
T_0:=T_0(\bar{\rho},\|\rho_0\|_{H^1(\Omega)},\|\mathbf{d}_0\|_{\mf{H}^3
(\Omega)},\underline{\rho},K,T,n)=\mathrm{min}\bigg\{\frac{K}{4h_4},\frac{\mm{ln}(5/4)}{K},
{\frac{K^2}{16h_5^2}},{\frac{25}{16\underline{\rho}^2 h_5^2}},T_d^{\tilde{K}}\bigg\},
\end{equation}
so that $T_0$ is nonincreasing in its first three variables and
nodecreasing in the fourth variable. With this choice we have
$$ h_4T_0\leq \frac{K}{4},\quad \underline{\rho}^{-1}e^{K T_0}\|\mathscr{P}\mathbf{m}_0\|_{\mathbf{X}_n}<\frac{K}{4},
\quad h_5\sqrt{T_0}(1+\|\mathscr{P}\mf{m}_0\|_{\mf{X}_n })< \frac{K}2. $$
Therefore, by virtue of (\ref{0522}), (\ref{0523}), (\ref{n0527}) and (\ref{n0528}), the mapping
\begin{eqnarray}\label{n0530}
&&\displaystyle\mathcal{T}:\ \mathbb{A}\rightarrow \mathbb{A},\nonumber\\
&&\label{jsn0530}\mathcal{T}(\mathbf{w}):=\mathscr{M}^{-1}_{[\mathscr{S}(\mathbf{w})]}
\big\{\mathscr{P}\mathbf{m}_0+\displaystyle\int_0^t[\mathscr{P}\mathscr{N}(\mathscr{S}(\mathbf{w}),
\mathbf{w},\mathscr{D}(\mathbf{w}))](s)\mathrm{d}s\big\},\end{eqnarray}
maps
\begin{equation*}B_{K,\tau_0}=\big\{\mathbf{w}\in \mathbb{A}~|~
 \|\mathbf{w}\|_{C(\bar{I}_{\tau_0},\mathbf{X}_n)}
+\|\partial_t\mathbf{w}\|_{L^2({I}_{\tau_0},\mathbf{X}_n)}\leq
K\big\}\end{equation*} into itself for any $0<\tau_0\leq T_0$, where we can take
%%% the both spaces
$\mf{X}_n=\mf{H}^1_0(\Omega)$.

 In the next step, we prove that $\mathcal{T}$ is contractive. Keeping in mind that
$$\mathscr{M}^{-1}_{\rho_1}(\mathbf{w}_1)-\mathscr{M}^{-1}_{\rho_2}(\mathbf{w}_2)
=(\mathscr{M}^{-1}_{\rho_1}-\mathscr{M}_{\rho_2}^{-1})
(\mathbf{w}_1)+\mathscr{M}^{-1}_{\rho_2}(\mathbf{w}_1-\mathbf{w}_2),$$
we get from (\ref{jsn0530}) that
\begin{equation}\label{n0533}\begin{aligned}
&\mathcal{T}(\mathbf{w}_1)-\mathcal{T}(\mathbf{w}_2)\\
&=(\mathscr{M}^{-1}_{\mathscr{S}(\mathbf{w}_1)}-\mathscr{M}^{-1}_{\mathscr{S}(\mathbf{w}_2)})
\left(\mathscr{P}\mathbf{m}_0+\int_0^t[\mathscr{P}\mathscr{N}(\mathscr{S}(\mathbf{w}_1
),\mathbf{w}_1,\mathscr{D}(\mathbf{w}_1))](s)\mathrm{d}s\right)\\
&\quad +\mathscr{M}^{-1}_{\mathscr{S}(\mathbf{w}_2)}\int_0^t\left(\mathscr{P}\mathscr{N}(\mathscr{S}(\mathbf{w}_1
),\mathbf{w}_1,\mathscr{D}(\mathbf{w}_1)))-\mathscr{P}\mathscr{N}(\mathscr{S}(\mathbf{w}_2
),\mathbf{w}_2,\mathscr{D}(\mathbf{w}_2)))\right)(s)\mathrm{d}s.
\end{aligned}\end{equation}
Recalling that $\underline{\rho}\leq \rho(0)\leq \bar{\rho}$, one has
\begin{equation} \label{wdh0432}
\|\mathscr{P}\mf{m}_0\|_{\mf{X}_n}\leq c(n)\bar{\rho}\|\mf{v}(0)\|_{\mf{X}_n}\leq c(n)\bar{\rho}K .
\end{equation}
We apply (\ref{0517}), (\ref{n0521}) and (\ref{wdh0432}) to bound the first term in
(\ref{n0533}), and use (\ref{0404}), (\ref{0409}), (\ref{n0456}), (\ref{0508}), (\ref{0520})
and H\"older's inequality to majorize the second term in (\ref{n0533}), and obtain
\begin{equation}\label{n0543}\begin{aligned}
&\|\mathcal{T}(\mathbf{w}_1) -\mathcal{T}(\mathbf{w}_2)\|_{\mf{X}_n}(t)   \\
&\leq h_6(\bar{\rho},\|\rho_0\|_{H^1(\Omega)},\|\mathbf{d}_0\|_{\mf{H}^2(\Omega)},
\underline{\rho},\|\nabla^3\mf{d}_0\|_{\mf{L}^2(\Omega)},K,T,n)
\sqrt{t}\|\mathbf{w_1}-\mathbf{w}_2 \|_{C^0(\bar{I}_t,\mf{X}_n)},
\end{aligned} \end{equation}
for $t\in [0,T_0]$, $\mathbf{w}_1,\mathbf{w}_2\in B_{K,T_0}$, where $h_6$ is
nondecreasing in its first three variables and nonincreasing in the fourth variable.

Similar to (\ref{n0533}), we also have the identity
\begin{equation*}\label{n0535}\begin{aligned}
&\partial_t\mathcal{T}(\mathbf{w}_1)-\partial_t\mathcal{T}(\mathbf{w}_2)  \\
&=(\mathscr{M}_{\mathscr{S}(\mf{w}_1)}^{-1}-\mathscr{M}_{\mathscr{S}(\mf{ w}_2)}^{-1})
\mathscr{M}_{\partial_t\mathscr{S}(\mf{w}_1)}\mathscr{M}_{\mathscr{S}(\mf{
w}_1)}^{-1}\bigg\{\mathscr{P}\mathbf{m}_0+\int_0^t
[\mathscr{P}\mathscr{N}(\mathscr{S}(\mathbf{w}_1),\mathbf{w}_1,\mathscr{D}(\mf{w}_1))](s)\mathrm{d}s\bigg\} \\
&\quad +\mathscr{M}_{\mathscr{S}(\mf{ w}_2)}^{-1}(\mathscr{M}_{\partial_t\mathscr{S}
(\mf{w}_1)}-\mathscr{M}_{\partial_t\mathscr{S}(\mf{w}_2)})\mathscr{M}_{\mathscr{S}(\mf{
w}_1)}^{-1}\bigg\{\mathscr{P}\mathbf{m}_0+\int_0^t
[\mathscr{P}\mathscr{N}(\mathscr{S}(\mathbf{w}_1),\mathbf{w}_1,\mathscr{D}(\mf{w}_1))](s)\mathrm{d}s\bigg\}\\
&\quad +\mathscr{M}_{\mathscr{S}(\mf{ w}_2)}^{-1} \mathscr{M}_{\partial_t\mathscr{S}(\mf{
w}_2)}(\mathscr{M}_{\mathscr{S}(\mf{w}_1)}^{-1}-\mathscr{M}_{\mathscr{S}
(\mf{w}_2)}^{-1})\bigg\{\mathscr{P}\mathbf{m}_0 +\int_0^t
[\mathscr{P}\mathscr{N}(\mathscr{S}(\mathbf{w}_1),\mathbf{w}_1,\mathscr{D}(\mf{w}_1))](s)\mathrm{d}s\bigg\}\\
&\quad +\mathscr{M}_{\mathscr{S}(\mf{ w}_2)}^{-1}
\mathscr{M}_{\partial_t\mathscr{S}(\mf{w}_2)}\mathscr{M}_{\mathscr{S}(\mf{ w}_2)}^{-1}\bigg\{\int_0^t
[\mathscr{P}\mathscr{N}(\mathscr{S}(\mathbf{w}_1),\mathbf{w}_1,\mathscr{D}(\mf{w}_1))\\
&\qquad -\mathscr{P}\mathscr{N}(\mathscr{S}(\mathbf{w}_2),\mathbf{w}_2,\mathscr{D}(\mf{w}_2))](s)
\mathrm{d}s\bigg\}+(\mathscr{M}_{\mathscr{S}(\mf{w}_1)}^{-1}-\mathscr{M}_{\mathscr{S}
(\mf{ w}_2)}^{-1})[\mathscr{P}\mathscr{N}(\mathscr{S}(\mathbf{w}_1),\mathbf{w}_1,\mathscr{D}(\mf{w}_1))]\\
&\quad +\mathscr{M}_{\mathscr{S}(\mf{ w}_2)}^{-1}[ \mathscr{P}\mathscr{N}(\mathscr{S}
(\mathbf{w}_1),\mathbf{w}_1,\mathscr{D}(\mf{w}_1))- \mathscr{P}\mathscr{N}(\mathscr{S}
(\mathbf{w}_2),\mathbf{w}_2,\mathscr{D}(\mf{w}_2))].\end{aligned}\end{equation*}
Now, employing (\ref{0404}), (\ref{n0456}), (\ref{0507}), (\ref{0508}),
(\ref{n0510}), (\ref{0517}), (\ref{0520}), (\ref{n0521}) and (\ref{wdh0432})
to control the six terms on the right-hand side of the equality above, we deduce that
\begin{equation}\begin{aligned}\label{n05361}
& \|\partial_t\mathcal{T}(\mathbf{w}_1) -\partial_t\mathcal{T}(\mathbf{w}_2)\|_{\mf{X}_n}(t) \\
& \leq h_6(\bar{\rho},\|\rho_0\|_{H^1(\Omega)},\|\mathbf{d}_0\|_{\mf{H}^2(\Omega)},
\underline{\rho},\|\nabla^3 \mf{d}_0\|_{\mf{L}^2(\Omega)},K,T,n)
\Big(\|\partial_t[\mathscr{S}(\mf{w}_1)-\mathscr{S}(\mf{w}_2)](t)\|_{L^1(\Omega)} \\
&\quad +(\|\partial_t\mathscr{S}(\mf{w}_1)(t)\|_{L^1(\Omega)}+1)
\|\mathbf{w_1}-\mathbf{w}_2 \|_{C^0(\bar{I}_t,\mathbf{X}_n)}
+\|\nabla (\mathbf{d}_1-\mathbf{d}_2)(t)\|_{\mathbf{L}^2(\Omega)}\Big),
\end{aligned} \end{equation}
for a.e. $t\in [0,T_0]$ and for any $\mathbf{w}_1,\mathbf{w}_2\in B_{K,T_0}$,
where $h_6$ is nondecreasing in its three variables and
nonincreasing in the fourth variable again.

%Using the formula
%$\mf{w}_1(t)-\mf{w}_2(t)=\int_0^t\partial_t(\mf{w}_1-\mf{w}_2)(s)\mm{d}s$, we have
%\begin{equation}\label{n0537}\|(\mf{w}_1-\mf{w}_2)(t)\|_{\mf{X}_n}\leq
%\sqrt{t}\|\partial_t(\mf{w}_1-\mf{w}_2)\|_{L^2(I_t,\mf{X}_n)}.\end{equation}
Consequently, substituting (\ref{0407}), (\ref{n04091}),
and (\ref{n0456}) into (\ref{n05361}), we find that
\begin{equation}\begin{aligned}\label{n0536}&\|\partial_t\mathcal{T}(\mathbf{w}_1)
-\partial_t\mathcal{T}(\mathbf{w}_2)\|_{L^2(I_t,\mf{X}_n)}  \\
&\leq h_6(\bar{\rho},\|\rho_0\|_{H^1(\Omega)},\|\mathbf{d}_0\|_{\mf{H}^2(\Omega)},
\underline{\rho},\|\nabla^3 \mf{d}_0\|_{\mf{L}^2(\Omega)},K,T,n)
\sqrt{t} \|\mathbf{w_1}-\mathbf{w}_2 \|_{C^0(\bar{I}_t,\mf{X}_n)}.
\end{aligned} \end{equation}

Adding (\ref{n0536}) to (\ref{n0543}), we finally get
\begin{equation*}\begin{aligned}&\|\mathcal{T}(\mathbf{w}_1)-\mathcal{T}(\mathbf{w}_2)\|_{\mf{X}_n}(t)
+\|\partial_t\mathcal{T}(\mathbf{w}_1)-\partial_t\mathcal{T}
(\mathbf{w}_2)\|_{L^2(I_{t},\mf{X}_n)}  \\
& \leq \tilde{h}_6(\bar{\rho},\|\rho_0\|_{H^1(\Omega)},\|\mathbf{d}_0\|_{\mf{H}^2(\Omega)},
\underline{\rho},\|\nabla^3 \mf{d}_0\|_{\mf{L}^2(\Omega)},K,T,n){t}\\
&\quad \left(\|\mathbf{w_1}-\mathbf{w}_2
\|_{C^0(\bar{I}_t,\mathbf{X}_n)}+\|\partial_t(\mathbf{w_1}-\mathbf{w}_2)
\|_{L^2(\bar{I}_t,\mathbf{X}_n)}\right)  \end{aligned} \end{equation*}
for any $t\in [0,T_0]$.

 If we take
$${T}_n^*=\mathrm{\min}\left\{T_0,\frac{1}{2\tilde{h}_6}\right\},$$
then $\mathcal{T}$ maps $B_{K,{T}_n^*}\subset C(\bar{I}_{{T}_n^*},\mathbf{X}_n)$
into itself and is contractive. Therefore, it possesses in $B_{K,{T}_n^*}$ a unique
fixed point $\mathbf{v}$ which satisfies (\ref{0513}). Thus, we have a solution
$(\rho=\mathscr{S}(\mathbf{v}),\mathbf{v},\mathscr{D}(\mathbf{v}))$
which is defined in $ Q_{{T}_n^*}$ and satisfies the
initial-boundary value problem (\ref{n301})--(\ref{n0306}) for each given
$n$. Moreover, we see that ${T}_n^*$ has the form
\begin{equation}\label{n0540}
0<{T}_n^*={h}({\bar{\rho},\|\rho_0\|_{H^1(\Omega)},\|\mathbf{d}_0\|_{\mf{H}^2(\Omega)},
\underline{\rho},\|\nabla^3 \mf{d}_0\|_{\mf{L}^2(\Omega)},K,T,n}) \leq T_0,
\end{equation}
where $h$ is nonincreasing in its first three variables and nondecreasing
 in the fourth variable. This means that we can find a unique maximal solution
$(\rho_n,\mf{v}_n,\mf{d}_n)$ defined in $[0,T_n)\times \Omega$ for
each given $n$, where $T_n\leq T$.

\subsection{Global existence}\label{sec:0404}

In order to show the maximal time $T_n=T$ for any $n$, it
suffices to derive uniform bounds for $\rho_n$, $\mathbf{v}_n$,
$\mathbf{d}_n$ and $\mathscr{P}_n(\rho_n\mathbf{v}_n)$. However, we
need to impose an additional (smallness) condition on the initial
approximate energy to get the uniform bound of
$\|\mf{d}_n\|_{L^\infty(\bar{I}_n,\mathbf{H}^2(\Omega))}$. We remark
that $T_h^*$ in (\ref{n0540}) depends on $\|\nabla^3
\mf{d}_0\|_{\mf{L}^2(\Omega)}$, since we have used the fact
$\mf{d}^{k-1}=\mf{d}_0$ on $\partial\Omega$ and the elliptic
estimate to obtain (\ref{nestf0429}), and thus we need an auxiliary
term $\|\nabla^3 \mf{d}_0\|_{\mf{L}^2(\Omega)}$. Such auxiliary term will
not changes, if only the initial data of $\mf{d}(\mf{x},t)$
changes and the boundary value of $\mf{d}(\mf{x},t)$ does not change.
This is why we do not estimate the uniform bound of
$\|\nabla ^3\mf{d}_n\|_{L^\infty(\bar{I}_n,\mathbf{L}^2(\Omega))}$.

For simplicity of notations, we denote
$$(\rho,\mathbf{v},\mathbf{d},\mathscr{P}\mathbf{m}):=(\rho_n,
\mathbf{v}_n, \mathbf{d}_n,\mathscr{P}_n(\rho_n\mathbf{v}_n)).$$
We mention that in the following estimates the letter $G(\ldots)$ will
denote various positive constants depending on its variables.
%% (may depends on other given physical parameters, but we omit them).

First, we derive energy estimates similar to Proposition
\ref{lem:0201}. Differentiating (\ref{0513}) with respect to $t$,
integrating by parts, employing (\ref{n301}) and the
regularity of $(\rho,\mf{v})$, we obtain
\begin{equation*}\begin{aligned}
& \frac{d}{dt}{{E}_\delta}(t) +\int_\Omega \left(\mu |\nabla
\mf{v}|^2+(\lambda+\mu)|\mm{div}\mf{v}|^2
+\varepsilon\delta\beta\rho^{\beta-2}|\nabla \rho|^2
\right)\mathrm{d}\mathbf{x}\leq \nu\int_\Omega (\nabla
\mf{d})^T\Delta \mf{d}\cdot\mf{v}\mathrm{d}\mathbf{x},
\end{aligned}\end{equation*}
where
$$\displaystyle{{E}}_\delta(t)
=\int_{\Omega}\left(\frac{1}{2}\frac{|\mathbf{m}|^2}{\rho}1_{\{\rho>0\}}+Q(\rho)
+\frac{\delta}{\beta-1}\rho^\beta\right)\mm{d}\mf{x}.$$
Noting that
$\mf{d}$ belongs to the function class (\ref{0465}) with $T_n$ in
place of $T_d^K$, $\mf{d}$ satisfies (\ref{0213}), and one further has
\begin{equation}\label{0542}\begin{aligned}&\frac{d}{dt}{\mathcal{E}}_\delta(t)
+\int_\Omega \left(\mu |\nabla
\mf{v}|^2+(\lambda+\mu)|\mm{div}\mf{v}|^2  +\nu\theta(|\Delta \mf{d}+|\nabla
\mf{d}|^2\mf{d}|^2)+\varepsilon\delta\beta\rho^{\beta-2}|\nabla
\rho|^2\right)\mathrm{d}\mathbf{x}\leq 0,
\end{aligned}\end{equation}
where
$$\displaystyle {\mathcal{E}}_\delta(t):={\mathcal{E}}_\delta(\rho,\mf{m},\mf{d}):
=\int_{\Omega}\left(\frac{1}{2}\frac{|\mathbf{m}|^2}{\rho}1_{\{\rho>0\}}+\frac{A}{\gamma-1}\rho^\gamma
+\frac{\delta}{\beta-1}\rho^\beta+\frac{\nu|\nabla
\mf{d}|^2}{2}\right)\mm{d}\mf{x}$$
with $\mf{m}=\rho\mf{v}$. Integrating (\ref{0542}) over $(0,t)$, we get
\begin{equation*}\begin{aligned}
& {\mathcal{E}}(t)+\int_0^t\int_\Omega \left(\mu |\nabla
\mf{v}|^2+(\lambda+\mu)|\mm{div}\mf{v}|^2  +\nu\theta(|\Delta
\mf{d}+|\nabla \mf{d}|^2\mf{d}|^2)\right)\mathrm{d}\mathbf{x}\mathrm{d}s\\
&\leq {\mathcal{E}}_\delta(0):={\mathcal{E}}_\delta(\rho_0,\mf{m}_0,\mf{d}_0)<
\bar{\mathcal{E}}_{\delta,0},
\end{aligned}\end{equation*}
where $\bar{\mathcal{E}}_{\delta,0}$ is a given positive constant. In particular,
\begin{eqnarray}
\label{n0544}&&\|\sqrt{\rho}\mathbf{v}\|_{L^\infty(I_n,\mathbf{L}^{2}(\Omega))}\leq
2\bar{\mathcal{E}}_{\delta,0},  \\
&&\label{n054412}    \|\nabla\mathbf{d}\|_{L^\infty(I_n,\mathbf{L}^{2}(\Omega))}\leq
2\bar{\mathcal{E}}_{\delta,0}/\nu,  \\
&& \label{n0544n}   \|\mathbf{v}\|_{L^2(I_n,\mathbf{H}^{1}(\Omega))}\leq
c(\Omega)\bar{\mathcal{E}}_{\delta,0}/\mu.
\end{eqnarray}

With the help of (\ref{n0544})--(\ref{n0544n}), we can deduce more
uniform bounds on $(\rho,\mathbf{v})$.  Using (\ref{0404}) and
(\ref{n0544n}), thanks to the norm of equivalence on
 $\mathbf{X}_n$ (see (\ref{0506}))), we find that
\begin{equation}\label{n0547}\begin{aligned}
G_1(\underline{\rho},\bar{\mathcal{E}}_{\delta,0},T,n)\leq
\rho\leq G_2(\bar{\rho},\bar{\mathcal{E}}_{\delta,0},T,n),
\end{aligned}\end{equation}
from which, (\ref{n0544}) and (\ref{0506}), it follows that
\begin{equation}\label{n0548}\begin{aligned}\|\mathbf{v}\|_{C^0(\bar{I}_{n},\mathbf{X}_n)}
\leq G(\underline{\rho},\bar{\mathcal{E}}_{\delta,0},T,n)
\end{aligned}\end{equation}
and
\begin{equation}\label{n0549}\begin{aligned}
\|\mathscr{P}(\rho\mathbf{v})\|_{C^0(\bar{I}_{n}, \mathbf{X}_n)}\leq
c(n)\|\rho\mathbf{v}\|_{C^0(\bar{I}_{n},\mathbf{L}^2(\Omega))}\leq
G(\bar{\rho},\bar{\mathcal{E}}_{\delta,0},T,n).\end{aligned}\end{equation}

Applying  (\ref{n0548}) to (\ref{0405}), one gets
\begin{equation}\label{nn0548}\begin{aligned}
\| \rho\|_{C^0(\bar{I}_n,H^1(\Omega))}\leq
G(\underline{\rho},\|\rho_0\|_{H^1(\Omega)},\bar{\mathcal{E}}_{\delta,0}, T,n).
\end{aligned}\end{equation}
%Moreover, in view of (\ref{n0409}) with $r=4$, we find that
%\begin{equation}\label{n0552}\|\partial_t\rho\|_{L^4(I_n,{L}^2\Omega))}
%\leq G(\underline{\rho},\|\rho_0\|_{H^1(\Omega)},
%\|\rho_0\|_{{H}^{\frac{3}{2}}(\Omega)},\bar{\mathcal{E}}_{\delta,0},T,n).\end{equation}
Utilizing (\ref{n054412}), (\ref{n0547})--(\ref{n0549}),
arguing similarly to that for (\ref{n0527}), we obtain from (\ref{n0413nn}) that
%\begin{equation}\label{n05551}\|\partial_t\mathbf{v}\|_{L^4(I_n,\mathbf{X}_n)}\leq
%G(\underline{\rho},\bar{\rho},\|\rho_0\|_{H^1(\Omega)},\|\rho_0\|_{{H}^{\frac{3}{2}}(\Omega)},
%\|\mf{d}_0\|_{\mf{H}^3(\Omega)},\bar{\mathcal{E}}_{\delta,0},T,n);\end{equation}
%and in particular, by H\"older's inequality,
\begin{equation}\label{n0555}
\|\partial_t\mathbf{v}\|_{L^2(I_n,\mathbf{X}_n)}\leq
G(\underline{\rho},\bar{\rho},\|\rho_0\|_{H^1(\Omega)},%\|\rho_0\|_{{H}^{\frac{3}{2}}(\Omega)},
\|\mf{d}_0\|_{\mf{H}^2(\Omega)},\bar{\mathcal{E}}_{\delta,0},T,n),
\end{equation}
%where any $r>2$ is sufficient, and we have taken $r=4$ for simplicity.

Hence, we have shown the uniform boundedness of $\underline{\rho}$,
$\bar{\rho}$, $\|\rho\|_{C^0(\bar{I}_n,\mathbf{X}_n)}$,
$\|\mathscr{P}(\rho\mathbf{v})\|_{C^0(\bar{I}_{n}, \mathbf{X}_n)}$,
$\|\mathbf{v}\|_{C^0(\bar{I}_n,\mathbf{X}_n)}$ and
$\|\partial_t\mathbf{v}\|_{L^2(I_n,\mathbf{X}_n)}$. It remains to
show the uniform boundedness of
$\|\mf{d}\|_{L^\infty(\bar{I}_n,\mf{H}^2(\Omega))}$, which can be
obtained by following the spirit of proof in Section \ref{030202}.
For the reader's convenience, we give the proof in the following.

First, we rewrite the energy inequality (\ref{0542}) in the following form as in (\ref{0219}):
$$ {\mathcal{E}}(t) +\int_0^t\int_\Omega \left[\mu |\nabla \mf{v}|^2
+(\lambda+\mu)|\mm{div}\mf{v}|^2 +\nu\theta|\Delta \mf{d}|^2
+\varepsilon\delta\beta\rho^{\beta-2}|\nabla
\rho|^2\right]\mathrm{d}\mathbf{x} \leq \bar{\mathcal{E}}_{\delta,0}+\nu\theta\||\nabla
\mf{d}|\|_{L^4(\Omega)}^4. $$
Let $\bar{\mathcal{E}}_{\delta,0}=\nu/4094$. Arguing in the same manner
as in the derivation of (\ref{0227}), we deduce that
\begin{equation}\label{n0452}
\begin{aligned}
&\|\nabla^2 \mf{d}\|_{\mf{L}^2(Q_T)}+\|\nabla
\mf{d}\|_{\mf{L}^4(Q_T)}+\|\partial_t\mathbf{d}\|_{L^{4/3}(I,\mathbf{L}^2(\Omega))}\leq
G(\|\mf{d}_0\|_{\mathbf{H}^{2}(\Omega)},\bar{\mathcal{E}}_{\delta,0}).
\end{aligned}\end{equation}

Now we proceed to derive uniform bounds of higher derivatives on
$\mathbf{d}$. Differentiating (\ref{n302}) with respect to $t$,
multiplying the resulting equations by $\partial_t\mf{d}$ in $L^2(\Omega )$,
recalling $|\mf{d}|=1$, we integrate by parts to infer that
\begin{equation}\label{n228}\begin{aligned}
&\frac{d}{dt}\int_\Omega|\partial_t\mf{d}|^2\mm{d}\mf{x}\\
&=2\int_\Omega\partial_t\mf{d}\cdot\left(\theta\Delta\partial_t\mf{d}+\theta|\nabla
\mf{d}|^2\partial_t\mf{d}+2\theta(\nabla \mf{d}:\nabla
\partial_t\mf{d})\mf{d}-\partial_t\mf{v}\cdot\nabla \mf{d}-\mf{v}\cdot\nabla
\partial_t\mf{d}\right)\mm{d}\mf{x}\\
& \leq -\theta\|\nabla
\partial_t\mf{d}\|_{\mathbf{L}^2(\Omega)}^2+(10\theta+1)\||\partial_t\mf{d}||\nabla
\mf{d}\|^2_{{L}^2(\Omega)}
+2\|\partial_t\mf{v}\|_{\mathbf{L}^2(\Omega)}^2
+\frac{1}{\theta}\|\mf{v}\|_{\mf{L}^\infty(Q_T)}^2\|\partial_t\mf{d}\|_{\mathbf{L}^2(\Omega)}^2,
\end{aligned}\end{equation}
where the second term on the right-hand side of (\ref{n228}) can be
bounded as follows, using (\ref{jjjw0432ww}), H\"older's and the triangle inequalities.
\begin{equation}\label{n2229}\begin{aligned}
(10\theta+1)\||\partial_t\mf{d}||\nabla \mf{d}\|^2_{{L}^2(\Omega)}
\leq & (10\theta+1)\||\partial_t\mf{d}|\|^2_{{L}^4(\Omega)}
\||\nabla
\mf{d}|\|^2_{{L}^4(\Omega)}\\
\leq &(10\theta+1)\sum_{i=1}^2\|\partial_t{d}_i
\|_{{L}^4(\Omega)}^2 \||\nabla \mf{d}|\|^2_{{L}^4(\Omega)}\\
\leq&(10\theta+1)\sum_{i=1}^2\sqrt{2}\|\partial_t{d}_i\|_{{L}^2(\Omega)}\|\nabla
\partial_t{d}_i\|_{{L}^2(\Omega)}\||\nabla
\mf{d}|\|^2_{{L}^4(\Omega)}\\
\leq &\frac{\theta}{2}\|\nabla
\partial_t\mf{d}\|_{\mathbf{L}^2(\Omega)}^2
+\frac{{}(10\theta+1)^2}{\theta}\|\partial_t\mf{d}\|_{\mathbf{L}^2(\Omega)}^2\||\nabla
\mf{d}|\|^4_{{L}^4(\Omega)}.
\end{aligned}\end{equation}
Inserting (\ref{n2229}) into (\ref{n228}), we conclude that
\begin{equation*}  \begin{aligned}
&\frac{d}{dt}\|\partial_t\mf{d}\|^2_{\mathbf{L}^2(\Omega)}+\frac{\theta}{2}\|\nabla
\partial_t\mf{d}\|_{\mathbf{L}^2(\Omega)}^2   \\
 &\leq \left(\frac{{}(10\theta+1)^2}{\theta}\||\nabla
\mf{d}|\|^4_{{L}^4(\Omega)}+\frac{1}{\theta}\|\mf{v}\|_{\mf{L}^\infty(Q_T)}^2\right)
\|\partial_t\mf{d}\|_{\mathbf{L}^2(\Omega)}^2+2\|\partial_t\mf{v}\|_{\mathbf{L}^2(\Omega)}^2.
\end{aligned}\end{equation*}
Thus, applying Gronwall's inequality, we have
\begin{equation*}  \begin{aligned}
\|\partial_t\mf{d}\|^2_{\mathbf{L}^2(\Omega)}\leq&\left(
\|\partial_t\mf{d}(0)\|_{\mathbf{L}^2(\Omega)}^2+2\|\partial_t\mf{v}\|_{\mathbf{L}^2(Q_T)}^2\right)
e^{\frac{(10\theta+1)}{2\theta}\||\nabla
\mf{d}|\|^4_{{L}^4(Q_T)}+\frac{T}{\theta}\|\mf{v}\|_{\mf{L}^\infty(Q_T)}^2}.
\end{aligned}\end{equation*}

Noting that
$$\|\partial_t\mf{d}(0)\|_{\mathbf{L}^2(\Omega)}=\|\theta(\Delta
\mathbf{d}_0+|\nabla \mathbf{d}_0|^2\mf{d}_0)-\mathbf{v}_0\cdot
\nabla \mathbf{d}_0\|_{\mf{L}^2(\Omega)},$$
we use (\ref{n0548}), (\ref{n0555}) and (\ref{n0452}) to arrive at
\begin{equation}\label{nnn0457}\begin{aligned}
\|\partial_t\mf{d}\|^2_{L^\infty(I_n,\mathbf{L}^2(\Omega))}+\|\nabla
\partial_t\mf{d}\|_{L^2(I_n,\mathbf{L}^2(\Omega))}^2 \leq &
G(\underline{\rho},\bar{\rho},\|\rho_0\|_{H^1(\Omega)},
\|\mf{d}_0\|_{\mf{H}^2(\Omega)},\bar{\mathcal{E}}_{\delta,0},T,n)
\end{aligned}\end{equation}

Recalling that $|\mf{d}|\equiv 1$, it follows from (\ref{n302}) that
\begin{equation}\label{n0458}\begin{aligned}
 \theta^2\int_\Omega|\Delta\mf{d}|^2\mm{d}\mf{x}\leq &
 3\theta^2\int_\Omega|\nabla\mf{d}|^4\mm{d}\mf{x}
+3\int_\Omega|\partial_t\mf{d}|^2\mm{d}\mf{x}+3\int_\Omega|\mf{v}|^2
|\nabla\mf{d}|^2\mm{d}\mf{x} .
\end{aligned}\end{equation}
Similarly to the derivation of (\ref{n0228j}), the first term on the right-hand side of (\ref{n0458})
can be estimated as follows.
\begin{equation}\label{n0459}\begin{aligned}
 & 3\theta^2\int_\Omega|\nabla  \mf{d}|^4\mm{d}\mf{x} \\
 &\quad \leq {384\theta^2} (c(\Omega)+12)^3\|\mf{d}_0\|_{\mathbf{H}^2(\Omega)}^2
\left(\|\mf{d}_0\|_{\mathbf{H}^2(\Omega)}^2+\frac{2\mathcal{E}_0}{\nu}
\right)+\frac{3072\theta^2\mathcal{E}_0}{\nu}\|\Delta\mf{d}\|_{L^2(\Omega)}^2.
\end{aligned}\end{equation}
Putting (\ref{n0458}) and (\ref{n0459}) together, using
(\ref{n0548}), (\ref{nnn0457}) and (\ref{n054412}), we
get\begin{equation}\label{nn0460}\begin{aligned}
 \frac{\theta^2}{4}\int_\Omega|\Delta\mf{d}|^2\mm{d}\mf{x}
 \leq& (c(\Omega)+12)^3\|\mf{d}_0\|_{\mathbf{H}^{2}(\Omega)}^2
\left(\|\mf{d}_0\|_{\mathbf{H}^{2}(\Omega)}^2+\frac{2\mathcal{E}_0}{\nu} \right) \\
& + 3\|\partial_t\mf{d}\|_{\mathbf{L}^2(\Omega)}^2+6\|\mf{v}\|_{\mf{L}^\infty(\Omega)}^2\|\nabla
\mf{d}\|_{\mathbf{L}^2(\Omega)}^2 \\
\leq & G(\underline{\rho},\bar{\rho},\|\rho_0\|_{H^1(\Omega)},
\|\mf{d}_0\|_{\mf{H}^2(\Omega)},\bar{\mathcal{E}}_{\delta,0},T,n).
\end{aligned}\end{equation}
Hence, by  the elliptic estimate (\ref{0221}), (\ref{n054412}),
(\ref{nn0460}) and the fact $|\mf{d}|=1$,
\begin{equation}\label{fzu0447}\begin{aligned}
 \|\mf{d}\|_{L^\infty(I_n,\mf{H}^2(\Omega))}\leq
& G(\underline{\rho},\bar{\rho},\|\rho_0\|_{H^1(\Omega)},
\|\mf{d}_0\|_{\mf{H}^2(\Omega)},\bar{\mathcal{E}}_{\delta,0},T,n).
\end{aligned}\end{equation}

The inequalities (\ref{n0547}),  (\ref{n0548}), (\ref{nn0548}),
  (\ref{n0555}) and (\ref{fzu0447}) furnish the desired estimates which, in combination with
(\ref{n0528}), (\ref{n0529}) and  (\ref{n0540}), give a
possibility to repeat the above fixed point argument to conclude
that $T_n=T$, and the global solution $(\rho_n,\mf{v}_n,\mf{d}_n)$
is unique. To end this section, we summarize our previous results on
the global existence of a unique solution
$(\rho_n,\mf{v}_n,\mathbf{d}_n)$ to the third approximate problem
(\ref{n301})--(\ref{n0306}) as follows.
%%%%%%%%%%%%%%%%%%%%%%%%%%%%%%%%%%%%%%%%%%%%%%%%%%%%%%%%%%%%%%%%%%%%%
\begin{pro}\label{thm:0301}
Let
\begin{equation}\label{condition0452}
\delta>0,\ \beta>0, \ \varepsilon>0,\mbox{ and }0<\underline{\rho}\leq \bar{\rho}<\infty. \end{equation}
Assume that $\Omega$ is a bounded $C^{2,\alpha}$-domain ($\alpha\in(0,1)$),
and the initial data $(\rho_0,\mf{m}_0,\mf{d}_0)$ satisfies
\begin{eqnarray}
&&\label{n0467}{\mathcal{E}_{\delta}}(\rho_0,\mf{m}_0,\mf{d}_0)<\bar{\mathcal{E}}_{\delta,0}:
=\nu /4094,  \\
&&\label{0920n} 0<\underline{\rho}\leq \rho_0\leq
\bar{\rho},\quad \rho_0\in W^{1,\infty}(\Omega),\quad  |\mf{d}_0|=1,\\
 &&\label{n0468}
 %\cap \widetilde{H}^{3/2}(\Omega),
 \mf{v}_0\in \mf{X}_n,\ \mathbf{d}_0\in
\mf{H}^{3}({\Omega}).
\end{eqnarray}
 Then there exists a unique triple $(\rho_n,\mf{v}_n,\mf{d}_n)$ with the following properties:
\begin{enumerate}
  \item[\quad \quad(1)] Regularity.
\begin{equation}\label{ren0454}
\left\{\begin{aligned}
&\rho_n\mbox{ satisfies the regularity as in
Proposition \ref{pro:0401} with }T\mbox{ in place of }T_d^K, \\
&\mf{v}_n\in C^0(\bar{I},\mf{X}_n),\ \partial_t\mf{v}_n\in
L^2(I,\mathbf{X}_n),\  \nabla \rho_n\in L^2(I,{E}_{0}^p(\Omega)),\
\rho_n\mf{v}_n\in C^0({I},{E}_0^p(\Omega)),
\end{aligned}\right.\end{equation}
\begin{equation*}\label{ren0455}
\mf{d}_n\mbox{ satisfies the regularity as in Proposition
\ref{pro:0402} with }T\mbox{ in place of }T_d^K.\qquad\qquad\qquad
\end{equation*}
%%%%%%%%%%%%%%%%%%%%%
  \item[\quad \quad(2)] $(\rho_n,\mf{v}_n,\mf{d}_n)
  $ solves (\ref{n301}) and (\ref{n302}) a.e. in $Q_T$, and
  satisfies (\ref{n303}) and $(\rho_n,\mf{v}_n,\mf{d}_n)|_{t=0}=(\rho_0,\mf{v}_0,\mf{d}_0)$.
%%\begin{eqnarray*}
%&& \int_\Omega\partial_t(\rho_n\mf{v}_n)\cdot
%\Phi\mm{d}\mf{x}+\int_\Omega \bigg[\partial_j(\rho_n\mf{v}_nv_n^j)-\mu
%\Delta \mf{v}_n-(\mu+\lambda)\nabla \mm{div}\mf{v}_n +\nabla A\rho^\gamma_n \\
%&&\; +\nu\mathrm{div}\left(\nabla
%\mathbf{d}_n\otimes\nabla\mathbf{d}_n -\frac{|\nabla
%\mf{d}_n|^2\mathbb{I}}{2}\right)+\delta\nabla \rho_n^\beta
 % +\varepsilon\nabla \rho_n\cdot \nabla \mf{v}_n\bigg]\cdot \mathbf{\Phi}\mm{d}\mf{x}=0,
  %\; t\in I,\Phi\in\mf{X}_n,  \\
%&& \partial_t\rho_n+\mm{div}(\rho_n\mf{v}_n)-\varepsilon\Delta\rho_n=0,\; \mbox{ a.e. in }Q_T,\\
%&& \partial_t\mathbf{d}_n+\mathbf{v}_n\cdot \nabla \mathbf{d}_n=
%\theta(\Delta \mathbf{d}_n+|\nabla \mathbf{d}_n|^2\mf{d}_n),\mbox{
%a.e. in }Q_T.
%\end{eqnarray*}
  \item[\quad \quad(3)] Finite and bounded energy inequalities:
\begin{equation}\begin{aligned}\label{0457enegy}
\frac{d}{dt}{\mathcal{E}}^n_\delta(t)+\mathcal{F}^n(t)+\int_\Omega
\varepsilon\delta\beta \rho_n^{\beta-2}|\nabla
\rho_n|^2(t)\mm{d}\mf{x}\leq 0\;\;\mbox{ in }\mathcal{D}'(I),
\end{aligned}\end{equation}
and
\begin{equation}\begin{aligned}\label{0458enegy}
{\mathcal{E}}_\delta^n(t)+\int_0^t\left(\mathcal{F}^n(s)+\int_\Omega
\varepsilon\delta\beta \rho_n^{\beta-2}|\nabla
\rho_n|^2(s)\mm{d}\mf{x}\right)\mm{d}s\leq
{\mathcal{E}_{\delta}}(\rho_0,\mf{m}_0,\mf{d}_0)\quad\mbox{a.e. in }I,
\end{aligned}\end{equation}
 where ${\mathcal{F}}^n(t):={\mathcal{F}}(\rho_n,\mf{v}_n,\mf{d}_n)$ and
${\mathcal{E}}_\delta^n(t):={\mathcal{E}}_\delta(\rho_n,\mf{m}_n,\mf{d}_n)$
with $\mf{m}_n=\rho_n\mf{v}_n$.
  \item[\quad \quad(4)] Additional uniform estimates.
\begin{eqnarray} && \label{n047712}
|\mf{d}_n|\equiv 1 \mbox{ in }\bar{Q}_T,
%&&\label{n0477123} \begin{aligned}
%& \|\mf{v}_n\|_{L^2(I,\mf{H}^{1}(\Omega))}+\|\rho_n\|_{L^\infty(I,L^{\gamma}(\Omega))}+\|\nabla
%\mathbf{d}_n\|_{L^\infty(I,\mathbf{L}^{2}(\Omega))} \\
%&+
%\|\rho_n|\mf{v}_n|^2\|_{L^\infty(I,L^{1}(\Omega))}
%\delta^{\frac{1}{\beta}}\|\rho_n|\mf{v}_n|^2\|_{L^\infty(I,L^{1}(\Omega))}%+\sqrt
%{\varepsilon\delta}\|\nabla
%\rho_n^{\frac{\beta}{2}}\|_{\mathbf{L}^2(Q_T)}
%\leq G,
%\end{aligned}
\\[1mm]
&&\label{n04771231}\|\nabla^2 \mf{d}_n\|_{\mf{L}^2(Q_T)}+\|\nabla
\mf{d}_n\|_{{L}^4(Q_T)}+\|\partial_t\mathbf{d}_n\|_{L^{4/3}(I,\mathbf{L}^2(\Omega))}\leq
G(\|\mf{d}_0\|_{\mathbf{H}^{{2}}(\Omega)}),\\
&&\label{n0477} \sqrt{\varepsilon}\|\nabla \rho_n\|_{\mf{L}^2(Q_T)}\leq G(\delta), \\
&& \label{n0478}
\|\rho_n\|_{L{^\frac{4\beta}{3}}(Q_T)}  \leq G(\varepsilon,\delta)
\end{eqnarray}
\end{enumerate}
(see \cite[ Section 7.7.5.2]{NASII04} for the proof of (\ref{n0477}) and (\ref{n0478})),
where $G$ is a positive constant which is independent of $n$ and nondecreasing in
its arguments. Moreover, if  $\varepsilon$ is not explicitly written
in the argument of $G$, then $G$ is independent of $\varepsilon$ as well.
\end{pro}
\begin{rem}Here we have used the notation
 ${E}_0^p(\Omega)$
 concerning the spaces of vector fields with
summable divergence introduced in \cite[Section 3.2]{NASII04}. For the
reader's convenience, we give the detailed definition. We set
$$E^{q,p}(\Omega)=\{\mf{g}\in (L^q(\Omega))^2~|~\mm{div}\mf{g}\in L^p(\Omega)\}
\hookrightarrow L^q(\Omega)$$ with the norm
$$\|\mf{g}\|_{E^{q,p}}=\|\mf{g}\|_{\mf{L}^q}+\|\mm{div}\mf{g}\|_{\mf{L}^p(\Omega)}.$$
For the sake of simplicity, ${E}^{p,p}(\Omega)$ is denoted by
${E}^{p}(\Omega)$. Then we define
$${E}^{q,p}_0(\Omega)=\overline{\mathcal{D}(\Omega)}^{E^{q,p}}
\mbox{ and }E_0^{p}=\overline{\mathcal{D}(\Omega)}^{E^p}.$$ The
notation ${E}^{q,p}_0(\Omega)$ will be used in Proposition
\ref{thm:0302}.
\end{rem}

\section{Proof of Theorem \ref{thm:0101}}\label{sec:03}

Once we have established Proposition \ref{thm:0301}, we can obtain
Theorem \ref{thm:0101} by the standard three-level approximation
scheme and the method of weak convergence as in \cite{LPLMTFM98,FENAPHOFJ35801}
for the compressible Naiver-Stokes equations. These arguments have also been
successfully used to establish the existence of weak solutions to
other models from fluid dynamics. Here we briefly describe how to
prove Theorem \ref{thm:0101} and state the existence of solutions to
the first and second approximate problems for the reader's convenience,
we refer to \cite{WDHYCGA} or \cite{FENAPHOFJ35801,NASII04} for the
details of the limit process.

In view of the proof in \cite{WDHYCGA}, it suffices to analyze the
convergence of the supercritical nonlinearity $|\nabla\mf{d}|^2\mf{d}$. In fact,
using the uniform bounds in (\ref{n047712}) and (\ref{n04771231}) on $\mf{d}_n$,
applying the Arzel\`a-Ascoli theorem and Aubin-Lions lemma, and taking
subsequences if necessary, we deduce that
\begin{eqnarray}\label{strong1}
\mathbf{d}_n\rightarrow \mathbf{d}\mbox{ strongly in }
C^0(I,\mathbf{L}^2(\Omega))\cap L^p(I, \mf{H}^{1}(\Omega))\cap
L^r(I, \mf{W}^{1,p}(\Omega))
\end{eqnarray}
for any $p\geq 1$ and $r\in [1,2)$, which implies that
\begin{eqnarray}\label{302}
|\nabla \mathbf{d}_n|^2\mf{d}_n\rightharpoonup |\nabla\mathbf{d}|^2\mf{d}
\;\;\mbox{ weakly in } \mf{L}^2(Q_T).
\end{eqnarray}
Moreover, using Vitali's convergence theorem, and
taking subsequences if necessary, we obtain
\begin{eqnarray}\label{strong2}
|\nabla \mathbf{d}_n|^2\mf{d}_n\rightarrow |\nabla
\mathbf{d}|^2\mf{d}\mbox{ strongly in } \mf{L}^r(Q_T)\;\;\mbox{ for any }r\in [1,2).
\end{eqnarray}
In addition, we also have the regularity $\mf{d}\in L^2(I,{\mf{H}^2(\Omega)})$
and $\partial_t\mathbf{d}\in {L^{2}(I,(\mf{H}^{1}(\Omega))^*)}$.
In view of \cite[Proposition 7.31]{NASII04}, we get consequently
\begin{eqnarray}\label{n0504regularity}
\mf{d}\in C^0(I,{\mf{H}^1(\Omega)}).\end{eqnarray}
 Putting (\ref{302}) and (\ref{n0504regularity}) into the proof of
\cite[Proposition 4.1]{WDHYCGA}, we immediately obtain a weak solution of the second
approximate problem which is the weak limit
$(\rho_\varepsilon,\mf{v}_\varepsilon,\mf{d}_\varepsilon)$ of the solution sequence
$(\rho_n,\mf{v}_n,\mf{d}_n)$ as $n\rightarrow \infty$
constructed in Proposition \ref{thm:0301}. Moreover, refereing to
the conclusions in \cite[Proposition 7.31]{NASII04}, the existence
of solutions to the second approximate problem reads as follows.
%%%%%%%%%%%%%%%%%%%%%%%%%%%%%%%%%%%%%%%%%%%%%%%%%%%%%%%%%%%%%%%%%%%%%%%%%
\begin{pro}\label{thm:0302}
Let $\delta>0$, $\beta\geq\max\{\gamma,8\}$, $\varepsilon>0$, and
$0<\underline{\rho}\leq \bar{\rho}<\infty$. Assume that $\Omega$ is
a bounded $C^{2,\alpha}$-domain ($\alpha\in (0,1)$), and the initial
data $(\rho_0,\mf{m}_0,\mf{d}_0)$ satisfies (\ref{n0467}),
(\ref{0920n}) and
\begin{eqnarray}
%&&\label{n0505jfei} 0<\underline{\rho}\leq \rho_0\leq \bar{\rho},\
%\rho\in W^{1,\infty}(\Omega),\ |\mf{d}_0|=1,\\
 \label{n0505}\mf{m}_0\in \mf{L}^2(\Omega),\ \mathbf{d}_0\in
\mf{H}^{2}({\Omega}).
%\mathbf{d}_0|_{\partial \Omega}\in \mf{H}^{\frac{3}{2}}({\partial\Omega})
\end{eqnarray}
 Then there exists a
unique triple $(\rho_\varepsilon,\mf{v}_\varepsilon,\mf{d}_\varepsilon)$ with the following
properties:
\begin{enumerate}
  \item[\quad \quad(1)] Regularity.
\begin{equation*}
\left\{\begin{aligned} &\rho_\varepsilon\geq 0\mbox{ a.e. in }Q_T,\quad
 {\rho}_\varepsilon\in C^0(\bar{I},L^\beta_{\mm{weak}}(\Omega))
\cap C^0(\bar{I},L^p(\Omega)),\ 1\leq p<\beta,\ \\
& {\rho}_\varepsilon^{\frac{\beta}{2}}\in L^2({I},H^1(\Omega)),
\quad \partial_t\rho_\varepsilon\in L^{\frac{5\beta-3}{4\beta}}(Q_T),\quad
\nabla^2\rho_\varepsilon\in \mf{L}^\frac{5\beta-3}{4\beta}(\Omega),\\
 &\nabla \rho_\varepsilon,\quad \rho_\varepsilon\mf{v}_\varepsilon\in
 L^2({I},{E}_0^{\frac{10\beta-6}{3\beta+3},\frac{5\beta-3}{4\beta}}(\Omega)),\quad
 \int_\Omega\rho_\varepsilon\mm{d}\mf{x}=\int_\Omega \rho_0\mm{d}\mf{x}, \\
& \mf{m}_\varepsilon:={\rho}_\varepsilon\mf{v}_\varepsilon \in
C^0(\bar{I},L^{\frac{2\beta}{\beta+1}}_{\mm{weak}}(\Omega)),\quad
\mathbf{d}_\varepsilon\in C^0(I,\mathbf{H}^1(\Omega)),\quad
\mf{d}|_{\partial \Omega}=\mf{d}_0\in
C^{0,1}(\partial\Omega).\end{aligned} \right.\end{equation*}
  \item[\quad \quad(2)] Weak-strong solutions.
\begin{eqnarray}
&& \label{n0507f}\partial_t\rho_\varepsilon+\mm{div}(\rho_\varepsilon\mf{v}_\varepsilon )
-\varepsilon\Delta \rho_\varepsilon =0\mbox{ in }\mathcal{D}'(Q_T),\\[1mm]
&&\label{n0508f}\begin{aligned}
&\partial_t(\rho_\varepsilon\mf{v}_\varepsilon)+
\partial_j(\rho_\varepsilon\mf{v}_\varepsilon v_\varepsilon^j)-\mu
\Delta \mf{v}_\varepsilon-(\mu+\lambda)\nabla \mm{div}\mf{v}_\varepsilon+\nabla A\rho^\gamma_\varepsilon
+\delta\nabla \rho_\varepsilon^\beta \\
&\quad +\nu\mathrm{div}\left(\nabla
\mathbf{d}_\varepsilon\otimes\nabla\mathbf{d}_\varepsilon-\frac{|\nabla
\mf{d}_\varepsilon|^2\mathbb{I}}{2}\right)+ \varepsilon
\nabla \rho_\varepsilon\cdot \nabla \mf{v}_\varepsilon=0\;\mbox{ in }
(\mathcal{D}'(Q_T))^2,\end{aligned}\\[1mm]
&& \label{n0509f}
\partial_t\mathbf{d}_\varepsilon+\mathbf{v}_\varepsilon\cdot \nabla
\mathbf{d}_\varepsilon= \theta(\Delta \mathbf{d}_\varepsilon+|\nabla
\mathbf{d}_\varepsilon|^2\mf{d}_\varepsilon),\mbox{ a.e. in }Q_T.
\end{eqnarray}
Moreover,
$(\rho_\varepsilon,\mf{v}_\varepsilon,\mf{d}_\varepsilon)(\mf{x},0)=(\rho_0,\mf{v}_0,\mf{d}_0)$.
  \item[\quad \quad(3)]
  $(\rho_\varepsilon,\mf{v}_\varepsilon,\mf{d}_\varepsilon)$
   satisfies the finite and bounded energy inequalities as in (\ref{0457enegy}) and
(\ref{0458enegy}), and the uniform estimates (\ref{n047712}) and (\ref{n04771231}).
  \item[\quad \quad(4)] Additional uniform estimates.
\begin{eqnarray*}
 &&\label{n0515}\begin{aligned} &
\|\rho_\varepsilon\mf{v}_\varepsilon\|_{L^\infty(I,\mf{L}^{\frac{2\beta}
{\beta+1} }(\Omega))} + \|\rho_\varepsilon\mf{v}_\varepsilon\|_{L^2(I,\mf{L}^{\frac{ 6\beta}
{\beta+6}}(\Omega))}+\|\rho_\varepsilon\mf{v}_\varepsilon\|_{\mf{L}^{\frac{10\beta-6} {3\beta+3}}
(Q_T)} \\
&\quad +\|\rho_\varepsilon|\mf{v}_\varepsilon|^2\|_{L^2(I,{L}^{\frac{6\beta} {4\beta+3}
}(\Omega))}+\|\rho_\varepsilon\|_{L^{\beta+1}(Q_T)} +\varepsilon\|\nabla
\rho_\varepsilon\|_{L^{\frac{10\beta-6}{3\beta+3}}(Q_T)} \\
& +\sqrt{\varepsilon}\|\nabla \rho_\varepsilon\|_{\mf{L}^2(Q_T)}+
\varepsilon\|\nabla \rho_\varepsilon \cdot\nabla
\mf{v}_\varepsilon\|_{\mathbf{L}^{\frac{5\beta-3}{4\beta}}(Q_T)} \leq G(\delta),
\end{aligned}
\end{eqnarray*}
\end{enumerate}
where $G$ is a positive constant independent of $\varepsilon$.
\end{pro}

Here we explain how the assumptions (\ref{n0468}) in Proposition
\ref{thm:0301} become (\ref{n0505}) in Proposition \ref{thm:0302}.
First, we have Proposition \ref{thm:0302} for fixed initial data
satisfying (\ref{n0468}). If they satisfy the assumption
(\ref{n0505}) only, then we can approximate $(\rho_0,\mf{m}_0,
\mf{d}_0)$ by a sequence $(\rho_{0m},\mf{m}_{0m},\mf{d}_{0m})$ which
satisfies (\ref{n0468}), (\ref{n0467}), and
\begin{eqnarray}\label{intialappr}
{\mathcal{E}_{\delta}}(\rho_{0m},\mf{m}_{0m},\mf{d}_{0m})\rightarrow
{\mathcal{E}_{\delta}}(\rho_0,\mf{m}_0,\mf{d}_0)\quad\mbox{as
}m\rightarrow +\infty;
\end{eqnarray}
moreover,
$\mathbf{d}_{0m}\rightarrow \mathbf{d}_0\mbox{ strongly in }\mf{H}^{2}(\Omega)\cap \mf{C}^{0,1}(\bar{\Omega})$.
Then it is sufficient to take $\mf{v}_{0m}=P_n (\mf{m}_0/\rho_0)\in
\mf{X}_n$. Obviously, Proposition \ref{thm:0302} holds with these
new initial data, and those uniform estimates in  Proposition
\ref{thm:0302} are independent of $m$. Therefore, repeating the
limit process as $m\rightarrow\infty$, we obtain Proposition \ref{thm:0302} as well.

Now, with Proposition \ref{thm:0302} in hand, following the proof
of \cite[Proposition 5.1]{WDHYCGA} or \cite[Proposition 7.27]{NASII04}, we can obtain a weak solution
$(\rho_\delta,\mf{v}_\delta,\mf{d}_\delta)$ of the first approximate
problem as the weak limit of the sequence
$(\rho_\varepsilon,\mf{v}_\varepsilon,\mf{d}_\varepsilon)$ as $\varepsilon\rightarrow0$
constructed in Proposition \ref{thm:0302}. Thus, we have the following existence result.
%%%%%%%%%%%%%%%%%%%%%%%%%%%%%%%%%%%%%%%%%%%%%%%%%%%%
\begin{pro}\label{thm:0303}
 Let $\delta>0$,
$\beta\geq\max\{\gamma,8\}$. Assume that $\Omega$ is a bounded
$C^{2,\alpha}$-domain ($\alpha\in(0,1)$), and
$(\rho_0,\mf{m}_0,\mf{d}_0)$ satisfies
(\ref{jfw0110})--(\ref{jfw0112}) with $\beta$ in place of $\gamma$,
and (\ref{n0467}).
%${\mathcal{E}_{\delta}}(\rho_0,\mf{m}_0,\mf{d}_0)< \nu/4094$.
 Then there exists a unique triple $(\rho_\delta,\mf{v}_\delta, \mf{d}_\delta)$ with the
following properties:
\begin{enumerate}
  \item[\quad \quad(1)] Regularity.
\begin{equation*}
\left\{\begin{aligned} & {\rho}_\delta\in C^0(\bar{I},L^\beta_{\mm{weak}}(\Omega))
\cap C^0(\bar{I},L^p(\Omega))\cap L^{\beta+1}(\mathbb{R}^3\times I), \, 1\leq p<\beta, \\
&\rho_\delta\geq 0\;\mbox{ a.e. in }Q_T,\quad \rho_\delta=0,
\quad \mf{v}_\delta=\mf{0}\mbox{ in }(\mathbb{R}^2\backslash\Omega)\times I, \\
 &\rho_\delta|\mf{v}_\delta|^2\in  L^2(I,L^{\frac{6\beta}{4\beta+3}}
(\mathbb{R}^2))\cap L^1(I,L^{\frac{3\beta}{\beta+3}}(\mathbb{R}^2)),\\
&\mf{m}_\delta:={\rho}_\delta\mf{v}_\delta \in
L^2(I,\mf{L}^{\frac{6\beta}{6+\beta}}(\mathbb{R}^2))\cap
C^0(\bar{I},\mathbf{L}^{\frac{2\beta}{\beta+1}}_{\mm{weak}}(\Omega)),\\
&\mathbf{d}_\delta\in C^0(I,\mathbf{H}^1(\Omega)),\quad
\mf{d}_\delta|_{\partial \Omega}=\mf{d}_0 .\end{aligned}
\right.\end{equation*}
  \item[\quad \quad(2)]
  $(\rho_\delta,\mf{v}_\delta,\mathbf{d}_\delta)$ solves
  (\ref{n0507f})--(\ref{n0509f}) with
  $\varepsilon=0$, and also satisfies
  the finite and bounded energy inequalities as in  (\ref{0457enegy}) and
(\ref{0458enegy}) with $\varepsilon=0$. Moreover,
$(\rho_\delta,\mf{v}_\delta,\mf{d}_\delta)|_{t=0}=(\rho_0,\mf{v}_0,\mf{d}_0)$.
%\begin{eqnarray*}
%&&\partial_t\rho_\delta +\mm{div}(\rho_\delta\mf{v}_\delta
%)=0\mbox{ in }\mathcal{D}'(\mathbb{R}^2\times I),\\[1mm]
%&&\begin{aligned} &\partial_t(\rho_\delta\mf{v}_\delta)+
%\partial_j(\rho_\delta\mf{v}_\delta v_\delta^j)-\mu
%\Delta \mf{v}_\delta-(\mu+\lambda)\nabla
%\mm{div}\mf{v}_\delta+\nabla
%A\rho^\gamma_\delta\\
%&\quad +\nu\mathrm{div}\left(\nabla \mathbf{d}_\delta \otimes\nabla
%\mathbf{d}_\delta-\frac{|\nabla
%\mf{d}_\delta|^2\mathbb{I}}{2}\right)+\delta\nabla \rho_\delta^\beta
%]=0\ \mbox{ in }
%(\mathcal{D}'(Q_T))^2,\end{aligned}\\[1mm]
%&& \partial_t\mathbf{d}_\delta+\mathbf{v}_\delta\cdot \nabla
%\mathbf{d}_\delta= \theta(\Delta \mathbf{d}_\delta+|\nabla
%\mathbf{d}_\delta|^2\mf{d}_\delta),\mbox{ a.e. in
%}Q_T.\end{eqnarray*}
\item[\quad \quad(3)] For any $b$ satisfying (\ref{0110}) and (\ref{0111}),
 the function $b(\rho_\delta)$ is in
 $C^0(\bar{I},L^{\frac{\beta}{\lambda_1+1}}(\Omega)) \cap C^0(\bar{I},
L^p(\Omega))$, $1\leq p<\frac{\beta}{\lambda_1+1}$. Moreover,
\begin{equation*}
\partial_tb(\rho_\delta)+\mathrm{div}[b(\rho_\delta)\mathbf{v}]+\left[\rho_\delta b'(\rho_\delta)-
b(\rho_\delta)\right]\mathrm{div}\mathbf{v}=0\;\mbox{ in }
\mathcal{D}'(\mathbb{R}^2\times I).\end{equation*} For any $b_k$ ($k>0$) defined by
\begin{equation*}b_k(s)=\left\{\begin{array}{ll}
                  b(s) &\mbox{ if }s\in [0,k),\\
                  b(k)&\mbox{ if }s\in [k,\infty),
                \end{array}\right.\quad (b_k)'_+(s)=\left\{\begin{array}{ll}
                  b'(s) &\mbox{ if }s\in [0,k),\\
                  0&\mbox{ if }s\in [k,\infty),
                \end{array}\right.
\end{equation*}
with $b$ satisfying (\ref{0110}) and (\ref{0111}), $b_k(\rho_\delta)
$ belongs to $C^0(\bar{I}, L^p(\Omega))$, $1\leq p<\infty$. Moreover,
\begin{equation*}
\partial_tb_k(\rho_\delta)+\mathrm{div}[b_k(\rho_\delta)\mathbf{v}]+
\left[\rho_\delta (b_k)'_+(\rho_\delta)-b_k(\rho_\delta)\right]\mathrm{div}\mathbf{v}
=0\;\mbox{ in }\mathcal{D}'(\mathbb{R}^2\times I).
\end{equation*}
 % \item[\quad \quad(4)] Finite and bounded energy inequalities.
%\begin{equation*}\begin{aligned}
%\frac{d}{dt}{\mathcal{E}}_\delta(t)+\int_\Omega (\nu|\nabla
%\mf{v}_\delta|^2+(\mu+\lambda)
%|\mm{div}\mf{v}_\delta|^2)\mm{d}\mf{x}\leq 0\ \mbox{ in
%}\mathcal{D}'(I),
%\end{aligned}\end{equation*} and
%\begin{equation*}\begin{aligned}
%{\mathcal{E}}_\delta(t)+ \int_0^t\int_\Omega (\nu|\nabla
%\mf{v}_\delta|^2+(\mu+\lambda)
%|\mm{div}\mf{v}_\delta|^2\mm{d}\mf{x}\mm{d}s\leq
%{\mathcal{E}_{\delta}}(\rho_0,\mf{m}_0),\quad \mbox{a.e. in }I,
%\end{aligned}\end{equation*} where ${\mathcal{E}}_\delta(t):=
%{\mathcal{E}}_\delta(\rho_\delta,\mf{m}_\delta ,\mf{d}_\delta)$.
  \item[\quad \quad(4)]
  $(\rho_\delta,\mf{v}_\delta,\mathbf{d}_\delta)$ satisfies  the uniform estimates (\ref{n047712}) and
(\ref{n04771231}).  Moreover,
\begin{eqnarray*}
&&\begin{aligned}&
\|\rho_\delta\mf{v}_\delta\|_{L^\infty(I,\mf{L}^{\frac{2\gamma}
{\gamma+1} }(\Omega))}+
\|\rho_\delta\mf{v}_\delta\|_{L^2(I,\mf{L}^{\frac{ 6\gamma}
{\gamma+6}
}(\Omega))}+\|\rho_\delta|\mf{v}_\delta|^2\|_{L^2(I,{L}^{\frac{
3\gamma} {\gamma+3}
}(\Omega))}\\
&\quad +\|\rho_\delta|\mf{v}_\delta|^2\|_{L^2(I,{L}^{\frac{ 6\gamma}
{4\gamma+3} }(\Omega))}
+\|\rho_\delta\|_{L^{\gamma+\theta}(Q_T)}+\delta^{\frac{1}{\beta+\theta}}
\| \rho_\delta \|_{L^{{\beta+\theta}}(Q_T)}\leq
G(\theta)\end{aligned}
\end{eqnarray*}
\end{enumerate}
for any constant $\theta\in (0,\gamma-1)$ (the uniform estimate on
$\|\rho_\delta\|_{L^{\gamma+\theta}(Q_T)}$ can be found in
\cite{JFTZWHAJ3132} for the two-dimensional case). Here $G$ is a
positive constant independent of $\delta$.
\end{pro}
\begin{rem} We assume that the initial density $\rho_0$ satisfies (\ref{0920n}) in Proposition \ref{thm:0302}.
However the condition (\ref{0920n}) on $\rho_0$ can be relaxed by
the condition ``$\rho_0\in L^\beta(\Omega)$ and $\rho_0\geq 0$" in
Proposition \ref{thm:0303}. This process can be dealt by the
approximation method (similarly as in \eqref{intialappr}), please
refer to \cite[Section 7.10.7]{NASII04} for the detailed proof.
\end{rem}
Finally, with the help of Proposition \ref{thm:0303}, we can follow the
proof of \cite[Theorem 2.1]{WDHYCGA} or \cite[Theorem 7.7]{NASII04} to
obtain a weak solution $(\rho,\mf{v},\mf{d})$ of the
original problem (\ref{0101})--(\ref{0105}) which is the weak limit as $\delta\to 0$
of the weak solution sequence $(\rho_\delta,\mf{v}_\delta,\mf{d}_\delta)$ constructed
in Proposition \ref{thm:0303}. This completes the proof of Theorem \ref{thm:0101}.

\section{Global existence of large solutions to  the Cauchy problem}\label{Cauchy}
%\subsection{Energy estimates to Cauchy problem}
%\subsection{Construction of approximate solutions}
In this section, we will briefly describe how to prove Theorem
\ref{thm:0102} on the Cauchy problem by modifying the proof
of Theorem \ref{thm:0101} and applying the domain expansion technique.

\subsection{Local solvability of the Cauchy problem on the direction vector}
First we establish the local existence of solutions to the
following  Cauchy problem on the direction vector:
%(i.e. the second approximate problem to the problem (\ref{0101})--(\ref{0105}))
\begin{eqnarray}
\label{0601}\partial_t\mathbf{d}+\mathbf{v}\cdot \nabla \mathbf{d}=
\theta(\Delta \mathbf{d}+|\nabla
\mathbf{d}|^2(\mf{d}+\mf{e}_2))\qquad \mbox{ in }\mathbb{R}^2\times I
\end{eqnarray}
 with initial data
\begin{eqnarray}\label{0602}
\mathbf{d}(\mf{x},0):=\mathbf{d}_0\in \mf{H}^2(\mathbb{R}^2),\;\
|\mf{d}_0(\mf{x})+\mf{e}_2|=1,\;\  {d}_{02}(\mf{x})+1\geq
\underline{{d}}_{02}>0\end{eqnarray} for some given constant
$\underline{{d}}_{02}$, where $d_{02}$
 denotes the second component of $\mf{d}_0$, and the known function
$$\mf{v}\in \tilde{\mathbb{V}}:=\{\mf{v}\in C^0
(\bar{I},\mathbf{H}^{1}(\mathbb{R}^2)\cap
\mf{L}^\infty(\mathbb{R}^2))~|~ \partial_t\mf{v}\in
L^{2}(I,\mathbf{H}^{1}(\mathbb{R}^2))\}.$$ % and
%$(\mf{v},\partial_t\mf{v})\equiv \mf{0}$ in $(\mathbb{R}^2\backslash
%B_R)\times I$ for some given positive integer $R\geq 0$.
The local existence can be
shown by following the proof of Proposition \ref{pro:0402} and applying the
domain extension technique. Next, we briefly describe the proof.

By virtue of (\ref{0602}), we have a sequence of functions
$\{{\mathbf{d}}_n^0(\mf{x})\}_{n\geq 1}\subset \mf{H} ^2(\mathbb{R}^2)$, such that
$$|{\mathbf{d}}_n^0(\mf{x})+\mf{e}_2|=1,\quad
{d}_{n2}^0(\mf{x})+1\geq  {d}_{02}/2, \quad
{\mathbf{d}}_n^0\rightarrow \mathbf{d}_0 \text{ in }
\mf{H}^2(\mathbb{R}^2)$$ and $\mm{supp}\,\mathbf{d}_n^0\subset
B_{R_{n}}:=\{\mf{x}\in \mathbb{R}^2~|~|\mf{x}|<R_n\}$ for some
$R_n\geq 1$.
%Let $\{S_{R_n^{-1}}\}_{n\geq 1}$ is a family of
%2-dimensional mollifiers. Then
%\begin{equation*}\begin{aligned}&
%S_{R_n^{-1}}(\mf{v})\in C^0
%(\bar{I},\mathbf{H}^{1}(\mathbb{R}^2)\cap
%\mf{L}^\infty(\mathbb{R}^2))\cap L^p(I,
%\mf{C}_0^\infty(\mathbb{R}^2)),\ \mbox{ for any }p>1,\\
%& \mm{supp}S_{R_n^{-1}}(\mf{v},\partial_t\mf{v})\subset
%B_{R_n}\mbox{ for a.e. }t\in I,\
% \|S_{R_n^{-1}}(\mf{v})\|_{\mathbb{V}}\leq   \|\mf{v}\|_{\mathbb{V}},
%\\ &(\partial_t S_{R_n^{-1}}(\mf{v})=S_{R_n^{-1}}(\partial_t\mf{v})
%,S_{R_n^{-1}}(\mf{v}))\rightarrow (\partial_t\mf{v},\mf{v}) \mbox{
%in
%}L^2(\bar{I},\mathbf{H}^{1}(\mathbb{R}^2)),\end{aligned}\end{equation*}
%where the norm $\|\cdot\|_{\mathbb{V}}$ is defined by (\ref{0415})
%with $\mathbb{R}^2$ in place of $\Omega$.
 In view of Proposition \ref{pro:0402}, there exists a unique local strong solution $\mf{d}$
to the following Direchlet problem:
\begin{eqnarray}\label{0603}
  \partial_t\mathbf{d}+%S_{R_{n}^{-1}}(
\mathbf{v}\cdot \nabla
\mathbf{d}= \theta(\Delta \mathbf{d}+|\nabla
\mathbf{d}|^2\mf{d})\qquad \mbox{ in }B_{R_n}\times I,
\end{eqnarray}
 with initial and boundary conditions
\begin{eqnarray}\label{0604}
\mathbf{d}(\mf{x},0)={\mathbf{d}}_n^0\in \mf{H}^2(B_{R_n}),\quad
\mathbf{d}(\mf{x},t)|_{\partial B_{R_n}}=\mf{0}\;\mbox{ for any } t\in
I.\end{eqnarray}
Moreover $d_2\geq c(\underline{d}_{02})$. However, the constants $h_1$ and
$C$ in Proposition \ref{pro:0402} depend on the bounded
domain $\Omega=B_{R_n}$, since we have used the interpolation
inequalities (\ref{jjjw0432}) and (\ref{njjjw0432}) (of course, the
embedding theorems $H^2(\Omega)\hookrightarrow L^\infty(\Omega)$ and
$H^1(\Omega)\hookrightarrow L^4(\Omega)$ used in (\ref{nsfs0333jwd})
and (\ref{jjwfsaf123}) can be replaced by such two interpolation
inequalities, respectively),  where the interpolation constants
depend on the bounded domain $\Omega=B_{R_n}$. Fortunately, we have
the following results:
\begin{equation}\label{infty0604}\|v\|_{L^\infty(B_{R_n})}^2\leq
{c}\| v\|_{H^{2}(B_{R_n})}\|v\|_{L^2(B_{R_n})}\mbox{ and
}\|v\|_{L^4(B_{R_n})}^4\leq {c}\|
v\|_{H^{1}(B_{R_n})}^2\|v\|_{L^2(B_{R_n})}^2
\end{equation}  for some
constant $c$ independent of $R_n\geq 1$, which can be deduced from
(\ref{jjjw0432}) and (\ref{njjjw0432}) with $B_1$ in place of
$\Omega$ by scaling the spatial variables. In addition,  the elliptic estimate
 used in (\ref{nestf0429}) is replaced by the
following special elliptic estimate:
if $v\in H^1_0(B_{R_n})$, $f\in H^1(B_{R_n})$ and $\Delta v=f$ in
$\mathcal{D}'(B_{R_n})$, then
\begin{equation}\label{n0605}
\|\nabla^3 v\|_{L^2(B_{R_n})}^2\leq {c}\| f\|_{H^{1}(B_{R_n})}\;\mbox{ for some constant } c
\mbox{ independent of }R_n\geq 1,
\end{equation}
  which can be deduced from the standard elliptic estimate on the domain
$B_1$  by scaling the spatial variables. With these facts, repeating the proof of
Proposition \ref{pro:0402} with a trivial modification, we can obtain
Proposition \ref{pro:0402} with $\Omega=B_{R_n}$. Moreover, the
constants $h_1$ and $C$ can be independent of $B_{R_n}$ and
$\|\nabla^3 \mf{d}_0\|_{\mf{L}^2(B_{R_n})}$. Then, using the domain expansion
technique (see \cite{CYKHOM}), we can obtain a unique strong solution
$\mf{d}$ of the Cauchy problem (\ref{0601}), (\ref{0602}) as the
weak limit of the sequence $\mf{d}_{R_n}$ with initial data
$\mf{d}_{R_n}|_{t=0}=\mathbf{d}_n^0$ as $R_n\rightarrow\infty$.
Thus, we have
%
% Moreover, the result of local
% well-posedness to the Cauchy problem can be read as follows.
%
\begin{pro}\label{pro:0602}
Proposition \ref{pro:0402} holds with $\mathbb{R}^2$ in place of $\Omega$, where
 $\mf{d}=\mathscr{D}^K_{\mathbf{d}_0}(\mathbf{v})$ solves the
 Cauchy problem (\ref{0601}), (\ref{0602}), the constants $h_1$
 and $C$ only depend on $K$ and $\|\mf{d}_0\|_{\mathbf{H}^2(\mathbb{R}^2)}$, and
 \begin{eqnarray}
\label{0617n1}d_2(\mf{x},t)+1\geq \underline{d}_{02}/2\;\mbox{ for
any }(\mf{x},t)\in \mathbb{R}^2\times I_d^K, \;\mbox{ if
}{d}_{02}(\mf{x})+1\geq \underline{d}_{02}.
\end{eqnarray}
\end{pro}

\subsection{Global solvability of the approximate problem}

Once we have established Proposition \ref{pro:0602}, we can see,
in view of the proof in Section \ref{sec:05}, that
  there exists a unique local solution to the
following approximate problem to the original problem (\ref{0101})--(\ref{0104})
with $\Omega=\mathbb{R}^2$:
\begin{eqnarray}&&\begin{aligned}\label{0610}&\int_{B_R}
(\rho\mathbf{v})(t)\cdot\mathbf{\mathbf{\Psi}}\mathrm{d}\mathbf{x}
-\int_{B_R}\mathbf{m}_0\cdot\mathbf{\mathbf{\Psi}}\mathrm{d}\mathbf{x} \\
& = \int_0^t\int_{B_R}\bigg[\mu\Delta\mathbf{v}+(\mu+\lambda)\nabla
\mm{div}\mf{v}-A\nabla \rho^\gamma -\delta \nabla
\rho^\beta-\varepsilon(\nabla\rho\cdot\nabla \mathbf{v})
-\mathrm{div}(\rho\mathbf{v}\otimes\mathbf{v})  \\
&\quad -\nu\mathrm{div}\left(\nabla \mathbf{d}\otimes\nabla \mathbf{d}-\frac{|\nabla
\mf{d}|^2\mathbb{I}}{2}\right)\bigg]\cdot\mathbf{\Psi}\mathrm{d}\mathbf{x}\mathrm{d}s
\quad\mbox{ for all }t\in I\mbox{ and any }\mathbf{\Psi}\in\mathbf{X}_n,
\end{aligned}\\
&& \label{0608}\partial_t\rho+\mathrm{div}(\rho\mathbf{v})=\varepsilon \Delta
\rho\qquad \mbox{ in }B_R\times I,\\
&&\label{0609}\partial_t\mathbf{d}+\mathbf{v}\cdot \nabla
\mathbf{d}= \theta(\Delta \mathbf{d}+|\nabla
\mathbf{d}|^2\mf{d})\qquad \mbox{ in }\mathbb{R}^2\times I,
\end{eqnarray}
 with boundary conditions
\begin{eqnarray} \nabla\rho\cdot\mf{n}|_{\partial B_R}=0,
\ \mathbf{v}|_{\partial B_R}=\mathbf{0},\
\mathbf{d}(\mf{x},t)-\mf{e}_2\to\mf{0}\mbox{ as }|\mf{x}|\rightarrow +\infty
\end{eqnarray}
and modified initial conditions
\begin{eqnarray}\label{0612}&&\rho(\mathbf{x},0)=
\rho_0\in W^{1,\infty}( B_R),\quad 0<\underline{\rho}\leq
\rho_0\leq \bar{\rho}<\infty, \\
&&\label{0613}
 \mathbf{v}(\mathbf{x},0)=\mf{v}_0\in \mathbf{X}_n,\quad
 \mathbf{d}(\mf{x},0)-\mf{e}_2\in \mf{H}^3(\mathbb{R}^2),
 \end{eqnarray}
where $(\mf{v},\partial_t\mf{v})\equiv \mf{0}$ in
$(\mathbb{R}^2\backslash B_R)\times I$ in (\ref{0609}). Here we have
used the fact that $\mf{v}\in L^2(I,\mf{H}^1_0(B_R))$ is equivalent
to $\mf{v}\in L^2(I,\mf{H}^1(\mathbb{R}^2))$ and $\mf{v}=\mf{0}$ in
$(\mathbb{R}^2\backslash B_R)\times I$. Hence one has $\mf{v}\in
C^0(I,\mathbf{H}^{1}(\mathbb{R}^2)\cap\mf{L}^\infty(\mathbb{R}^2))$
and $\partial_t\mf{v}\in L^{2}(I,\mathbf{H}^{1}(\mathbb{R}^2))$ by zero extension,
if $\mf{v}\in C(\bar{I}, \mathbf{X}_n)$ and $\partial_t\mf{v}\in L^2(I,\mf{X}_n)$
with $\Omega=B_R$.

Similar to the arguments in Section \ref{sec:0404}, we shall derive the
uniform-in-time energy estimates in order to extend the local solution globally in time.
This can be done by modifying the derivation of the uniform estimates in Section
\ref{sec:0404}. The only difference in arguments arises from the energy estimates
for $\mf{d}$ in the whole space, recalling that we have shown the uniform estimates on
$\|\Delta\mf{d}\|_{L^2(I,\mathbf{L}^2(\Omega))}$ and
$\|\Delta\mf{d}\|_{L^\infty(I,\mathbf{L}^2(\Omega))}$ in (\ref{n0452}) and (\ref{fzu0447})
in a bounded domain. In the case of $\Omega=\mathbb{R}^2$, however,
if ${d}_2$ satisfies some geometric angle condition,
then we can use the rigidity theorem to deduce the uniform estimates on
$\|\Delta\mf{d}\|_{L^2(I,\mathbf{L}^2(\mathbb{R}^2))}$ and
$\|\Delta\mf{d}\|_{L^\infty(I,\mathbf{L}^2(\mathbb{R}^2))}$. The
rigidity theorem, which was recently established in \cite{LZLDZXY}, reads as follows.
%%%%%%%%%%%%%%%%%%%%%%%%%%%%%%%%%%%%%%%%%%%%%%%%%%%
\begin{pro}\label{pro:0603} Let $\underline{d}_{2}>0$, $c_0>0$. Then there exists a
positive constant $\varpi_0=\varpi(\underline{d}_{2},c_0)$, such that the following holds.

If $\mf{d}:=(d_1,d_2):\mathbb{R}^2\rightarrow \mathbb{S}^1$, $\nabla
\mf{d}\in \mf{H}^1(\mathbb{R}^2)$ with $\|\nabla
\mf{d}\|_{\mf{L}^2(\mathbb{R}^2)}\leq c_0$ and $d_2\geq
\underline{d}_{2}$, then
\begin{equation*}  \|\nabla \mf{d}\|_{\mf{L}^4(\mathbb{R}^2)}^4\leq
(1-\varpi_0)\|\Delta\mf{d}\|_{\mf{L}^2(\mathbb{R}^2)}^2.
\end{equation*}

 Consequently, for such a map the associated harmonic energy is coercive, i.e.,
\begin{equation}\label{0614}\|\Delta\mf{d}+|\nabla
\mf{d}|^2\mf{d}\|_{\mathbf{L}^2(\mathbb{R}^2)}^2\geq
\frac{\varpi_0}{2}(\|\Delta\mf{d}\|_{\mathbf{L}^2(\mathbb{R}^2)}^2+\|\nabla
\mf{d}\|_{\mathbf{L}^4(\mathbb{R}^2)}^4).
\end{equation}
\end{pro}

On the other hand, we have used the fact $\mf{d}\equiv 1$ and the
boundedness of domain $\Omega$ to get the uniform estimate
$\|\mf{d}\|_{L^\infty(I,\mf{L}^2(\Omega))}$ in (\ref{fzu0447}). For
the case of the whole space, we need the uniform estimate
$\|\mf{d}-\mf{e}_2\|_{L^\infty(I,\mf{L}^2(\mathbb{R}^2))}$, which
can be deduced from $\mf{d}\equiv 1$  and (\ref{0601}). More
precisely, we have
\begin{equation*}\begin{aligned}\label{ffww1314b}
&\frac{1}{2}\frac{d}{dt}\|\mf{d}-\mf{e}_2\|^2_{\mf{L}^2(\mathbb{R}^2)}+\theta\|\nabla
\mf{d}\|^2_{\mf{L}^2(\mathbb{R}^2)}
%&=\theta\int_0^t\int_{\mathbb{R}^2}[|\nabla
%\mf{d}|^2(1-d_2)-\mf{v}\cdot \nabla \mf{d}\cdot
% ( \mf{d}-\mf{e}_2)]\mm{d}\mf{x}\mm{d}s+\frac{1}{2}\|\mf{d}_0-\mf{e}_2\|^2_{\mf{L}^2(\mathbb{R}^2)}\\
=\theta\int_{\mathbb{R}^2}\left[|\nabla \mf{d}|^2(1-d_2)+\frac{1}{2}|
 \mf{d}-\mf{e}_2|^2\mm{div}\mf{v}\right]\mm{d}\mf{x} \\
 %+\frac{1}{2}\|\mf{d}_0-\mf{e}_2\|^2_{\mf{L}^2(\mathbb{R}^2)}
&\leq \left(2\theta+\frac{\theta^2}{4} \|
\mm{div}\mf{v}\|_{{L}^2(\mathbb{R}^2)}^2\right)\|\nabla
\mf{d}\|_{\mf{L}^2(\mathbb{R}^2)}^2+\frac{1}{2}\|\mf{d}-\mf{e}_2\|_{\mf{L}^2(\mathbb{R}^2)}^2,
%+\frac{1}{2}\|\mf{d}_0-\mf{e}_2\|^2_{\mf{L}^2(\mathbb{R}^2)}.
   \end{aligned}\end{equation*}
which, together with  Gronwall's inequality, yields
 \begin{equation}\begin{aligned}\label{ffww1314}&\|\mf{d}-\mf{e}_2\|^2_{\mf{L}^2(\mathbb{R}^2)}
\leq \|\mf{d}_0-\mf{e}_2\|^2_{\mf{L}^2(\mathbb{R}^2)}+\left(4\theta T+\frac{\theta^2}{2} \|
\mm{div}\mf{v}\|_{L^2(I,{L}^2(\mathbb{R}^2))}^2\right)e^T\|\nabla
\mf{d}\|_{L^\infty(I,\mf{L}^2(\mathbb{R}^2))}^2.
   \end{aligned}\end{equation}

%Thirdly  we have used the  elliptic estimates on bounded domain to
%deduce the uniform estimates  $\|\nabla^2\mf{d}
%\|_{L^\infty(I,\mf{H}^1(\Omega))}$  in (\ref{fzn0448}) and
%(\ref{fzu0447}).  For the case of whole space,  we shall use another
%version of elliptic estimates: let $k\geq 0$ be a integer, $v\in
%H^1_0(\mathbb{R}^2)$, $f\in H^k(\mathbb{R}^2)$ and $\Delta v=f$ in
%$\mathcal{D}'(\mathbb{R}^2)$, then
%\begin{equation}\label{n0605}\|\nabla ^{k+2} v\|_{L^2(\mathbb{R}^2)}^2\leq
%{c}\| f\|_{H^{k}(\mathbb{R}^2)}\mbox{ for some constant } c ,
%\end{equation}
% please refer to \cite[Lemma 15]{CHJCHJKHU}.

%Finally, we have used $H^1(\Omega)\hookrightarrow L^6(\Omega)$ to
%get $\|\nabla \mf{d}\|_{\mf{L}^6(\Omega)}^6\leq c(\Omega)\|\nabla
%\mf{d}\|_{\mf{H}^1(\Omega)}^6$ in (\ref{fzn0448}).  For the case of
%whole space, using (\ref{infty0604}), the term $\|\nabla
%\mf{d}\|_{\mf{L}^6(\mathbb{R}^2)}^6$ can be bounded as follows:
%\begin{equation}\begin{aligned}\label{fzun0620}&\|\nabla
%\mf{d}\|_{\mf{L}^6(\mathbb{R}^2)}^6\leq  \|\nabla
%\mf{d}\|_{\mf{L}^\infty(\mathbb{R}^2)}^2\|\nabla
%\mf{d}\|_{\mf{L}^4(\mathbb{R}^2)}^4\leq c \|\nabla
%\mf{d}\|_{\mf{H}^{2}(\mathbb{R}^2)}\|\nabla
%\mf{d}\|_{\mf{L}^2(\mathbb{R}^2)}\|\nabla
%\mf{d}\|_{\mf{L}^4(\mathbb{R}^2)}^4\\
%\leq &c (\|\nabla \mf{d}\|_{\mf{H}^{1}(\mathbb{R}^2)}\|\nabla
%\mf{d}\|_{\mf{L}^2(\mathbb{R}^2)}\|\nabla
%\mf{d}\|_{\mf{L}^4(\mathbb{R}^2)}^4+\|\nabla
%\mf{d}\|_{\mf{L}^2(\mathbb{R}^2)}^2\|\nabla
%\mf{d}\|_{\mf{L}^4(\mathbb{R}^2)}^8)+\frac{1}{4}\|\nabla^3
%\mf{d}\|_{\mf{L}^{2}(\mathbb{R}^2)}^2.
%\end{aligned}\end{equation}

Consequently, plugging (\ref{0614}), (\ref{ffww1314}) and
Proposition \ref{pro:0602} with $\underline{d}_{02}>0$ into the
deduction of the global estimates in Section \ref{sec:0404}, we immediately
establish the following global existence of a unique solution $(\rho_n,\mf{v}_n,\mathbf{d}_n)$
to the approximate problem (\ref{0610})--(\ref{0613}):
%%%%%%%%%%%%%%%%%%%%%%%%%%%%%%%%%%%%%%%
\begin{pro}\label{pro:0604}
Let $\Omega=B_R$ with $R>0$, $\underline{d}_{02}$ be a positive
constant, $( \beta,\varepsilon,\underline{\rho},\bar{\rho})$ satisfy
(\ref{condition0452}) and $\delta\in (0,1]$. Assume $\rho_0$
satisfies (\ref{0920n}), $\mf{v}_0\in\mf{X}_n$ and
\begin{eqnarray}\label{finallyok}
{d}_{02}\geq \underline{d}_{02},\ |\mf{d}_0(\mf{x})|\equiv1,\
\mf{d}_0(\mf{x})-\mf{e}_2\in \mf{H}^2(\mathbb{R}^2).
\end{eqnarray}
 Then there exists a unique triple $(\rho_n,\mf{v}_n,\mf{d}_n)$ defined on
$\mathbb{R}^2\times I$ with the following properties:
\begin{enumerate}
  \item[\quad \quad(1)] $(\rho_n,\mf{v}_n)$ satisfies the regularity
  (\ref{ren0454}), and $\mf{d}_n$ satisfies
\begin{equation*}
\mf{d}_n-\mf{e}_2\in C^0(\bar{I},\mf{H}^2(\mathbb{R}^2)) \cap
L^2({I},\mf{H}^3(\mathbb{R}^2)),\
\partial_t\mf{d}_n\in C^0(\bar{I},\mathbf{L}^{2}(\mathbb{R}^2))\cap
{L}^{2}(I,\mf{H}^1(\mathbb{R}^2)).
\end{equation*}
  \item[\quad \quad(2)] $(\rho_n,\mf{v}_n,\mf{d}_n)$
solves (\ref{0608}) a.e. in $Q_T$ and (\ref{0609}) a.e. in $\mathbb{R}^2\times I$,
and satisfies (\ref{0610}).
   \item[\quad \quad(3)]   $(\rho_n,\mf{v}_n,\mf{d}_n)$ satisfies the finite
  and bounded energy inequalities as in Proposition \ref{thm:0301}
   with $\mathbb{R}^2$ in place of $\Omega$, where
$(\rho_n,\mf{v}_n,\nabla \rho_n)=\mf{0}$ in $(\mathbb{R}^2\backslash\Omega)\times I$.
Moreover, we have the uniform estimates (\ref{n0477}), (\ref{n0478}), and
\begin{eqnarray*} &&
|\mf{d}_n|\equiv 1 \mbox{ in }\mathbb{R}^2\times I, \\[1mm]
&&  \|\nabla^2 \mf{d}_n\|_{\mf{L}^2(\mathbb{R}^2\times
I)}+\|\nabla \mf{d}_n\|_{{L}^4(\mathbb{R}^2\times I)}\leq
G(\underline{d}_{02},\sup_{\delta\in
(0,1]}{\mathcal{E}_{\delta}}(\rho_0,\mf{q}_0,\mf{d}_0)). \\[1mm]
&&  \|\partial_t\mathbf{d}_n\|_{L^{4/3}(I,\mathbf{L}^2(\Omega))}
+\|\partial_t\mathbf{d}_n\|_{L^{2}(I,(\mf{H}^{1}(\Omega))^*)}\leq
G(\underline{d}_{02},\sup_{\delta\in
(0,1]}{\mathcal{E}_{\delta}}(\rho_0,\mf{q}_0,\mf{d}_0),\Omega).
\end{eqnarray*}
\end{enumerate}
% Here $G$ is a positive constant which is, in particular,
%independent of $n$, and nondecreasing in its argument. Moreover, if
%$\varepsilon$ is not explicitly written in the argument of $G$, then
%$G$ is independent of $\varepsilon$ as well.
\end{pro}
As a consequence of Proposition \ref{pro:0604}, similarly to Section \ref{sec:03}, we can use the
standard three-level approximation scheme and the method of weak
convergence based on Proposition \ref{pro:0604} to establish the
existence of weak solutions to the following problem:
\begin{eqnarray}
&& \label{0621}\partial_t\rho+\mathrm{div}(\rho\,\mathbf{v})=0\qquad \mbox{ in }B_R\times I,\\
 &&  \label{0622} \partial_t(\rho
\mathbf{v})+\mathrm{div}(\rho\mathbf{v}\otimes\mathbf{v})+ \nabla
P(\rho)\nonumber\\
&&\quad=\mu \Delta\mathbf{v}+(\mu+\lambda)\nabla
\mm{div}\mf{v}-\nu\mm{div}\left(\nabla \mathbf{d}\odot\nabla
\mathbf{d}-\frac{1}{2}|\nabla \mf{d}|^2\mathbb{I}\right)\quad \mbox{ in }B_R\times I,\\
&&  \label{0623}\partial_t\mathbf{d}+\mathbf{v}\cdot \nabla
\mathbf{d}= \theta(\Delta \mathbf{d}+|\nabla
\mathbf{d}|^2\mf{d})\qquad \mbox{ in }\mathbb{R}^2\times I.
\end{eqnarray}
More precisely, we have the following conclusion:
\begin{pro}\label{pro:0605}
Let the initial data $(\rho_0,\,\mathbf{m}_0)$ satisfy
(\ref{jfw0110}) and (\ref{fzu0109}) with $B_R$ in place of $\Omega$,
and $\mathbf{d}_0$ satisfy
\begin{equation}\label{thank}{d}_{02}\geq \underline{d}_{02},\;\
|\mf{d}_0|=1\;\mbox{ a.e. in }\mathbb{R}^2,\ \
   \mathbf{d}_0(\mathbf{x})-\mf{e}_2\in
   \mathbf{H}^1({\mathbb{R}^2}).\end{equation}
 Then the initial-boundary value problem (\ref{0621})--(\ref{0623}) has a
global weak solution $(\rho,\mathbf{v},\mathbf{d})$ defined on
$\mathbb{R}^2\times I$ for any given $T>0$, such that
\begin{enumerate}
            \item[\quad (1)]
$(\rho,\mf{v})$ satisfies the regularity (\ref{fz0114}) and
(\ref{fz0115}) with $\mathbb{R}^2$ in place of $\Omega$, and
$\mf{d}$ satisfies the regularity (\ref{fz0125}). Moreover,
$(\rho,\mf{v})=\mf{0}$ in $(\mathbb{R}^2\backslash B_R)\times I$.
    \item[\quad (2)]  Equations (\ref{0621}), (\ref{0622}) hold in
    $\mathcal{D}'(B_R\times I)$, the equation (\ref{0623}) holds a.e. in
$\mathbb{R}^2\times I$, the equation (\ref{0621}) is satisfied in
the sense of renormalized solutions, and the solution satisfies the
energy equality (\ref{eniq1}) with $\mathbb{R}^2$ in place of
$\Omega$.
 \item[\quad (3)] Additional estimate:
\begin{equation*}  \begin{aligned}
&\int_0^T\int_{\mathbb{R}^2}\frac{\theta\varpi_0}{2}(|\Delta \mf{d}|^2+|\nabla
\mf{d}|^4)\mathrm{d}\mathbf{x}\mm{d}s\leq \mathcal{E}_0,
\end{aligned}\end{equation*}
where the constant $\varpi_0$ depends on $\underline{d}_{02}$, $\nu$ and
$$\mathcal{E}_0=\int_{B_R}\left(\frac{1}{2}\frac{|\mf{m}_0|^2}{\rho_0}1_{\{\rho_0>0\}}+
Q(\rho_0)\right)\mm{d}\mf{x}+\int_{\mathbb{R}^2}\frac{\nu|\nabla
\mf{d}_0|^2}{2}\mm{d}\mf{x}.$$
\end{enumerate}
 \end{pro}
\begin{rem} To replace (\ref{finallyok}) by (\ref{thank}), we have used the fact that
there exists a sequence of approximate functions
    $\{\mf{d}_{m}^0\}_{m=1}^\infty\subset \mf{H}^1(\mathbb{R}^2)+\mf{e}_2$, such that
$$ |\mf{d}_{m}^0|=1,\;\ \mf{d}_{m}^0-\mf{d}_{0}\rightarrow
\mf{0}\mbox{ in }\mf{H}^1(\mathbb{R}^2), \mbox{ and } {d}_{m2}^0\geq
\underline{d}_{02}/2.$$
\end{rem}

Finally, we can follow the arguments in \cite[Section 7.11]{NASII04} for the
compressible Navier-Stokes equations to prove Theorem
\ref{thm:0102} by using Proposition \ref{pro:0605} on the bounded
invading domain $B_R$ and letting $R\rightarrow \infty$, where one should
use the embedding theorem
$$H^1(\Omega') \hookrightarrow L^p(\Omega')\;\mbox{ for any }p\geq 1\mbox{ and any
bounded Lipschiz domain } \Omega'\subset\mathbb{R}^2, $$
to replace $L^6(\mathbb{R}^3)\hookrightarrow D^{1,2}(\mathbb{R}^3)$ in
\cite[Section 7.11]{NASII04} in the treatment of the velocity $\mf{v}$.
Here we omit the proof, since the additional limit process on
$\mf{d}$ is trivial. We mention that, in the limit process on
$\mf{d}$, we also have the weak convergence  as in (\ref{302}), and the
strong convergence as in (\ref{strong1}) and (\ref{strong2}) with
any bounded space $\Omega'\subset \mathbb{R}^2$.
%%%%%%

\newcommand\ack{\section*{Acknowledgement}}
\newcommand\acks{\section*{Acknowledgements}}
\acks
 The authors are grateful to Prof. Xianpeng Hu for his suggestions
 which improved the presentation of this paper. The research of Fei Jiang
was supported by the NSFC (Grant No. 11101044),  the research of
Song Jiang by NSFC (Grant No. 11229101) and the National Basic
Research Program under the Grant 2011CB309705, and the research of
Dehua Wang by the National Science Foundation under Grant DMS-0906160
and the Office of Naval Research under Grant N00014-07-1-0668.

%The first author express his approbations to Prof. Song Jiang for
%his kind supports.
%The authors would like to thank the anonymous referee for invaluable
%suggestions.
%
\renewcommand\refname{References}
\renewenvironment{thebibliography}[1]{%
\section*{\refname}
\list{{\arabic{enumi}}}{\def\makelabel##1{\hss{##1}}\topsep=0mm
\parsep=0mm
\partopsep=0mm\itemsep=0mm
\labelsep=1ex\itemindent=0mm
\settowidth\labelwidth{\small[#1]}%
\leftmargin\labelwidth \advance\leftmargin\labelsep
\advance\leftmargin -\itemindent
\usecounter{enumi}}\small
\def\newblock{\ }
\sloppy\clubpenalty4000\widowpenalty4000
\sfcode`\.=1000\relax}{\endlist}
\bibliographystyle{model1b-num-names}

\end{document}